\documentclass[a4paper,12pt,authoryear]{elsarticle}

\usepackage[utf8]{inputenc}
\usepackage[T1]{fontenc}

\usepackage[english]{babel}
\usepackage{fouriernc}
\usepackage{amsthm,amssymb,amsfonts,amsmath,color,graphicx,float,url,bbm,multirow,caption,subfig,placeins,xspace,xr-hyper}

\usepackage[pagebackref = false]{hyperref}
\hypersetup{
  colorlinks = true,
  urlcolor = blue, 
  linkcolor = blue,
  citecolor = blue,
  pdftitle = {COBRA: A Combined Regression Strategy - \today},
  pdfauthor = {G. Biau, A. Fischer, B. Guedj and J. D. Malley},
  pdfsubject = {Regression, Combining estimators, COBRA}
}

\addto\extrasenglish{
  
}

\numberwithin{equation}{section}

\newtheorem{lem}{Lemma}[section]

\newtheorem{pro}{Proposition}[section]
\newtheorem{theo}{Theorem}[section]

\newtheorem{model}{\textsf{Model}}

\newcommand{\bX}{\mathbf{X}}
\newcommand{\bx}{\mathbf{x}}
\newcommand{\bz}{\mathbf{z}}

\newcommand{\br}{\mathbf{r}}
\newcommand{\by}{\mathbf{y}}
\newcommand{\bu}{\mathbf{u}}
\newcommand{\bv}{\mathbf{v}}
\newcommand{\e}{\varepsilon}

\def\e{\varepsilon}
\def\E{\mathbb{E}}
\def\1{\mathbf{1}}
\def\R{\mathbb{R}}
\def\var{\mathbb{V}}
\def\ie{\emph{i.e.}}

\def\cobra{\texttt{COBRA}\xspace}

\begin{document}

\begin{frontmatter}
  \title{COBRA: A Combined Regression Strategy}
  
  \author[a,b]{G\'erard Biau}
%  \ead{gerard.biau@upmc.fr}
  
  \author[c]{Aur\'elie Fischer}
%  \ead{aurelie.fischer@univ-paris-diderot.fr}
  
  \author[d]{Benjamin Guedj\corref{cor1}}
  \ead{benjamin.guedj@inria.fr}
  
  \author[e]{James D. Malley}
%  \ead{jmalley@mail.nih.gov}
  
  \address[a]{Universit\'e Pierre et Marie Curie, France}
  \address[b]{Institut universitaire de France}
  \address[c]{Universit\'e Paris Diderot, France}
  \address[d]{Inria, France}
  \address[e]{National Institutes of Health, USA}
  
  \cortext[cor1]{Corresponding author}

  \begin{abstract}\noindent
    A new method for combining several initial estimators of the
    regression function is introduced. Instead of building a linear or
    convex optimized combination over a collection of basic estimators
    $r_1,\dots,r_M$, we use them as a collective indicator of the
    proximity between the training data and a test observation. This
    local distance approach is model-free and very fast. More specifically, the resulting nonparametric/nonlinear combined estimator is shown to perform asymptotically at least as well in the $L^2$ sense as the best combination of the basic estimators in the collective. A companion R package called \cobra (standing for COmBined Regression Alternative) is
    presented (downloadable on
    \url{http://cran.r-project.org/web/packages/COBRA/index.html}). Substantial
    numerical evidence is provided on both
    synthetic and real data sets to assess the excellent performance and
    velocity of our method
    in a large variety of prediction problems.
    \medskip
    
    \noindent\emph{Index terms} --- Combining estimators, Consistency,
    Nonlinearity, Nonparametric regression, Prediction.
    
    \medskip
    
    \noindent\emph{2010 Mathematics Subject Classification}: 62G05, 62G20.
    
  \end{abstract}
  
\end{frontmatter}

\section{Introduction}

Recent years have witnessed a growing interest in combined
statistical procedures, supported by a considerable research and
extensive empirical evidence. Indeed, the increasing number of available
estimation and prediction methods (hereafter denoted \emph{machines})
in a wide range of modern statistical problems naturally suggests
using some efficient strategy for combining procedures and
estimators. Such an approach would be a valuable research and development tool, for example when dealing with high or infinite dimensional data.
\medskip

There exists an extensive literature on linear aggregation of
estimators, in a wide range of statistical models: A review of these
methods may be found for example in \citet{Gir2014}. Our contribution relies on a
nonparametric/nonlinear approach based on an original proximity
criterion to combine estimators. In that sense, it is different
from existing techniques.
\medskip

Indeed, the present article investigates a novel point of view,
motivated by the sense that nonlinear, data-dependent techniques are a
source of analytic flexibility. Instead of forming a linear combination of
estimators, we propose an original nonlinear method for combining the
outcomes over some list of candidate procedures. We call this
combined scheme a regression collective over the given basic machines. 
We consider the problem of building a new estimator by
combining  $M$ estimators of the regression function, thereby exploiting an
idea proposed in the context of supervised classification by \citet{Moj1999}.
Given a set of preliminary estimators $r_1,\dots,r_M$, the
idea behind this combining method is a ``unanimity''
concept, which is based on the values predicted by
$r_1,\dots,r_M$ for the data and for a new observation $\bx$. In a nutshell, a data point is considered to be ``close'' to $\bx$, and consequently, reliable for contributing to the estimation of this new observation, if all estimators 
predict values which are close to each
other for $\bx$ and this data item, \ie, not more distant than a prespecified threshold $\e$. The
predicted value corresponding to this query point $\bx$ is then set to the average of the responses of
the  selected observations. Let us stress here that the average is over the original outcome values of the selected observations, and \emph{not} 
over the estimates provided by the several machines for these observations.
\medskip

\begin{figure}[h!]
  \caption{A toy example: Combining two primal estimators.}
  \label{toy}
  \subfloat[How should we predict the response for the query point $\bx$ (dotted line)?]{\includegraphics[width=.49\textwidth]{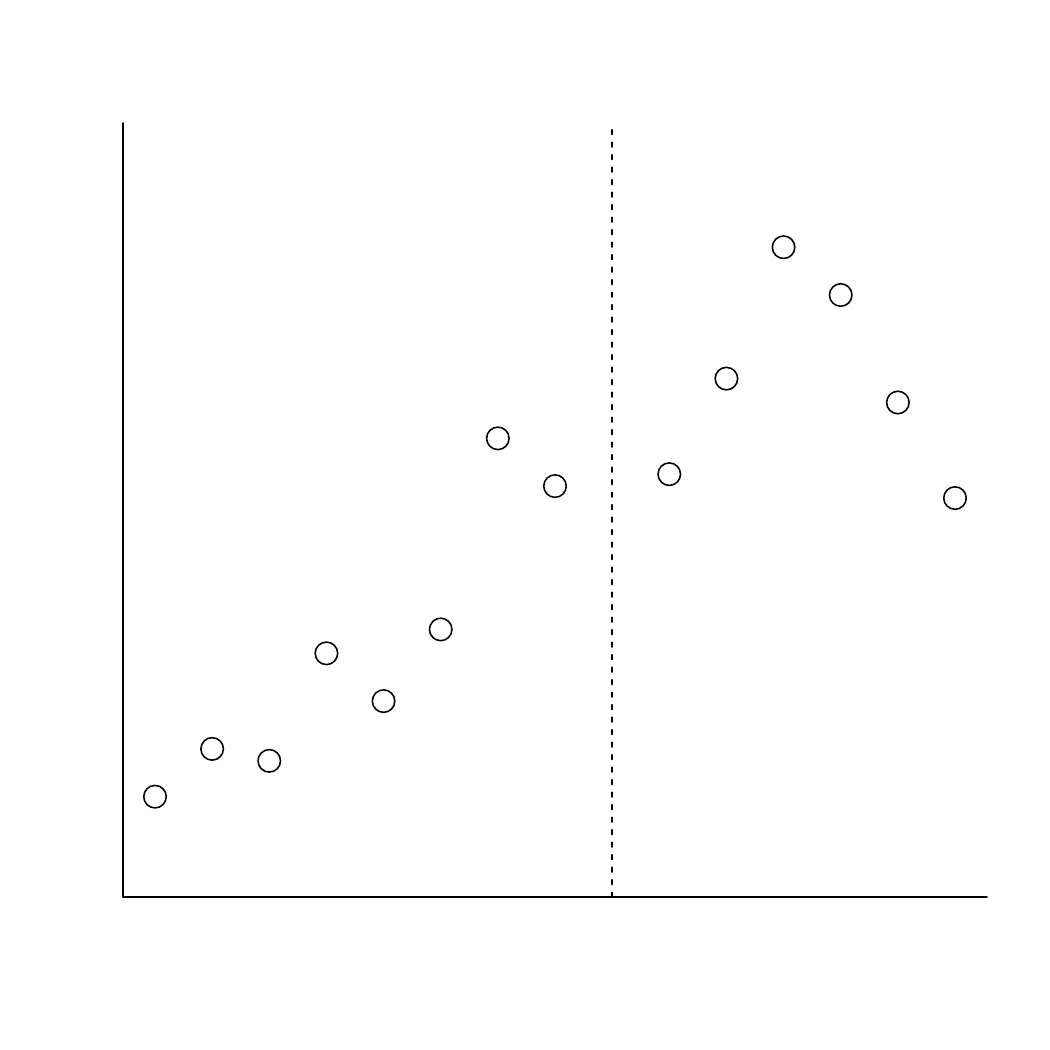}}
  \hfill
  \subfloat[The two primal estimators $r_1$ and $r_2$.]{\includegraphics[width=.49\textwidth]{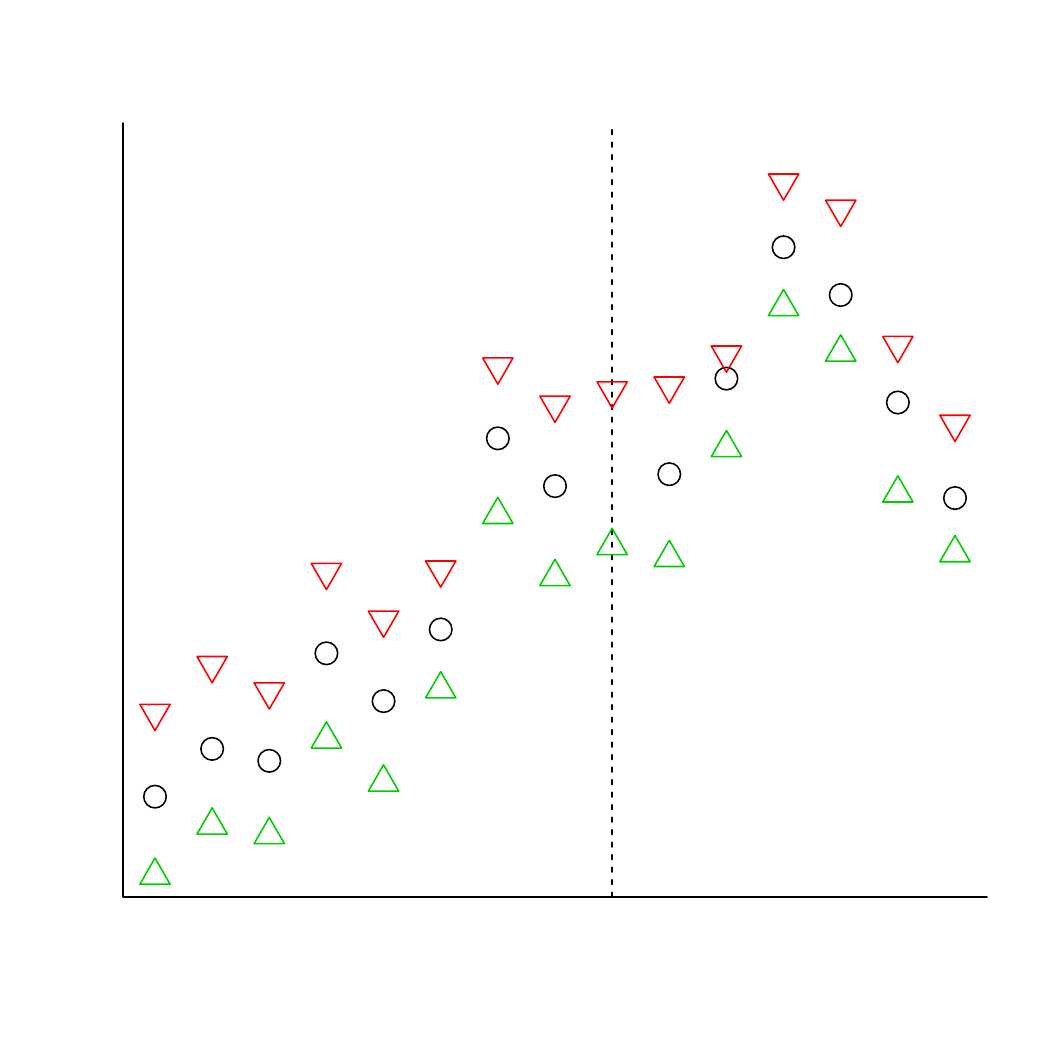}} \\
  \hfill
  \subfloat[The collective operates.]{\includegraphics[width=.49\textwidth]{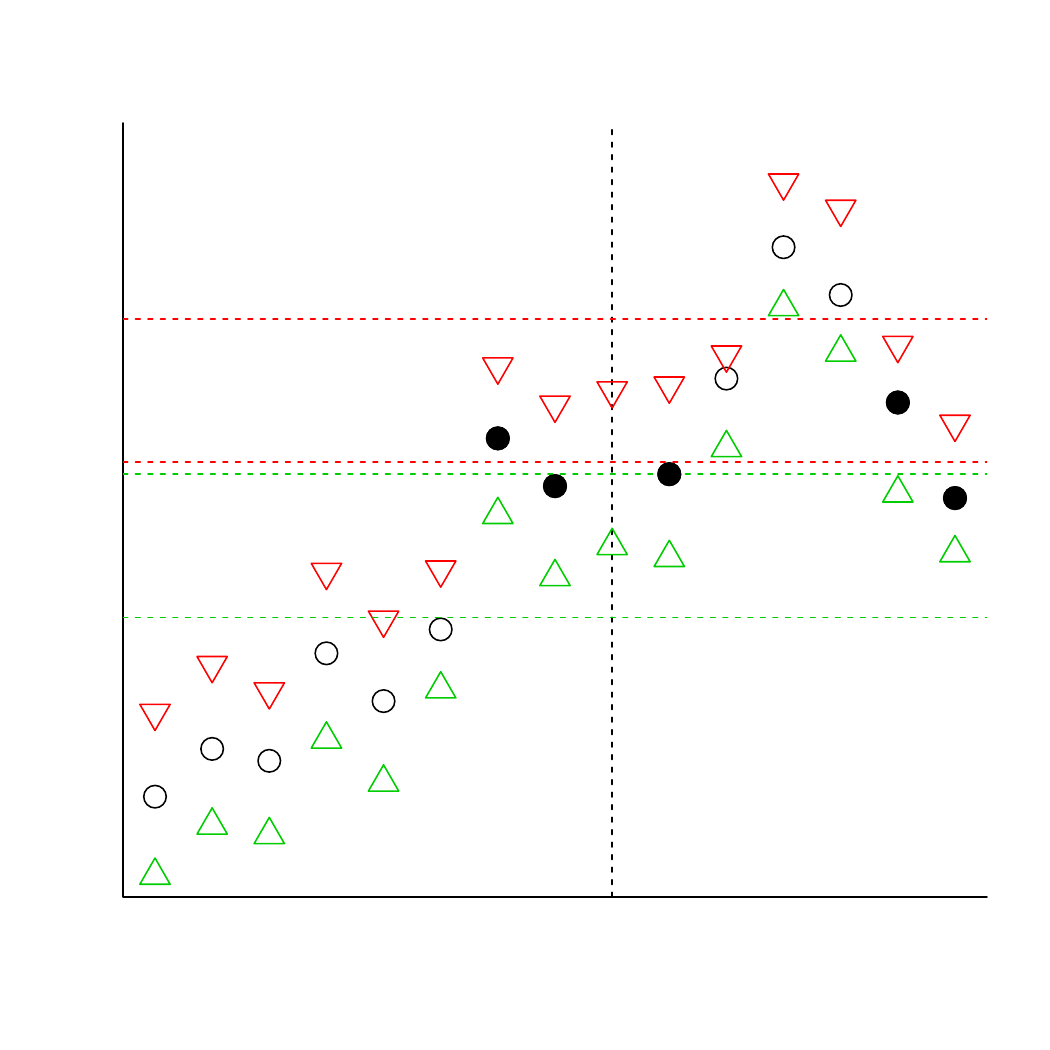}}
  \hfill
  \subfloat[Predicted value (diamond) for the query point $\bx$.]{\includegraphics[width=.49\textwidth]{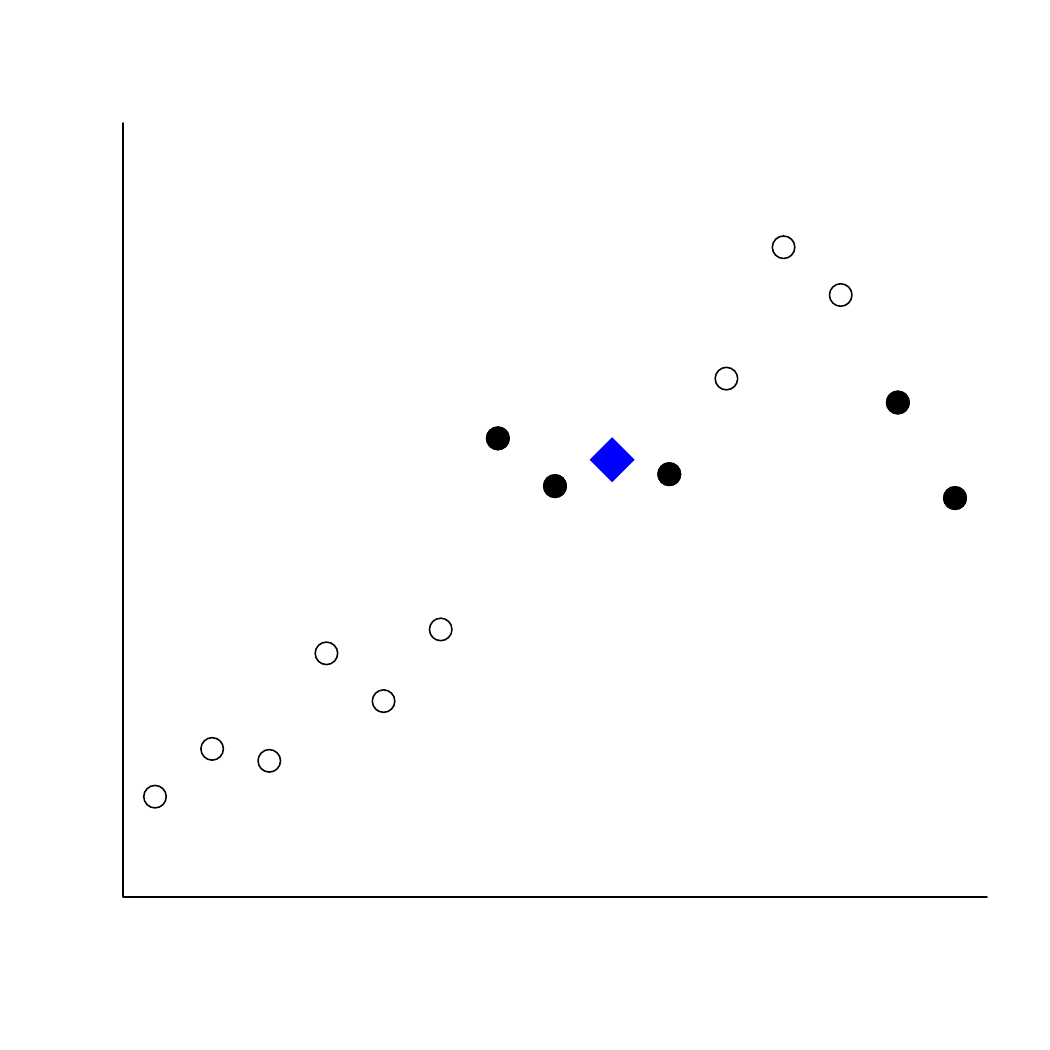}}
\end{figure}

To make the concept clear, consider the following toy
example illustrated by \autoref{toy}. Assume we are given the observations plotted in circles, and the values predicted by two known machines $r_1$
and $r_2$ (triangles pointing up and down, respectively). The goal is to predict the response
for the new point $\bx$ (along the dotted line). Setting a threshold $\e$, the black solid
circles are the data points $(\bx_i,y_i)$ within the two dotted
intervals, \ie, such that for $m=1,2$,
$|r_m(\bx_i)-r_m(\bx)|\leq\e$. Averaging the corresponding $y_i$'s yields the
prediction for $\bx$ (diamond).
\medskip

We stress that the central and original idea behind our approach is
that the resulting regression predictor is a nonlinear, nonparametric, data-dependent
function of the basic predictors $r_1,\dots,r_M$, where the
predictors are used to determine a local distance between a new test
instance and the original training data. To the best of our
knowledge there exists no formalized procedure in the machine learning
and aggregation literature that operates as ours does. In particular, note that the original nonparametric nature of our combined estimator opens up new perspectives of research.
\medskip

Indeed, though we have in mind a batch setting where the data collected consists in an $n$-sample of i.i.d. replications of some variable $(\bX,Y)$, our procedure may be linked to other situations. For example, consider the case of functional data analysis (see \citealp[][]{FV2006}, and \citealp[][]{BSGV2014}, for a survey on recent developments). Even though our method is fitted for finite dimensional data, it may be naturally extended to functional data after a suitable preprocessing of the curves. For example, this can be achieved using an expansion of the curves on an appropriate functional dictionary, and/or \emph{via} a variable selection approach, as in \citet{AV2014}. Note that in a recent work, \citet{CFKL2015} adapts our procedure in a classification setting, also in a functional example.
\medskip

Along with this paper, we release the
software \cobra \citep{Gue2013} which implements
the method as an additional package to the statistical
software R \citep[see][]{CTea2012}. \cobra is freely downloadable on the
CRAN
website\footnote{\url{http://cran.r-project.org/web/packages/COBRA/index.html}}. 
As detailed in \autoref{section:cobra}, we undertook a lengthy series of
numerical experiments, over which \cobra proved extremely successful. These stunning results lead us to believe
that regression collectives can provide valuable insights on a wide
range of prediction problems. 
Further, these same results demonstrate
that \cobra has remarkable speed in terms of CPU timings. In the
context of high-dimensional (such as genomic) data, such velocity is
critical, and in fact \cobra can natively take advantage of multi-core
parallel environments.
\medskip

The paper is organized as follows. In \autoref{section:theorie}, we describe the
combined estimator---the regression collective---and derive a
nonasymptotic risk bound. Next we present the main result, that is, the
collective is asymptotically at least as good as any functional of the basic
estimators. We also provide a rate of convergence for our procedure. \autoref{section:cobra} is devoted to the companion R package
\cobra and presents benchmarks of its excellent performance on both
simulated and real data sets, including high-dimensional models. We also show that \cobra compares
favorably with two competitors, Super Learner \citep[][]{vdLPH2007} and the exponentially
weighted aggregate \citep[see for example][]{Gir2014}, in that it performs
similarly in most situations, much better in some, while it is
consistently faster than the Super Learner in every case. Finally, for ease of
exposition, proofs and additional simulation results (figures and tables with (SM) as suffix) are postponed to a Supplementary Material.

\section{The combined estimator}\label{section:theorie}

\subsection{Notation}

Throughout the article, we assume that we are given a training sample
denoted by $\mathcal D_n=\{(\bX_1,Y_1), \dots,
(\bX_n,Y_n)\}$. $\mathcal{D}_n$ is composed
of i.i.d. random variables taking
their values in $\mathbb R^d\times \mathbb R$, and distributed as an
independent prototype pair $(\bX, Y)$ satisfying $\mathbb E Y^2 <
\infty$ (with the notation $\bX=(X_1,\dots,X_d)$). The space $\mathbb
R^d$ is equipped with the standard Euclidean metric. Our goal is to consistently estimate the regression
function $r^{\star}(\bx)=\mathbb E[Y|\bX=\bx]$, $\bx \in \mathbb R^d$, using the data $\mathcal D_n$. 
\medskip

To begin with, the original data set $\mathcal D_n$ is split into two data
sequences $\mathcal{D}_k=\{(\bX_1,Y_1),\dots,(\bX_k,Y_k)\}$ and $\mathcal{D}_\ell=\{(\bX_{k+1},Y_{k+1}),\dots,(\bX_n,Y_n)\}$, with $\ell=n-k\geq 1$. For ease of notation, the elements of $\mathcal{D}_\ell$ are renamed $\{(\bX_1,Y_1),\dots,(\bX_\ell,Y_\ell)\}$. There is a slight abuse of notation here, as the same letter is used for both subsets $\mathcal D_k$ and $\mathcal D_{\ell}$---however, this should not cause any trouble since the context is clear.
\medskip

Now, suppose that we are given a collection of $M\geq 1$ competing
candidates $r_{k,1},\dots,r_{k,M}$ to estimate $r^{\star}$.
These basic estimators---basic machines---are assumed to be generated
using only the first subsample $\mathcal{D}_k$. These machines can be
any among the researcher's favorite toolkit, such as linear
regression, kernel smoother, SVM, Lasso, neural networks, naive Bayes, or
random forests. They could equally well be any ad hoc regression rules
suggested by the experimental context. The essential idea is that
these basic machines can be parametric, nonparametric, or
semi-parametric, with possible tuning rules. All that is asked for is
that each of the $r_{k,m}(\bx)$, $m=1, \dots, M$, is able to provide
an estimation of $r^{\star}(\bx)$ on the basis of $\mathcal D_k$
alone. Thus, any collection of model-based or model-free machines are
allowed, and our way of combining such a collection is here called the regression
collective. Let us emphasize that the number of basic machines
$M$ is considered as fixed throughout this paper. Hence, the number of machines is not expected to grow and is typically of a
reasonable size ($M$ is chosen on the order of $10$ in \autoref{section:cobra}).
\medskip

Given the collection of basic machines
$\br_k=(r_{k,1},\dots,r_{k,M})$, we define the collective estimator
$T_n$ to be
$$ T_n\left (\br_k(\bx)\right)=\sum_{i=1}^\ell W_{n,i}(\bx) Y_i,\quad \bx \in \mathbb R^d,$$ 
where the random weights $W_{n,i}(\bx)$ take the form
\begin{equation}W_{n,i}(\bx)=\frac{\1_{\bigcap_{m=1}^M\{|r_{k,m}(\bx)-r_{k,m}(\bX_i)|\leq \varepsilon_{\ell}\}}}{\sum_{j=1}^\ell \1_{\bigcap_{m=1}^M\{|r_{k,m}(\bx)-r_{k,m}(\bX_j)|\leq {\varepsilon}_\ell\}}}.\label{eq:weights}\end{equation}
In this definition, ${\varepsilon}_\ell$ is some positive parameter and, by convention, $0/0=0$. 
\medskip

The weighting scheme used in our regression collective is distinctive
but not obvious. Starting from \citet{DGL1996} and \citet{GKKW2002}, we see that $T_n$ is
a local averaging estimator in the following sense: The predicted value for
$r^\star(\bx)$, that is, the estimated outcome at the query point $\bx$, is the
unweighted average over those $Y_i$'s such that $\bX_i$ is ``close'' to the query
point. More precisely, for each $\bX_i$ in the sample $\mathcal{D}_\ell$, ``close'' means that the output at the query
point, generated from each basic machine, is within an $\e_\ell$-distance of the output
generated by the same basic machine at $\bX_i$. If a basic machine evaluated at $\bX_i$ is close to the basic machine
evaluated at the query point $\bx$, then the corresponding outcome $Y_i$ is
included in the average, and not otherwise. Also, as a further note of
clarification: ``Closeness'' of the $\bX_i$'s is not here to be
understood in the Euclidean
sense. It refers to closeness
of the primal estimators outputs at the query point as compared to the outputs over all points in the training data. Training points
$\bX_i$ that are close, in this sense, to the corresponding
outputs at the query point contribute to the indicator
function for the corresponding outcome $Y_i$. This alternative approach is
motivated by the fact that a major issue in learning problems consists
of devising a metric that is suited to the data \citep[see,
\emph{e.g.}, the monograph by][]{PD2005}.
\medskip

In this context,
$\varepsilon_{\ell}$ plays the role of a smoothing parameter: Put
differently, in order
to retain $Y_i$, all basic estimators $r_{k,1},\dots,r_{k,M}$ have to
deliver
predictions for the query point $\bx$ which are in a
$\e_\ell$-neighborhood of the predictions
$r_{k,1}(\bX_i),\dots,r_{k,M}(\bX_i)$. Note that the greater
$\e_\ell$, the more tolerant the process. It turns out that the practical performance of $T_n$ strongly
relies on an appropriate choice of $\varepsilon_{\ell}$. This
important question will be discussed in
\autoref{section:cobra}, where we devise an automatic (\ie, data-dependent) selection strategy of $\varepsilon_{\ell}$.
\medskip

Next, we note that the subscript $n$ in $T_n$ may be a little
confusing, since $T_n$ is a weighted average of the $Y_i$'s in
$\mathcal D_{\ell}$ only. However, $T_n$ depends on the entire data
set $\mathcal D_n$, as the rest of the data is used to set up the
original machines $r_{k,1}, \dots, r_{k,M}$. Most
importantly, it should be noticed that the combined estimator $T_n$ is
nonlinear with respect to the basic estimators $r_{k,m}$. As such, it is inspired by the preliminary work of
\citet{Moj1999} in the supervised classification context.
\medskip

In addition, let us mention that, in the definition of the weights 
\eqref{eq:weights}, all original estimators are invited to have the
same, equally valued
opinion on the importance of the observation $\bX_i$ (within the range
of $\varepsilon_{\ell}$) for the corresponding $Y_i$ to be integrated
in the combination $T_n$. However, this unanimity constraint may be
relaxed  by imposing, for example, that a fixed fraction $\alpha\in \{1/M,2/M,\dots,1\}$ of the machines agrees on the importance of $\bX_i$. In that case, the weights take the more sophisticated form
$$W_{n,i}(\bx)=\frac{\1_{\{\sum_{m=1}^M \1_{\{|r_{k,m}(\bx)-r_{k,m}(\bX_i)|\leq \varepsilon_{\ell}\}}\geq M\alpha\}}}{\sum_{j=1}^\ell {\1_{\{\sum_{m=1}^M \1_{\{|r_{k,m}(\bx)-r_{k,m}(\bX_j)|\leq \varepsilon_{\ell}\}}\geq M\alpha\}}}}.
$$
It turns out that adding the parameter $\alpha$ does not change the
asymptotic properties of $T_n$, provided $\alpha\to 1$. Thus, to keep
a sufficient degree of clarity in the mathematical statements and
subsequent proofs, we have decided to consider only the case
$\alpha=1$ (\emph{i.e.}, unanimity). Extension of the results to more general
values of $\alpha$ is left for future work. On the other hand, as highligthed by \autoref{section:cobra}, $\alpha$ has a nonnegligible impact on the performance of the combined estimator. Accordingly, we will discuss in \autoref{section:cobra} an automatic procedure to select this extra parameter.
\subsection{Theoretical performance}
This section is devoted to the study of some asymptotic and nonasymptotic properties of the combined estimator $T_n$, whose quality will be assessed by the quadratic risk 
$$\mathbb E\left |T_n\left (\br_k(\bX)\right)-r^\star(\bX)\right|^2.$$

Here and later, $\mathbb E$ denotes the expectation with respect to both $\bX$ and the sample $\mathcal D_n$.
Everywhere in the document, it is assumed that $\E
|r_{k,m}(\bX)|^2 <\infty$ for all $m=1,\dots,M$.
\medskip

For any $m=1,\dots,M$, let $r_{k,m}^{-1}$ denote the
inverse image of machine $r_{k,m}$. Assume that for any $m=1,\dots,M$,
\begin{equation}\label{bound}
  r_{k,m}^{-1}((t,+\infty)) \underset{t\uparrow +\infty}{\searrow}\emptyset \quad\mbox{and}\quad r_{k,m}^{-1}((-\infty,t)) \underset{t\downarrow -\infty}{\searrow}\emptyset.
\end{equation}
It is stressed that this is a mild assumption which is met, for example, whenever the machines are bounded.
 Throughout, we let
$$T\left (\br_k(\bX)\right)=\mathbb E\left [Y|\br_k(\bX)\right]$$
and note that, by the very definition of the $L^2$ conditional expectation,
\begin{equation}
\mathbb E\left |T(\br_{k}(\bX))-Y\right|^2\leq \inf_{f}\ \mathbb E\left |f(\br_{k}(\bX))-Y\right|^2,
\end{equation}
where the infimum is taken over all square integrable functions of $\br_k(\bX)$.
\medskip

Our first result is a nonasymptotic inequality, which states that the
combined estimator behaves as well as the best one in the original
list, within a term measuring how far $T_n$ is from $T$.

\begin{pro}
\label{theo:orac} Let $\br_k=(r_{k,1}, \dots, r_{k,M})$ be the
collection of basic estimators, and let $T_n(\br_k(\bx))$ be the
combined estimator. Then, for all distributions of $(\bX,Y)$ with $\mathbb E Y^2<\infty$,
\begin{align*}
\mathbb E|&T_n(\br_k(\bX))-r^\star(\bX)|^2 \\ &\leq \mathbb E|T_n(\br_k(\bX))-T(\br_k(\bX))|^2+\inf_{f}\ \mathbb E|f(\br_k(\bX))-r^\star(\bX)|^2,
\end{align*}
where the infimum is taken over all square integrable functions of $\br_k(\bX)$. In particular,
\begin{align*}
&\mathbb E|T_n(\br_k(\bX))-r^\star(\bX)|^2\\
& \quad \leq \min_{m=1,\dots,M}\mathbb E|r_{k,m}(\bX)-r^\star(\bX)|^2 +\mathbb E|T_n(\br_k(\bX))-T(\br_k(\bX))|^2.
\end{align*}
\end{pro}

\autoref{theo:orac} guarantees the performance of $T_n$ with
respect to the basic machines, whatever the distribution of $(\bX,Y)$
is and regardless of which initial estimator is actually the best. The
term $\min_{m=1,\dots,M}\E|r_{k,m}(\bX)-r^\star(\bX)|^2$ may be
regarded as a bias term, whereas the term $\mathbb E |T_n
(\br_k(\bX))-T(\br_k(\bX))|^2$ is a variance-type term, which can be
asymptotically neglected, as shown by the following result.
 \begin{pro}
 \label{prop:cons}
Assume that $\varepsilon_{\ell}\to 0$ and $ \ell \varepsilon_{\ell}^M
\to \infty$ as $\ell \to \infty.$
Then
$$\mathbb E\left |T_n\left (\br_{k}(\bX)\right)-T\left(\br_{k}(\bX)\right)\right|^2 \to 0 \quad \mbox{as } \ell \to \infty,$$ 
for all distributions of $(\bX,Y)$ with $\mathbb EY^2<\infty$.
Thus,
$$\limsup_{\ell \to \infty} \mathbb E \left |T_n\left
(\br_k(\bX)\right)-r^\star(\bX)\right|^2 \leq \inf_{f}\ \mathbb E\left
|f(r_{k,m}(\bX))-r^\star(\bX)\right|^2.$$
In particular,
$$\limsup_{\ell \to \infty} \mathbb E \left |T_n\left (\br_k(\bX)\right)-r^\star(\bX)\right|^2 \leq \min_{m=1,\dots,M}\mathbb E\left |r_{k,m}(\bX)-r^\star(\bX)\right|^2.$$
\end{pro}
This result is remarkable, for two reasons. Firstly, it shows
that, in terms of predictive quadratic risk, the combined estimator
does asymptotically at least as well as the best primitive
machine. Secondly, the result is nearly universal, in the sense that it is
true for all distributions of $(\bX,Y)$ such that $\E Y^2<\infty$.
\medskip

This is especially interesting because the performance of any
estimation procedure eventually depends upon some model and smoothness
assumptions on the observations. For example, a linear regression fit
performs well if the distribution is truly linear, but may behave
poorly otherwise. Similarly, the Lasso procedure is known to do a good
job for non-correlated designs, with no clear guarantee however in
adversarial situations. Likewise, performance of
nonparametric procedures such as the $k$-nearest neighbor method,
kernel estimators and random forests dramatically deteriorate as the ambient dimension
increases, but may be significantly improved if the true underlying
dimension is reasonable. Note that this phenomenon is thoroughly analyzed for
the random forests algorithm in \citet{Bia2012}.
\medskip

The result exhibited in \autoref{prop:cons} holds
under a minimal regularity assumption on the basic machines. However, this
universality comes at a price since we have no guarantee on the rate of
convergence of the variance term. Nevertheless,
assuming some light additional smoothness conditions, one has the
following result, which is the central statement of the paper.
\begin{theo}\label{pro:vitesse}
Assume that $Y$ and the basic machines $\br_k$ are bounded by some constant
$R$. Assume moreover that there exists a constant $L\geq 0$ such that, for every $k\geq 1$, 
$$ 
 |T(\br_k(\bx))-T(\br_k(\by))|\leq L|\br_k(\bx)-\br_k(\by)|,\quad \bx,\by\in\R^d.$$
Then, with the choice $\e_\ell\propto \ell^{-\frac{1}{M+2}}$, one has
 \begin{equation*}
    \mathbb E \left |T_n\left (\br_k(\bX)\right)-r^\star(\bX)\right|^2 \leq \min_{m=1,\dots,M}\mathbb E\left |r_{k,m}(\bX)-r^\star(\bX)\right|^2
    +C\ell^{-\frac{2}{M+2}},
  \end{equation*}
for some positive constant $C=C(R,L)$, independent of $k$.
\end{theo}
\autoref{pro:vitesse} offers an oracle-type inequality with
leading constant $1$ (\emph{i.e.}, sharp oracle inequality), stating that
the risk of the regression collective is bounded by the lowest risk
among those of the basic machines, \emph{i.e.}, our procedure mimics
the performance of the oracle over the set $\{r_{k,m} \colon m=1,\dots,M\}$, plus a remainder term
of the order of $\ell^{-2/(M+2)}$ which is the price to pay for
combining $M$ estimators. In our setting, it is important to observe that this term has a
limited impact. As a matter of fact, since the number of basic machines $M$ is
assumed to be fixed and not too large (the implementation presented in
\autoref{section:cobra} considers $M$ at most $6$), the remainder term
is negligible compared to the standard
nonparametric rate $\ell^{-2/(d+2)}$ in dimension $d$. While the rate
$\ell^{-2/(d+2)}$ is affected by the curse of dimensionality when $d$
is large, this is not the case for the term
$\ell^{-2/(M+2)}$. That way, our procedure appears well armed to face high dimensional problems. When $d\gg n$, many methods deteriorate and suffer from the curse of dimensionality. However, it is important to note here that even if some of the basic machines $r_{k,1},\dots,r_{k,M}$ might be less performant in that context, this does not affect in any way our combining procedure. Indeed, forming the regression collective $T_n$ does not require any additional effort if $d$ grows. Obviously, when $d$ is large, the best choice would be to include as basic machines methods and models which are adapted to the high dimensional setting. This is an interesting track for future research, which is connected to functional data analysis and dimension-reduction models \citep[see][]{GV2014}.
\medskip

Obviously, under the assumption that the distribution of
$(\bX,Y)$ might be described parametrically and that one of the initial
estimators is adapted to this distribution, faster rates of the
order of $1/\ell$ could emerge in the bias term.
Nonetheless, the regression
collective is designed for much more adversarial regression problems,
hence the rate exhibited in \autoref{pro:vitesse} appears
satisfactory. We stress that our
approach carries no assumption on the random design and mild ones over the
primal estimators, in line with our attempt to design a procedure
which is as model-free as
possible.
\medskip

The central motivation for our method
is that model and smoothness assumptions are usually unverifiable,
especially in modern high-dimensional and large scale data sets. To
circumvent this difficulty, researchers often try many different methods
and retain the one exhibiting the best empirical (\emph{e.g.}, cross-validated) results. Our
combining strategy offers a nice alternative, in the sense that if
one of the initial estimators is  consistent for a given class $\mathcal M$ of distributions, then, under light smoothness assumptions, $T_n$ inherits the same
property. To be more precise, assume that the initial pool of estimators
includes a consistent estimator, \emph{i.e.}, that one of the original estimators, say $r_{k,m_0}$, satisfies
$$\mathbb E \left |r_{k, m_0}(\bX)-r^\star(\bX)\right|^2\to 0 \quad \mbox{as } k\to \infty,$$
for all distributions of $(\bX,Y)$ in some class
$\mathcal{M}$. 
Then, under the assumptions of \autoref{pro:vitesse}, with the choice $\e_\ell\propto \ell^{-\frac{1}{M+2}}$, one has
\begin{equation*}
\lim_{k,\ell\to\infty}\mathbb E \left |T_n\left
    (\br_k(\bX)\right)-r^\star(\bX)\right|^2 = 0.
\end{equation*}

\section{Implementation and numerical studies}\label{section:cobra}

This section is devoted to the implementation of the described method.
Its excellent performance is then assessed in a series of experiments. The companion R package
\cobra (standing for COmBined Regression Alternative) is
available on the CRAN website\footnote{\url{http://cran.r-project.org/web/packages/COBRA/index.html}}, for
Linux, Mac and Windows platforms \citep[see][]{Gue2013}. Note that in a will to favor its execution speed, \cobra includes a \texttt{parallel}
option, allowing for improved performance on multi-core computers \citep[from][]{Kna2010}.
\medskip

As raised in the previous section, a precise calibration of the smoothing parameter
$\e_\ell$ is crucial. Clearly, a value that is too small will discard many machines
and most weights will be zero. Conversely, a large value sets all
weights to $1/\Sigma$ with $$\Sigma=\sum_{j=1}^\ell
\1_{\bigcap_{m=1}^M\{|r_{k,m}(\bx)-r_{k,m}(\bX_j)|\leq
  {\varepsilon}_\ell\}},$$ giving the naive predictor that does not
account for any new data point and predicts the mean over the sample
$\mathcal{D}_\ell$. We also consider a relaxed version of the unanimity
constraint: Instead of requiring global agreement over the implemented
machines, consider some $\alpha\in (0,1]$ and keep observation $Y_i$ in the
construction of $T_n$ if and only if at least a proportion $\alpha$ of the machines
agrees on the importance of $\bX_i$. This parameter requires some
calibration. To understand this better, consider the following toy example: On some data set,
assume most machines but one have nice predictive performance. For any
new data point, requiring global agreement will fail
since the pool of machines is heterogeneous. In this regard, $\alpha$
should be seen as a measure of homogeneity: If a small value is
selected, it may be an indicator that some
machines perform (possibly much) better than some others. Conversely,
a large value indicates that the predictive abilities of the machines
are close.
\medskip

A natural measure of the risk in the prediction
context is the empirical
quadratic loss, namely
$$\hat{R}(\hat{\mathbf{Y}})=\frac{1}{p}\sum_{j=1}^p
(\hat{Y}_j-Y_j)^2,$$
where $\hat{\mathbf{Y}}=(\hat{Y}_1,\dots,\hat{Y}_p)$ is the vector
of predicted values for the responses $Y_1,\dots,Y_p$ and
$\{(\bX_j,Y_j)\}_{j=1}^p$ is a testing sample.
We adopted the following protocol: Using a simple data-splitting
device, $\e_\ell$ and $\alpha$ are chosen by minimizing the empirical risk
$\hat{R}$ over the set
$\{\e_{\ell,\mathrm{min}},\dots,\e_{\ell,\mathrm{max}}\}\times\{1/M,\dots,1\}$,
where $\e_{\ell,\mathrm{min}}=10^{-300}$ and
$\e_{\ell,\mathrm{max}}$ is proportional to the largest absolute difference between two
predictions of the pool of machines.
\medskip

In the package, the number
$\#\{\e_{\ell,\mathrm{min}},\dots,\e_{\ell,\mathrm{max}}\}$ of
evaluated values may be
modified by the user, otherwise the default value $200$ is
chosen. It is also possible to choose either a linear or a logistic
scale. \autoref{epscal} illustrates the discussion about the choice
of $\e_\ell$ and $\alpha$.
\medskip

By default, \cobra includes the following classical packages dealing
with regression estimation and prediction. However, note that
the user has the choice to modify this list to her/his own convenience:
\begin{itemize}
\item Lasso \citep[R package \texttt{lars}, see][]{HE2012}.
\item Ridge regression \citep[R package \texttt{ridge}, see][]{Cul2012}.
\item $k$-nearest neighbors \citep[R package \texttt{FNN}, see][]{Li2012}.
\item CART algorithm \citep[R package \texttt{tree}, see][]{Rip2012}.
\item Random Forests algorithm \citep[R package \texttt{randomForest}, see][]{LW2002}.
\end{itemize}
First, \cobra is benchmarked on synthetic data. For each of the
following eight models, two designs are considered: Uniform over
$(-1,1)^d$ (referred to as ``Uncorrelated'' in \autoref{tabMachines},
\autoref{tabSL} and \autoref{tabTime}), and Gaussian with mean $0$ and covariance matrix
$\Sigma$ with $\Sigma_{ij}=2^{-|i-j|}$ (``Correlated''). Models
considered cover a wide spectrum of contemporary regression problems. Indeed,
\autoref{m1} is a toy example, \autoref{m2} comes from \citet{vdLPH2007}, \autoref{m3} and \autoref{m4} appear in
\citet{MvdGB2009}. \autoref{m5} is somewhat a classic setting. \autoref{m6} is about predicting labels, \autoref{m7}
is inspired by high-dimensional sparse regression problems. Finally,
\autoref{m8} deals with probability estimation, forming a link with
nonparametric model-free approaches such as in \citet{MKD+2012}. In the
sequel, we let $\mathcal{N}(\mu,\sigma^2)$ denote a Gaussian random variable
with mean $\mu$ and variance $\sigma^2$. In the simulations, the
training data set was usually set to $80\%$ of the whole sample, then
split into two equal parts corresponding to $\mathcal{D}_k$ and $\mathcal{D}_\ell$.
\begin{model}\label{m1}
  $n=800$, $d=50$, $Y = X_1^2+\exp(-X_2^2)$.
\end{model}
\begin{model}\label{m2}
  $n=600$, $d=100$, $Y = X_1X_2+X_3^2-X_4X_7+X_8X_{10}-X_6^2+\mathcal{N}(0,0.5)$.
\end{model}
\begin{model}\label{m3}
  $n=600$, $d=100$, $Y = -\sin(2X_1)+X_2^2+X_3-\exp(-X_4)+\mathcal{N}(0,0.5)$.
\end{model}
\begin{model}\label{m4}
  $n=600$, $d=100$, $Y = X_1+(2X_2-1)^2+\sin(2\pi X_3)/(2-\sin(2\pi
  X_3))+\sin(2\pi X_4)+2\cos(2\pi X_4)+3\sin^2(2\pi X_4)+4\cos^2(2\pi X_4)+\mathcal{N}(0,0.5)$.
\end{model}
\begin{model}\label{m5}
  $n=700$, $d=20$, $Y = \1_{\{X_1>0\}}+X_2^3+\1_{\{X_4+X_6-X_8-X_9>1+X_{14}\}}+\exp(-X_2^2)+\mathcal{N}(0,0.5)$.
\end{model}
\begin{model}\label{m6}
  $n=500$, $d=30$, $Y = \sum_{k=1}^{10}\1_{\{X^3_k<0\}}-\1_{\{\mathcal{N}(0,1)>1.25\}}$.
\end{model}
\begin{model}\label{m7}
  $n=600$, $d=300$, $Y = X_1^2+X_2^2X_3\exp(-|X_4|)+X_6-X_8+\mathcal{N}(0,0.5)$.
\end{model}
\begin{model}\label{m8}
  $n=600$, $d=50$, $Y = \1_{\{X_1+X_4^3+X_9+\sin(X_{12}X_{18})+\mathcal{N}(0,0.1)>0.38\}}$.
\end{model}
\autoref{tabMachines} presents the empirical mean quadratic error and standard
deviation over $100$ independent replications, for each model and
design. Bold numbers identify the lowest error, \ie, the apparent best
competitor. Boxplots of errors are presented in \autoref{boxplot-U} and \autoref{boxplot-C}. Further,
\autoref{pred-U} and \autoref{pred-C} show the predictive capacities of \cobra,
and \autoref{func} depicts its ability to reconstruct the
functional dependence over the covariates in the context of additive regression, assessing the striking performance of our approach in a wide
spectrum of statistical settings. A persistent and notable
fact is that \cobra performs at least as well as the best machine, especially so in \autoref{m3}, \autoref{m5} and
\autoref{m6}.
\medskip

Next, since more and more problems in contemporary statistics
involve high-dimensional data, we have tested the abilities of \cobra in
that context. As highlighted by \autoref{tab-HD} and \autoref{box-HD},
the main message is that \cobra is perfectly able to
deal with high-dimensional data, provided that it is generated over
machines, at least some of which are
known to perform well in such situations (possibly at the price of a
sparsity assumption). In that context, we conducted
$200$ independent replications for the three following models:
\begin{model}\label{mHD1}
  $n=500$, $d=1000$, $Y=X_1+3X_3^2-2\exp(-X_5)+X_6$. Uncorrelated design.
\end{model}
\begin{model}\label{mHD2}
  $n=500$, $d=1000$, $Y=X_1+3X_3^2-2\exp(-X_5)+X_6$. Correlated design.
\end{model}
\begin{model}\label{mHD3}
  $n=500$, $d=1500$, $Y=\exp(-X_1)+\exp(X_1)+\sum_{j=2}^dX_j^{j/100}$. Uncorrelated design.
\end{model}

A legitimate question that arises is where one should cut the initial sample
$\mathcal{D}_n$? In other words, for a given data set of size $n$, what is the optimal value for $k$? A naive approach is to cut the initial sample in two
halfs (\emph{i.e.}, $k=n/2$): This appears to be satisfactory provided that
$n$ is large enough, which may be too much of an unrealistic
assumption in numerous experimental settings. A more involved choice is to adopt a random cut scheme, where
$k$ is chosen uniformly in $\{1,\dots,n\}$. \autoref{box-Stab}
presents the boxplots of errors of the five default machines and
\cobra with that random cutting strategy, and also shows
the risk of \cobra with respect to $k$. To illustrate this phenomenon,
we tested a thousand random cuts on the following \autoref{mSTAB}. As showed in
\autoref{box-Stab}, for that particular model, the best value seems to
be near $3n/4$.
\begin{model}\label{mSTAB}
  $n=1200$, $d=10$, $Y=X_1+3X_3^2-2\exp(-X_5)+X_6$. Uncorrelated design.
\end{model}
The average risk of \cobra on
a thousand replications of \autoref{mSTAB} is $0.3124$. Since this
delivered a thousand prediction vectors, a natural idea is to take their
mean or median. The risk of the mean is $0.2306$, and the median has
an even better risk ($0.2184$). Since a random cut scheme may generate
some unstability, we advise practitioners to compute a few \cobra
estimators, then compute the mean or median vector of their predictions.
\medskip

Next, we compare \cobra to the Super Learner algorithm
\citep[][]{PvdL2012}. This widely used algorithm was first described in
\citet{vdLPH2007} and extended in \citet{PvdL2010}. Super Learner is used in this section as the key
competitor to our method.
In a nutshell, the Super Learner
trains basic machines $r_1,\dots,r_M$ on the whole sample
$\mathcal{D}_n$. Then, following a $V$-fold cross-validation procedure,
Super Learner adopts a $V$-blocks partition of the set $\{1,\dots,n\}$ and
computes the matrix $$H=(H_{ij})_{1\leq i\leq n}^{1\leq j\leq M},$$
where $H_{ij}$ is the prediction for the query
point $\bX_i$ made by machine $j$ trained on all remaining $V-1$
blocks, \emph{i.e.}, excluding the block containing $\bX_i$. The Super Learner estimator is then
\begin{equation*}
SL=\sum_{j=1}^M\hat{\alpha}_jr_j,
\end{equation*}
where
\begin{equation*}
  \hat{\alpha}\in\underset{\alpha\in\Lambda^M}{\arg\inf}\sum_{i=1}^n|Y_i-(H\alpha)_i|^2,
\end{equation*}
with $\Lambda^M$ denoting the
simplex $$\Lambda^M=\left\{\alpha\in\R^M\colon
  \sum_{j=1}^M\alpha_j=1,\ \alpha_j\geq 0 \textrm{ for any }
  j=1,\dots,M\right\}.$$
This convex aggregation scheme is significantly
different from our collective approach. Yet, we feel close to the
philosophy carried by the \texttt{SuperLearner} package, in that both
methods allow the user to aggregate as many
machines as desired, then combining them to deliver predictive
outcomes. For that reason, it is
reasonable to deploy Super Learner as a benchmark in our study of our
collective approach.
\medskip

\autoref{tabSL} 
summarizes the performance of
\cobra and \texttt{SuperLearner} (used with \texttt{SL.randomForest},
  \texttt{SL.ridge} and \texttt{SL.glmnet}, for the fairness of the comparison) through the described protocol. Both
methods compete on similar terms in most models, although \cobra
proves much more efficient on correlated design in \autoref{m2} and
\autoref{m4}. This already remarkable result is to be stressed by the
flexibility and velocity showed by \cobra. Indeed, as emphasized in
\autoref{tabTime}
, without even using the \texttt{parallel}
option, \cobra obtains similar or better results than
\texttt{SuperLearner} roughly five times faster. Note also that \cobra
suffers from a disadvantage: \texttt{SuperLearner} is built on the
whole sample $\mathcal{D}_n$ whereas \cobra only uses $\ell<n$ data
points. Finally, observe that the algorithmic cost of
computing the random weights on $n_{\mathrm{test}}$ query points is
$\ell\times M\times n_{\mathrm{test}}$ operations. In the package, those calculations
are handled in C language for optimal speed performance.
\medskip

Super Learner is a natural competitor on the implementation
side. However, on the theoretical side, we do not assume that it
should be the only benchmark. Thus, we compared \cobra to the popular
exponentially weighted aggregate estimator \citep[EWA, see][]{Gir2014}. We implemented the
following version of the EWA: For all preliminary estimators
$r_{k,1},\dots,r_{k,M}$, their empirical risks
$\hat{R}_1,\dots,\hat{R}_M$ are computed on a subsample of $\mathcal{D}_\ell$ and the EWA is
\begin{equation*}
  \mathrm{EWA}_\beta\colon \bx\mapsto
  \sum_{j=1}^M\hat{w}_jr_{k,j}(\bx),\quad \bx\in\R^d,
\end{equation*}
where
\begin{equation*}
  \hat{w}_j = \frac{\exp(-\beta \hat{R}_j)}{\sum_{i=1}^M \exp(-\beta \hat{R}_i)}, \quad j=1,\dots,M.
\end{equation*}
The temperature parameter $\beta>0$ is selected by minimizing the
empirical risk of $\mathrm{EWA}_\beta$ over a data-based grid, in the same
spirit as the selection of $\e_\ell$ and $\alpha$. We conducted $200$ independent replications, on Models \ref{mHD1} to \ref{mSTAB}. The
conclusion is that \cobra outperforms the EWA estimator in some
models, and delivers similar performance in others, as shown in
\autoref{box-EWA} and \autoref{tab-EWA}.
\medskip

Finally, \cobra is used to process the following real-life data sets:
\begin{itemize}
\item Concrete Slump Test\footnote{\url{http://archive.ics.uci.edu/ml/datasets/Concrete+Slump+Test}.} \citep[see][]{Yeh2007}.
\item Concrete Compressive
  Strength\footnote{\url{http://archive.ics.uci.edu/ml/datasets/Concrete+Compressive+Strength}.}
  \citep[see][]{Yeh1998}.
\item Wine
  Quality\footnote{\url{http://archive.ics.uci.edu/ml/datasets/Wine+Quality}.}
  \citep[see][]{CCA+2009}. We point out that the Wine Quality data set
  involves supervised classification and leads naturally to a line of
  future research using \cobra over
  probability machines \citep[see][]{MKD+2012}.
\end{itemize}
The good predictive performance of \cobra is summarized in
\autoref{real} and errors are presented in \autoref{box-real}. For
every data set, the sample is divided into a training
set ($90\%$) and a testing set ($10\%$) on which the predictive
performance is evaluated. Boxplots are obtained by randomly shuffling the data
points a hundred times.
\medskip

As a conclusion to this thorough experimental protocol, it is our
belief that \cobra sets a
new high standard of reference, a benchmark procedure, both in terms of performance and velocity, for prediction-oriented problems in the context of regression, including high-dimensional problems.

 \section*{Acknowledgements}

 The authors thank the Editor and two anonymous referees for
 providing constructive and helpful remarks, thus greatly improving the paper.

\begin{table}[t]
  \caption{Quadratic errors of the implemented machines and \cobra. Means
    and standard deviations over
    $100$ independent replications.}
  \label{tabMachines}
  \begin{center}
    \begin{tabular}{cc*{6}{c}}
      \textbf{Uncorr.} & & \texttt{lars} & \texttt{ridge} & \texttt{fnn} & \texttt{tree} & \texttt{rf} & \cobra \\ \hline\hline
      \multirow{2}{*}{\autoref{m1}} & m.
      & 0.1561 & 0.1324 & 0.1585 & 0.0281 & 0.0330 & \textbf{0.0259} \\ & sd. & 0.0123 & 0.0094 & 0.0123 & 0.0043 & 0.0033 & 0.0036 \\
      \multirow{2}{*}{\autoref{m2}} & m.  & 0.4880 & 0.2462 & 0.3070 & 0.1746 & \textbf{0.1366} & 0.1645 \\ & sd.  & 0.0676 & 0.0233 & 0.0303 & 0.0270 & 0.0161 & 0.0207 \\
      \multirow{2}{*}{\autoref{m3}} & m. & 0.2536 & 0.5347 & 1.1603 & 0.4954 & 0.4027 & \textbf{0.2332} \\ & sd.  & 0.0271 & 0.4469 & 0.1227 & 0.0772 & 0.0558 & 0.0272 \\
      \multirow{2}{*}{\autoref{m4}} & m.  & 7.6056 & 6.3271 & 10.5890 & 3.7358 & 3.5262 & \textbf{3.3640} \\ & sd.  & 0.9419 & 1.0800 & 0.9404 & 0.8067 & 0.3223 & 0.5178 \\
      \multirow{2}{*}{\autoref{m5}} & m.  & 0.2943 & 0.3311 & 0.5169 & 0.2918 & 0.2234 & \textbf{0.2060} \\ & sd.  & 0.0214 & 0.1012 & 0.0439 & 0.0279 & 0.0216 & 0.0210 \\
      \multirow{2}{*}{\autoref{m6}} & m.  & 0.8438 & 1.0303 & 2.0702 & 2.3476 & 1.3354 & \textbf{0.8345} \\ & sd.  & 0.0916 & 0.4840 & 0.2240 & 0.2814 & 0.1590 & 0.1004 \\
      \multirow{2}{*}{\autoref{m7}} & m. & 1.0920 & 0.5452 & 0.9459 & 0.3638 & 0.3110 & \textbf{0.3052} \\ & sd. & 0.2265 & 0.0920 & 0.0833 & 0.0456 & 0.0325 & 0.0298  \\
      \multirow{2}{*}{\autoref{m8}} & m. & 0.1308 & 0.1279 & 0.2243 & 0.1715 & 0.1236 & \textbf{0.1021} \\ & sd. & 0.0120 & 0.0161 & 0.0189 & 0.0270 & 0.0100 & 0.0155 \\ \\
      \textbf{Corr.} & & \texttt{lars} & \texttt{ridge} & \texttt{fnn} & \texttt{tree} & \texttt{rf} & \cobra \\ \hline\hline
      \multirow{2}{*}{\autoref{m1}} & m. & 2.3736 &
      1.9785 & 2.0958 & 0.3312 & 0.5766 & \textbf{0.3301} \\ & sd.  &
      0.4108 & 0.3538 & 0.3414 & 0.1285 & 0.1914 & 0.1239 \\
      \multirow{2}{*}{\autoref{m2}} & m.  & 8.1710 & 4.0071 & 4.3892 & \textbf{1.3609} & 1.4768 & 1.3612 \\ & sd.  & 1.5532 & 0.6840 & 0.7190 & 0.4647 & 0.4415 & 0.4654 \\
      \multirow{2}{*}{\autoref{m3}} & m.  & 6.1448 & 6.0185 & 8.2154 & 4.3175 & 4.0177 & \textbf{3.7917} \\ & sd.  & 11.9450 & 12.0861 & 13.3121 & 11.7386 & 12.4160 & 11.1806 \\
      \multirow{2}{*}{\autoref{m4}} & m.  & 60.5795 & 42.2117 & 51.7293 & \textbf{9.6810} & 14.7731 & 9.6906 \\ & sd.  & 11.1303 & 9.8207 & 10.9351 & 3.9807 & 5.9508 & 3.9872 \\
      \multirow{2}{*}{\autoref{m5}} & m.  & 6.2325 & 7.1762 & 10.1254 & 3.1525 & 4.2289 & \textbf{2.1743} \\ & sd.  & 2.4320 & 3.5448 & 3.1190 & 2.1468 & 2.4826 & 1.6640 \\
      \multirow{2}{*}{\autoref{m6}} & m.  & 1.2765 &
      1.5307 & 2.5230 & 2.6185 & 1.2027 & \textbf{0.9925} \\ & sd.  &
      0.1381 & 0.9593 & 0.2762 & 0.3445 & 0.1600 & 0.1210 \\
      \multirow{2}{*}{\autoref{m7}} & m. & 20.8575 & 4.4367 & 5.8893 & 3.6865 & \textbf{2.7318} & 2.9127 \\ & sd. & 7.1821 & 1.0770 & 1.2226 & 1.0139 & 0.8945 & 0.9072 \\
      \multirow{2}{*}{\autoref{m8}} & m. & 0.1366 & 0.1308 & 0.2267 & 0.1701 & 0.1226 & \textbf{0.0984} \\ & sd. & 0.0127 & 0.0143 & 0.0179 & 0.0302 & 0.0102 & 0.0144 \\
    \end{tabular}
  \end{center}
\end{table}

\begin{table}[t]
  \begin{minipage}[t]{.45\textwidth}
    \caption{Quadratic errors of \texttt{SuperLearner} and \cobra. Means and standard deviations over
      $100$ independent replications.}
    \label{tabSL}
    \begin{center}
      \begin{tabular}{cccc}
        \textbf{Uncorr.} & & \texttt{SL} & \cobra \\ \hline\hline
        \multirow{2}{*}{\autoref{m1}} & m. & 0.0541 & \textbf{0.0320} \\ & sd. & 0.0053 & 0.0104 \\
        \multirow{2}{*}{\autoref{m2}} & m. & \textbf{0.1765} & 0.3569
        \\ & sd. & 0.0167 & 0.8797 \\
        \multirow{2}{*}{\autoref{m3}} & m. &
        \textbf{0.2081} & 0.2573
        \\ & sd. & 0.0282 & 0.0699 \\
        \multirow{2}{*}{\autoref{m4}} & m. & 4.3114 & \textbf{3.7464}
        \\ & sd. & 0.4138 & 0.8746 \\
        \multirow{2}{*}{\autoref{m5}} & m. & \textbf{0.2119} & 0.2187
        \\ & sd.  & 0.0317 & 0.0427\\
        \multirow{2}{*}{\autoref{m6}} & m. & \textbf{0.7627} & 1.0220
        \\ & sd. & 0.1023 & 0.3347 \\
        \multirow{2}{*}{\autoref{m7}} & m. & \textbf{0.1705} & 0.3103 \\ & sd. & 0.0260 & 0.0490  \\        
        \multirow{2}{*}{\autoref{m8}} & m. & 0.1081 & \textbf{0.1075} \\ & sd. & 0.0121 & 0.0235 \\
        \\
        \textbf{Corr.} & & \texttt{SL} & \cobra \\ \hline\hline
        \multirow{2}{*}{\autoref{m1}} & m. & 0.8733 & \textbf{0.3262}
        \\ & sd.& 0.2740 & 0.1242 \\
        \multirow{2}{*}{\autoref{m2}} & m. & 2.3391 & \textbf{1.3984}
        \\ & sd. & 0.4958 & 0.3804 \\
        \multirow{2}{*}{\autoref{m3}} & m. & \textbf{3.1885} & 3.3201 \\ & sd. & 1.5101 & 1.8056 \\
        \multirow{2}{*}{\autoref{m4}} & m. & 25.1073 & \textbf{9.3964} \\ & sd. & 7.3179 & 2.8953 \\
        \multirow{2}{*}{\autoref{m5}} & m. & 5.6478 & \textbf{4.9990} \\ & sd. & 7.7271 & 9.3103 \\
        \multirow{2}{*}{\autoref{m6}} & m. & \textbf{0.8967} & 1.1988 \\ & sd. & 0.1197 & 0.4573 \\
        \multirow{2}{*}{\autoref{m7}} & m. & \textbf{3.0367} & 3.1401 \\ & sd. & 1.6225 & 1.6097 \\
        \multirow{2}{*}{\autoref{m8}} & m. & 0.1116 & \textbf{0.1045} \\ & sd. & 0.0111 & 0.0216 \\
        \end{tabular}
    \end{center}
  \end{minipage}
  \hfill
  \begin{minipage}[t]{.45\textwidth}
    \caption{Average CPU-times in seconds. No parallelization. Means and
    standard deviations over $10$ independent replications.}
    \label{tabTime}
    \begin{center}
      \begin{tabular}{cccc}
        \textbf{Uncorr.} & & \texttt{SL} & \cobra \\ \hline\hline
        \multirow{2}{*}{\autoref{m1}} & m. & 53.92 & \textbf{10.92} \\ & sd. & 1.42 & 0.29 \\
        \multirow{2}{*}{\autoref{m2}} & m. & 57.96 & \textbf{11.90} \\ & sd. & 0.95 & 0.31 \\
        \multirow{2}{*}{\autoref{m3}} & m. & 53.70 & \textbf{10.66} \\ & sd. & 0.55 & 0.11 \\
        \multirow{2}{*}{\autoref{m4}} & m. & 55.00 & \textbf{11.15} \\ & sd. & 0.74 & 0.18 \\
        \multirow{2}{*}{\autoref{m5}} & m. & 28.46 & \textbf{5.01} \\ & sd. & 0.73 & 0.06 \\
        \multirow{2}{*}{\autoref{m6}} & m. & 22.97 & \textbf{3.99} \\ & sd. & 0.27 & 0.05 \\
        \multirow{2}{*}{\autoref{m7}} & m. & 127.80 & \textbf{35.67} \\ & sd. & 5.69 & 1.91 \\
        \multirow{2}{*}{\autoref{m8}} & m. & 32.98 & \textbf{6.46} \\ & sd. & 1.33 & 0.33 \\
        \\
        % \hline
        \textbf{Corr.} & & \texttt{SL} & \cobra \\ \hline\hline
        \multirow{2}{*}{\autoref{m1}} & m. & 61.92 & \textbf{11.96} \\ & sd. & 1.85 & 0.27 \\
        \multirow{2}{*}{\autoref{m2}} & m. & 70.90 & \textbf{14.16} \\ & sd. & 2.47 & 0.57 \\
        \multirow{2}{*}{\autoref{m3}} & m. & 59.91 & \textbf{11.92} \\ & sd. & 2.06 & 0.41 \\
        \multirow{2}{*}{\autoref{m4}} & m. & 63.58 & \textbf{13.11} \\ & sd. & 1.21 & 0.34 \\
        \multirow{2}{*}{\autoref{m5}} & m. & 31.24 & \textbf{5.02} \\ & sd. & 0.86 & 0.07 \\
        \multirow{2}{*}{\autoref{m6}} & m. & 24.29 & \textbf{4.12} \\ & sd. & 0.82 & 0.15 \\
        \multirow{2}{*}{\autoref{m7}} & m. & 145.18 & \textbf{41.28} \\ & sd. & 8.97 & 2.84 \\
        \multirow{2}{*}{\autoref{m8}} & m. & 31.31 & \textbf{6.24} \\ & sd. & 0.73 & 0.11 \\
        % \hline
      \end{tabular}
    \end{center}
  \end{minipage}
\end{table}

\FloatBarrier

\bibliographystyle{elsarticle-harv}
\bibliography{biblio}

\clearpage

\renewcommand{\thesection}{\Alph{section}}

\begin{center}
\textbf{\large \textsc{Supplementary Material} \\ COBRA: A
  Combined Regression Strategy \\ by G. Biau, A. Fischer, B. Guedj
  and J. D. Malley}
\end{center}

\setcounter{equation}{0}
\setcounter{section}{0}
\setcounter{page}{1}

\makeatletter 
\renewcommand{\thefigure}{\@arabic\c@figure\ (SM)}
\makeatother

\makeatletter 
\renewcommand{\thetable}{\@arabic\c@table\ (SM)}
\makeatother

\section{Proofs}\label{section:proof}

\subsection{Proof of \autoref{theo:orac}}

We have
\begin{align*}
  \mathbb E|T_n(\br_k(\bX))-r^\star(\bX)|^2 &= \mathbb
  E|T_n(\br_k(\bX))-T(\br_k(\bX))|^2\\ & \quad +\mathbb E|T(\br_k(\bX))-r^\star(\bX)|^2\\
& \quad -2 \mathbb E[(T_n(\br_k(\bX))-T(\br_k(\bX)))(T(\br_k(\bX))-r^\star(\bX))].
\end{align*}
As for the double product, notice that
\begin{align*}
&\mathbb E[(T_n(\br_k(\bX))-T(\br_k(\bX)))(T(\br_k(\bX))-r^\star(\bX))] \\
& \quad =\mathbb E \left[ \mathbb E\left[(T_n(\br_k(\bX))-T(\br_k(\bX)))(T(\br_k(\bX))-r^\star(\bX)) | \br_k(\bX), \mathcal D_n\right]\right]\\
&  \quad =\mathbb E \left[(T_n(\br_k(\bX))-T(\br_k(\bX))) \mathbb E\left[T(\br_k(\bX))-r^\star(\bX) | \br_k(\bX), \mathcal D_n\right]\right].
\end{align*}
But 
\begin{align*}
\mathbb E [r^\star(\bX) | \br_k(\bX), \mathcal D_n] &= \mathbb E [r^\star(\bX) | \br_k(\bX)]\\
& \quad (\mbox{by independence of $\bX$ and $\mathcal D_n$})\\
&=\mathbb E [\mathbb E[Y|\bX]|  \br_k(\bX)]\\
&=\mathbb E[Y|\br_k(\bX)]\\
& \quad (\mbox{since } \sigma(\br_k(\bX)) \subset  \sigma(\bX))\\
& =T(\br_k(\bX)).
\end{align*}
Consequently,
$$\mathbb E[(T_n(\br_k(\bX))-T(\br_k(\bX)))(T(\br_k(\bX))-r^\star(\bX))]=0$$
and 
$$\mathbb E|T_n(\br_k(\bX))-r^\star(\bX)|^2 = \mathbb E|T_n(\br_k(\bX))-T(\br_k(\bX))|^2+\mathbb E|T(\br_k(\bX))-r^\star(\bX)|^2.$$
Thus, by definition of the conditional expectation, and using the fact that $T(\br_k(\bX))=\mathbb E [r^\star(\bX) | \br_k(\bX)]$,
$$
\mathbb E|T_n(\br_k(\bX))-r^\star(\bX)|^2 \leq \mathbb E|T_n(\br_k(\bX))-T(\br_k(\bX))|^2+\inf_{f}\mathbb E|f(\br_k(\bX))-r^\star(\bX)|^2,
$$
where the infimum is taken over all square integrable functions of $\br_k(\bX)$. In particular,
\begin{align*}
&\mathbb E|T_n(\br_k(\bX))-r^\star(\bX)|^2\\
& \quad \leq \min_{m=1,\dots,M}\mathbb E|r_{k,m}(\bX)-r^\star(\bX)|^2 +\mathbb E|T_n(\br_k(\bX))-T(\br_k(\bX))|^2,
\end{align*}
as desired.

\subsection{Proof of \autoref{prop:cons}}

Note that the second statement is an immediate consequence of the
first statement and \autoref{theo:orac}, therefore we only have to
prove that $$\mathbb E\left |T_n\left (\br_{k}(\bX)\right)-T\left(\br_{k}(\bX)\right)\right|^2 \to 0 \quad \mbox{as } \ell \to \infty.$$ 
We start with a technical lemma, whose proof can be found in the monograph by \citet{GKKW2002}.
\begin{lem}\label{lem:binom}Let $B(n,p)$ be a  binomial random variable with parameters $n\geq 1$ and $p>0$.
Then
$$\mathbb E\left[\frac 1{1+B(n,p)}\right]\leq \frac 1 {p(n+1)}$$
and
$$\mathbb E\left[\frac {\1_{\{B(n,p)>0\}}}{B(n,p)}\right]\leq \frac 2 {p(n+1)}.$$
\end{lem}
\medskip

For all distributions of $(\bX,Y)$, using the elementary inequality $(a+b+c)^2\leq 3(a^2+b^2+c^2)$, note that
\begin{align}
&\mathbb E|T_n(\br_k(\bX))-T(\br_k(\bX))|^2\nonumber
\\&\quad=
\mathbb E\left|\sum_{i=1}^\ell W_{n,i}(\bX) \left(Y_i-T(\br_k(\bX_i))+T(\br_k(\bX_i))-T(\br_k(\bX))+T(\br_k(\bX))\right)\right.\nonumber\\&\qquad\Bigg.-T(\br_k(\bX))\Bigg|^2\nonumber
\\&\quad\leq 3\mathbb E\left|\sum_{i=1}^\ell W_{n,i}(\bX) (T(\br_k(\bX_i))-T(\br_k(\bX)))\right|^2 \label{eq:An}
\\&\qquad+3\mathbb E\left|\sum_{i=1}^\ell W_{n,i}(\bX) (Y_i-T(\br_k(\bX_i)))\right|^2\label{eq:Bn}
\\&\qquad+3\mathbb E\left|\left(\sum_{i=1}^\ell W_{n,i}(\bX)-1\right) T(\br_k(\bX))\right|^2.\label{eq:Cn}
\end{align}

Consequently, to prove the proposition, it suffices to establish that \eqref{eq:An}, \eqref{eq:Bn} and \eqref{eq:Cn} tend to $0$ as $\ell$ tends to infinity. This is done, respectively, in \autoref{prop:A}, \autoref{prop:B} and \autoref{prop:C} below.

\begin{pro}\label{prop:A} Under the assumptions of \autoref{prop:cons}, 
$$\lim_{\ell \to \infty} \mathbb E\left|\sum_{i=1}^\ell W_{n,i}(\bX) (T({\br}_k(\bX_i))-T({\br}_k(\bX)))\right|^2=0.$$
\end{pro}
\begin{proof}[Proof of \autoref{prop:A}]
By the Cauchy-Schwarz inequality,  \begin{align*}
&\mathbb E\left|\sum_{i=1}^\ell W_{n,i}(\bX) (T({\br}_k(\bX_i))-T({\br}_k(\bX)))\right|^2\\
&\quad = \mathbb E\left|\sum_{i=1}^\ell\sqrt{W_{n,i}(\bX)}\sqrt{W_{n,i}(\bX)}\left( 
T({\br}_k(\bX_i))-T({\br}_k(\bX))\right)\right|^2\\
&\quad\leq \mathbb E\left[\sum_{j=1}^\ell W_{n,j}(\bX)\sum_{i=1}^\ell W_{n,i}(\bX)\left| 
T({\br}_k(\bX_i))-T({\br}_k(\bX))\right|^2\right]\\
&\quad=
\mathbb E\left[\sum_{i=1}^\ell W_{n,i}(\bX)\left| 
T({\br}_k(\bX_i))-T({\br}_k(\bX))\right|^2\right]\\&\quad:=A_n.
\end{align*}

The function $T$ is such that $\mathbb
E[T^2({\br}_k(\bX))]<\infty$. Therefore, it can be approximated in an
$L^2$ sense by a
continuous function with compact support, say $\tilde{T}$ \citep[see, \emph{e.g.}, Theorem
A.1 in][]{GKKW2002}. More precisely, for any $\eta>0$, there exists a function $\tilde{T}$ such that 
$$\mathbb E \left|T({\br}_{k}(\bX))-\tilde{T}({\br}_{k}(\bX))\right|^2<\eta.$$
Consequently, we obtain
\begin{align*}
A_n&=\mathbb E \left[\sum_{i=1}^\ell W_{n,i}(\bX)|T({\br}_{k}(\bX_i))-T({\br}_{k}(\bX))|^2\right]
\\
&\leq 3 \mathbb E \left[\sum_{i=1}^\ell W_{n,i}(\bX)|T({\br}_{k}(\bX_i))-\tilde{T}({\br}_{k}(\bX_i))|^2\right]\\&\quad+3
 \mathbb E \left[\sum_{i=1}^\ell W_{n,i}(\bX)|\tilde{T}({\br}_{k}(\bX_i))-\tilde{T}({\br}_{k}(\bX))|^2\right]\\&\quad+3
 \mathbb E \left[\sum_{i=1}^\ell W_{n,i}(\bX)|\tilde{T}({\br}_{k}(\bX))-T({\br}_{k}(\bX))|^2\right]\\&:=3A_{n1}+3A_{n2}+3A_{n3}.
\end{align*}
\paragraph{Computation of $A_{n3}$}
Thanks to the approximation of $T$ by $\tilde{T}$,
 \begin{align*}
A_{n3}&=\mathbb E \left[\sum_{i=1}^\ell W_{n,i}(\bX)|T({\br}_{k}(\bX))-\tilde{T}({\br}_{k}(\bX))|^2\right]\\&\leq
\mathbb E \left|T({\br}_{k}(\bX))-\tilde{T}({\br}_{k}(\bX))\right|^2<\eta.
 \end{align*} 

\paragraph{Computation of $A_{n1}$}
Denote by $\mu$ the distribution of $\bX$. Then,
\begin{align*}
A_{n1} &=\mathbb E \left[\sum_{i=1}^\ell W_{n,i}(\bX)|\tilde{T}({\br}_{k}(\bX_i))-T({\br}_{k}(\bX_i))|^2\right]
\\&=\ell \mathbb E \left[\frac{{\1}_{\bigcap_{m=1}^M\{|r_{k,m}(\bX)-r_{k,m}(\bX_1)|\leq \e_\ell\}}}{\sum_{j=1}^\ell {\1}_{\bigcap_{m=1}^M\{|r_{k,m}(\bX)-r_{k,m}(\bX_j)|\leq \e_\ell\}}}|\tilde{T}({\br}_{k}(\bX_1))-T({\br}_{k}(\bX_1))|^2\right].
\\
&=\ell \E\Bigg\{\int_{\R^d} |\tilde{T}({\br}_{k}(\bu))-T({\br}_{k}(\bu))|^2 \\
\times\mathbb E
&\Bigg[\int_{\R^d}\frac{{\1}_{\bigcap_{m=1}^M\{|r_{k,m}(\bx)-r_{k,m}(\bu)|\leq
      \e_\ell\}}}{{\1}_{\bigcap_{m=1}^M\{|r_{k,m}(\bx)-r_{k,m}(\bu)|\leq
      \e_\ell\}}+\sum_{j=2}^\ell
    {\1}_{\bigcap_{m=1}^M\{|r_{k,m}(\bx)-r_{k,m}(\bX_j)|\leq
      \e_\ell\}}}\mu(\textrm{d}\bx) \\ &\qquad\qquad \Bigg| \mathcal{D}_k\Bigg] \mu(\textrm{d}\bu)\Bigg\}.
\end{align*}
Letting
\begin{align*}
A_{n1}'&=\mathbb E
\Bigg[\int_{\R^d}\frac{{\1}_{\bigcap_{m=1}^M\{|r_{k,m}(\bx)-r_{k,m}(\bu)|\leq
      \e_\ell\}}}{{\1}_{\bigcap_{m=1}^M\{|r_{k,m}(\bx)-r_{k,m}(\bu)|\leq
      \e_\ell\}}+\sum_{j=2}^\ell
    {\1}_{\bigcap_{m=1}^M\{|r_{k,m}(\bx)-r_{k,m}(\bX_j)|\leq
      \e_\ell\}}}\mu(\textrm{d}\bx) \\ &\qquad\qquad\Bigg| \mathcal{D}_k\Bigg],
\end{align*}
let us prove that $A_{n1}'\leq \frac{2^M}{\ell}$.
To this aim, observe that
\begin{align*}
A_{n1}'&= \mathbb E \left[\int_{\R^d}\frac{{\1}_{\{\bx\in\bigcap_{m=1}^Mr_{k,m}^{-1}([r_{k,m}(\bu)-\e_\ell,r_{k,m}(\bu)+\e_\ell])\}}}{1+\sum_{j=2}^\ell {\1}_{\{\bX_j\in\bigcap_{m=1}^Mr_{k,m}^{-1}([r_{k,m}(\bx)-\e_\ell,r_{k,m}(\bx)+\e_\ell])\}}}\mu(\textrm{d}\bx)\Bigg| \mathcal{D}_k\right]\\&=\mathbb E \left[\int_{\R^d}\frac{{\1}_{\{\bx\in\bigcup_{(a_1,\dots,a_M)\in \{1,2\}^M}r_{k,1}^{-1}(I_{n,1}^{a_1}(\bu))\cap\dots\cap r_{k,M}^{-1}(I_{n,M}^{a_M}(\bu))\}}}{1+\sum_{j=2}^\ell {\1}_{\{\bX_j\in\bigcap_{m=1}^Mr_{k,m}^{-1}([r_{k,m}(\bx)-\e_\ell,r_{k,m}(\bx)+\e_\ell])\}}}\mu(\textrm{d}\bx)\Bigg| \mathcal{D}_k\right]\\
&\leq\sum_{p=1}^{2^M}\mathbb E \left[\int_{\R^d}\frac{{\1}_{\{\bx\in R^p_n(\bu)\}}}{1+\sum_{j=2}^\ell {\1}_{\{\bX_j\in\bigcap_{m=1}^Mr_{k,m}^{-1}([r_{k,m}(\bx)-\e_\ell,r_{k,m}(\bx)+\e_\ell])\}}}\mu(\textrm{d}\bx)\Bigg| \mathcal{D}_k\right].
\end{align*}
Here, $I^1_{n,m}(\bu)=[r_{k,m}(\bu)-\e_\ell,r_{k,m}(\bu)]$,   $I^2_{n,m}(\bu)=[r_{k,m}(\bu),r_{k,m}(\bu)+\e_\ell]$, and $R_n^p(\bu)$ is the $p$-th set of the form $r_{k,1}^{-1}(I_{n,1}^{a_1}(\bu))\cap\dots\cap r_{k,M}^{-1}(I_{n,M}^{a_M}(\bu))$ assuming that they have been ordered using the lexicographic order of $(a_1,\dots,a_M)$.
\medskip

Next, note that \begin{equation*}\bx\in
  R_n^p(\bu) \Rightarrow R_n^p(\bu)\subset \bigcap_{m=1}^Mr_{k,m}^{-1}([r_{k,m}(\bx)-\e_\ell,r_{k,m}(\bx)+\e_\ell]).\end{equation*}
To see this, just observe that, for all $m=1,\dots, M$, if $r_{k,m}(\bz)\in[r_{k,m}(\bu)-\e_\ell,r_{k,m}(\bu)]$, \emph{i.e.}, $r_{k,m}(\bu)-\e_\ell\leq r_{k,m}(\bz)\leq r_{k,m}(\bu)$, then, as  $r_{k,m}(\bu)-\e_\ell\leq r_{k,m}(\bx)\leq r_{k,m}(\bu)$, one has $r_{k,m}(\bx)-\e_\ell\leq r_{k,m}(\bz)\leq r_{k,m}(\bx)+\e_\ell$.
Similarly, if $r_{k,m}(\bu)\leq r_{k,m}(\bz)\leq r_{k,m}(\bu)+\e_\ell$, then  $r_{k,m}(\bu)\leq r_{k,m}(\bx)\leq r_{k,m}(\bu)+\e_\ell$ implies $r_{k,m}(\bx)-\e_\ell\leq r_{k,m}(\bz)\leq r_{k,m}(\bx)+\e_\ell$. Consequently,
\begin{align*}
A_{n1}'&\leq
\sum_{p=1}^{2^M}\mathbb E \left[\int_{\R^d}\frac{{\1}_{\{\bx\in R^p_n(\bu)\}}}{1+\sum_{j=2}^\ell {\1}_{\{\bX_j\in R^p_n(\bu)\}}}\mu(\textrm{d}\bx)\Bigg| \mathcal{D}_k\right]
\\&=\sum_{p=1}^{2^M}\mathbb E\left[\frac{\mu{\{R^p_n(\bu)\}}}{1+\sum_{j=2}^\ell {\1}_{\{\bX_j\in R^p_n(\bu)\}}}\Bigg|\mathcal{D}_k\right]\\&\leq
\sum_{p=1}^{2^M}\mathbb E \left[\frac{\mu{\{R^p_n(\bu)\}}}{\ell \mu{\{R^p_n(\bu)\}}}\Bigg| \mathcal{D}_k\right]\
\\&\leq\frac{2^M}{\ell}
\end{align*}
(by the first statement of \autoref{lem:binom}).
Thus, returning to  $A_{n1}$, we obtain $$A_{n1}\leq 2^M\mathbb E \left|\tilde{T}({\br}_{k}(\bX)-T({\br}_{k}(\bX)))\right|^2<2^M\eta.$$
\paragraph{Computation of $A_{n2}$}
For any $\delta>0$, write
\begin{align}
A_{n2}&=\mathbb E\left[ \sum_{i=1}^\ell W_{n,i}(\bX)|\tilde{T}(\br_{k}(\bX_i))-\tilde{T}(\br_{k}(\bX))|^2\right]\nonumber\\
&=\mathbb E\left[ \sum_{i=1}^\ell
  W_{n,i}(\bX)|\tilde{T}(\br_{k}(\bX_i))-\tilde{T}(\br_{k}(\bX))|^2{\1}_{\bigcup_{m=1}^M\{|r_{k,m}(\bX)-r_{k,m}(\bX_i)|>\delta
    \}}\right]\nonumber\\&\quad +\mathbb E\left[ \sum_{i=1}^\ell
  W_{n,i}(\bX)|\tilde{T}(\br_{k}(\bX_i))-\tilde{T}(\br_{k}(\bX))|^2{\1}_{\bigcap_{m=1}^M\{|r_{k,m}(\bX)-r_{k,m}(\bX_i)|\leq
    \delta\}}\right]\nonumber\end{align}
from which we get that
\begin{align}A_{n2}&\leq 4\sup_{\bu\in\mathbb R^d}|\tilde{T}(\br_{k}(\bu))|^2\mathbb E\left[ \sum_{i=1}^\ell W_{n,i}(\bX){\1}_{\bigcup_{m=1}^M\{|r_{k,m}(\bX)-r_{k,m}(\bX_i)|> \delta\}}\right]\label{eq:delta1}\\&
\qquad +\bigg(\sup_{\bu,\bv\in \mathbb R^d,\bigcap_{m=1}^M\{|r_{k,m}(\bu)-r_{k,m}(\bv)|\leq \delta\}}|\tilde{T}(\br_{k}(\bv))-\tilde{T}(\br_{k}(\bu))|\bigg)^2.\label{eq:delta2}
\end{align}
With respect to the term $\eqref{eq:delta1}$, if $\delta>\e_\ell$, then
\begin{align*}
& \sum_{i=1}^\ell W_{n,i}(\bX){\1}_{\bigcup_{m=1}^M\{|r_{k,m}(\bX)-r_{k,m}(\bX_i)|> \delta\}}
\\&\quad=\sum_{i=1}^\ell\frac{{\1}_{\bigcap_{m=1}^M\{|r_{k,m}(\bX)-r_{k,m}(\bX_i)|\leq \e_\ell\}}{\1}_{\bigcup_{m=1}^M\{|r_{k,m}(\bX)-r_{k,m}(\bX_i)|> \delta\}}}{\sum_{j=1}^\ell {\1}_{\bigcap_{m=1}^M\{|r_{k,m}(\bX)-r_{k,m}(\bX_j)|\leq \e_\ell\}}}  \\&\quad=0.
\end{align*}
It follows that, for all $\delta>0$, this term converges to 0 as $\ell$ tends to infinity. On the other hand, letting $\delta \to 0$, we see that the  term $\eqref{eq:delta2}$ tends to 0 as well, by uniform continuity of $\tilde{T}$. Hence, $A_{n2}$ tends to 0 as $\ell$ tends to infinity.
%the term \eqref{eq:An}
Letting finally $\eta$ go to 0, we conclude that $A_n$ vanishes as $\ell$ tends to infinity.
\end{proof}

\begin{pro}\label{prop:B}Under the assumptions of \autoref{prop:cons}, 
$$\lim_{\ell \to \infty} \mathbb E\left|\sum_{i=1}^\ell W_{n,i}(\bX) (Y_i-T(\br_k(\bX_i)))\right|^2=0.$$
\end{pro}
\begin{proof}[Proof of \autoref{prop:B}]%Evaluation de $A_n$
\begin{align*}
&\mathbb E\left|\sum_{i=1}^\ell W_{n,i}(\bX) (Y_i-T(\br_k(\bX_i)))\right|^2\\
&\quad=\sum_{i=1}^\ell\sum_{j=1}^\ell\mathbb E [W_{n,i}(\bX)W_{n,j}(\bX) (Y_i-T(\br_k(\bX_i)))(Y_j-T(\br_k(\bX_j)))]
\\&\quad=\mathbb E\left[\sum_{i=1}^\ell W_{n,i}^2(\bX) |Y_i-T(\br_k(\bX_i))|^2\right]\\&\quad=\mathbb E\left[\sum_{i=1}^\ell W_{n,i}^2(\bX) \sigma^2(\br_k(\bX_i))\right],
\end{align*}where 
$$\sigma^2(\br_k(\bx))=\mathbb E[|Y-T(\br_k(\bX))|^2|\br_k(\bx)].$$
\medskip

For any $\eta>0$, $\sigma^2$ can be approximated in an $L^1$ sense by a continuous function with compact support $\tilde{\sigma}^2$, \emph{i.e.}, $$\mathbb E|\tilde{\sigma}^2(\br_k(\bX))-\sigma^2(\br_k(\bX))|<\eta.$$
Thus
\begin{align*}
&\mathbb E\left[\sum_{i=1}^\ell W_{n,i}^2(\bX)
  \sigma^2(\br_k(\bX_i))\right]\\&\quad\leq \mathbb
E\left[\sum_{i=1}^\ell
  W^2_{n,i}(\bX)\tilde{\sigma}^2(\br_k(\bX_i))\right] \\ &\qquad+\mathbb E\left[\sum_{i=1}^\ell W^2_{n,i}(\bX)|
\sigma^2(\br_k(\bX_i))-\tilde{\sigma}^2(\br_k(\bX_i))|\right]\\
&\quad\leq \sup_{\bu\in\mathbb R^d}|\tilde{\sigma}^2(\br_k(\bu))|\mathbb E\left[\sum_{i=1}^\ell W_{n,i}^2(\bX)\right]\\
& \qquad +\mathbb E\left[\sum_{i=1}^\ell W_{n,i}(\bX)|
\sigma^2(\br_k(\bX_i))-\tilde{\sigma}^2(\br_k(\bX_i))|\right].
\end{align*}
With the same argument as for $A_{n1}$, we obtain
$$\mathbb E\left[\sum_{i=1}^\ell W_{n,i}(\bX)|
\sigma^2(\br_k(\bX_i))-\tilde{\sigma}^2(\br_k(\bX_i))|\right]\leq 2^M\eta.$$
Therefore, it remains to prove that $\mathbb E\left[\sum_{i=1}^\ell W_{n,i}^2(\bX)\right] \to 0$ as $\ell \to \infty$. 
To this aim, fix $\delta>0$, and note that
\begin{align*}
\sum_{i=1}^\ell W_{n,i}^2(\bX)
%&=\frac{\sum_{i=1}^\ell{\1}_{\bigcap_{m=1}^M\{|r_{k,m}(\bX)-r_{k,m}(\bX_i)|\leq \e_\ell\}}^2}{\left(\sum_{j=1}^\ell {\1}_{\bigcap_{m=1}^M\{|r_{k,m}(\bX)-r_{k,m}(\bX_j)|\leq \e_\ell\}}\right)^2}
&=
\frac{\sum_{i=1}^\ell{\1}_{\bigcap_{m=1}^M\{|r_{k,m}(\bX)-r_{k,m}(\bX_i)|\leq \e_\ell\}}}{\left(\sum_{j=1}^\ell {\1}_{\bigcap_{m=1}^M\{|r_{k,m}(\bX)-r_{k,m}(\bX_j)|\leq \e_\ell\}}\right)^2}\\
&\leq \min\left\{\delta, \frac 1{\sum_{i=1}^\ell {\1}_{\bigcap_{m=1}^M\{|r_{k,m}(\bX)-r_{k,m}(\bX_i)|\leq \e_\ell\}}}\right\}\\
&\leq \delta + \frac {{\1}_{\left\{\sum_{i=1}^\ell {\1}_{\bigcap_{m=1}^M\{|r_{k,m}(\bX)-r_{k,m}(\bX_i)|\leq \e_\ell\}}>0\right\}}}{\sum_{i=1}^\ell {\1}_{\bigcap_{m=1}^M\{|r_{k,m}(\bX)-r_{k,m}(\bX_i)|\leq \e_\ell\}}}
.\end{align*}
To complete the proof, we have to establish that the expectation of
the right-hand term tends to $0$. Denoting by $I$ a bounded interval on the real line, we have
%the following chain of inequalities is valid: 
\begin{align*}
&\mathbb E \left[\frac{{\1}_{\left\{\sum_{i=1}^\ell {\1}_{\left\{\bX_i\in\bigcap_{m=1}^Mr_{k,m}^{-1}([r_{k,m}(\bX)-\e_\ell,r_{k,m}(\bX)+\e_\ell])\right\}}>0\right\}} }{\sum_{i=1}^\ell {\1}_{\left\{\bX_i\in\bigcap_{m=1}^Mr_{k,m}^{-1}([r_{k,m}(\bX)-\e_\ell,r_{k,m}(\bX)+\e_\ell])\right\}}}\right]
\\&\leq \mathbb E \left[\frac{{\1}_{\left\{\sum_{i=1}^\ell {\1}_{\left\{\bX_i\in\bigcap_{m=1}^Mr_{k,m}^{-1}([r_{k,m}(\bX)-\e_\ell,r_{k,m}(\bX)+\e_\ell])\right\}}>0\right\}}
{\1}_{\left\{\bX\in \bigcap_{m=1}^Mr_{k,m}^{-1}(I)\right\}}}{\sum_{i=1}^\ell {\1}_{\left\{\bX_i\in\bigcap_{m=1}^Mr_{k,m}^{-1}([r_{k,m}(\bX)-\e_\ell,r_{k,m}(\bX)+\e_\ell])\right\}}}\right]\\&\quad
\quad +\mu\Big(\bigcup_{m=1}^Mr_{k,m}^{-1}(I^c)\Big)\\&
= \mathbb E \left[\mathbb E\left[
\frac{{\1}_{\left\{\sum_{i=1}^\ell {\1}_{\left\{\bX_i\in\bigcap_{m=1}^Mr_{k,m}^{-1}([r_{k,m}(\bX)-\e_\ell,r_{k,m}(\bX)+\e_\ell])\right\}}>0\right\}}
{\1}_{\left\{\bX\in \bigcap_{m=1}^Mr_{k,m}^{-1}(I)\right\}}}{\sum_{i=1}^\ell {\1}_{\left\{\bX_i\in\bigcap_{m=1}^Mr_{k,m}^{-1}([r_{k,m}(\bX)-\e_\ell,r_{k,m}(\bX)+\e_\ell])\right\}}}
\right. \right.\\
&\qquad\left. \left. \Big|\mathcal{D}_k,\bX\right]\right] +\mu\Big(\bigcup_{m=1}^Mr_{k,m}^{-1}(I^c)\Big)
\\&
\leq \frac 2{(\ell+1)}\mathbb E\left[\frac{{\1}_{\left\{\bX\in \bigcap_{m=1}^Mr_{k,m}^{-1}(I)\right\}}}{\mu(\bigcap_{m=1}^Mr_{k,m}^{-1}([r_{k,m}(\bX)-\e_\ell,r_{k,m}(\bX)+\e_\ell]))}\right]\\
& \quad +\mu\Big(\bigcup_{m=1}^Mr_{k,m}^{-1}(I^c)\Big).
\end{align*}
The last inequality arises from the second statement of
\autoref{lem:binom}. By an appropriate choice of $I$, according to the
technical statement \eqref{bound}, the second term on the right-hand side can be made as small as desired.
Regarding the first term, there exists a finite number $N_\ell$ of points  $\bz_1,\dots,\bz_{N_\ell}$ such that $$\bigcap_{m=1}^Mr_{k,m}^{-1}(I)\subset\bigcup_{(j_1,\dots,j_M)\in \{1,\dots,N_\ell\}^M}r_{k,1}^{-1}(I_{n,1}(\bz_{j_1}))\cap\dots\cap
r_{k,M}^{-1}(I_{n,M}(\bz_{j_M})), $$ where $I_{n,m}(\bz_j)=[\bz_j-\e_\ell/2,\bz_j+\e_\ell/2]$. 
Suppose, without loss of generality, that the sets $$r_{k,1}^{-1}(I_{n,1}(\bz_{j_1}))\cap\dots\cap
r_{k,M}^{-1}(I_{n,M}(\bz_{j_M}))$$ are ordered, and denote by  $R^p_n$ the $p$-th among the $N_\ell^M=(\lceil|I|/\e_\ell\rceil)^M$ sets. Here $|I|$ denotes the length of the interval $I$ and $\lceil x\rceil$ denotes the smallest integer greater than $x$.
For all $p$, 
\begin{equation*}\bx\in
  R_n^p \Rightarrow R_n^p\subset \bigcap_{m=1}^Mr_{k,m}^{-1}([r_{k,m}(\bx)-\e_\ell,r_{k,m}(\bx)+\e_\ell]).\end{equation*}
% if $\bX\in R_n^p$, this set is included in each $r_{k,m}^{-1}([r_{k,m}(\bX)-\e_\ell,r_{k,m}(\bX)+\e_\ell])$, $m=1,\dots,M$.
Indeed, if $\bv\in R_n^p$, then,  for all $m=1,\dots,M$, there exists  $j\in \{1,\dots,N_\ell\}$ such that   $r_{k,m}(\bv)\in [\bz_j-\e_\ell/2,\bz_j+\e_\ell/2]$, that is $\bz_j-\e_\ell/2\leq r_{k,m}(\bv)\leq \bz_j+\e_\ell/2$. Since we also have $\bz_j-\e_\ell/2\leq r_{k,m}(\bX)\leq \bz_j+\e_\ell/2$, we obtain
$r_{k,m}(\bX)-\e_\ell\leq r_{k,m}(\bv)\leq r_{k,m}(\bX)+\e_\ell$. In conclusion, 
\begin{align*}
&\mathbb E\left[\frac{{\1}_{\left\{\bX\in \bigcap_{m=1}^Mr_{k,m}^{-1}(I)\right\}}}{\mu(\bigcap_{m=1}^M
r_{k,m}^{-1}([r_{k,m}(\bX)-\e_\ell,r_{k,m}(\bX)+\e_\ell]))}\right]\\&
\quad \leq\sum_{p=1}^{N_\ell^M}\mathbb E\left[\frac{{\1}_{\left\{\bX\in R_{n}^{p}\right\}}}{\mu(\bigcap_{m=1}^M
r_{k,m}^{-1}([r_{k,m}(\bX)-\e_\ell,r_{k,m}(\bX)+\e_\ell]))}\right]
\\&
\quad \leq\sum_{p=1}^{N_\ell^M}\mathbb E\left[\frac{{\1}_{\left\{\bX\in R_{n}^{p}\right\}}}{\mu(R_n^p)}\right]\\&
\quad =N_\ell^M\\&
\quad =\left\lceil\frac{|I|}{\e_\ell}\right\rceil^M.
\end{align*}
The result follows from the assumption $\lim_{\ell \to \infty} \ell\e_\ell^M=\infty$.
\end{proof}

\begin{pro}\label{prop:C}Under the assumptions of \autoref{prop:cons}, 
$$\lim_{\ell \to \infty}\mathbb E\left|\left(\sum_{i=1}^\ell W_{n,i}(\bX)-1\right) T(\br_k(\bX))\right|^2 =0.$$
\end{pro}

\begin{proof}[Proof of \autoref{prop:C}]%Evaluation de $C_n$

Since $|\sum_{i=1}^\ell W_{n,i}(\bX)-1|\leq 1$, one has $$\left|\left(\sum_{i=1}^\ell W_{n,i}(\bX)-1\right) T(\br_k(\bX))\right|^2\leq T^2(\br_k(\bX)).$$ Consequently, by
Lebesgue's dominated convergence theorem, to prove the proposition, it suffices to show that  $W_{n,i}(\bX)$ tends to 1 almost surely.
Now,
\begin{align*}
&\mathbb P\left(\sum_{i=1}^\ell W_{n,i}(\bX)\neq 1\right)\\&\quad=
\mathbb P\left(\sum_{i=1}^\ell{\1}_{\bigcap_{m=1}^{M}\{|r_{k,m}(\bX)-r_{k,m}(\bX_i))|\leq \e_\ell\}}=0\right)\\
&\quad=
\mathbb P\left(\sum_{i=1}^\ell{\1}_{\left\{\bX_i\in\bigcap_{m=1}^{M}r_{k,m}^{-1}\left([r_{k,m}(\bX)-\e_\ell,r_{k,m}(\bX)+\e_\ell]\right)\right\}}=0\right)\\
&\quad=\int_{\R^d} \mathbb P\left( \forall i=1,\dots,\ell,{\1}_{\left\{\bX_i\in \bigcap_{m=1}^Mr_{k,m}^{-1}\left([r_{k,m}(\bx)-\e_\ell,r_{k,m}(\bx)+\e_\ell]\right)\right\}} =0\right) \mu(\textrm{d}\bx)\\
&\quad=\int_{\R^d} \left[ 1 - \mu(\cap_{m=1}^M
  r_{k,m}^{-1}\left([r_{k,m}(\bx)-\e_\ell,r_{k,m}(\bx)+\e_\ell]\right))\right]^\ell
\mu(\textrm{d}\bx).
\end{align*}
Denote by $I$ a bounded interval. Then,
\begin{align*}
&\mathbb P\left(\sum_{i=1}^\ell W_{n,i}(\bX)\neq 1\right)\\&\quad\leq\int_{\R^d} {\exp\left( - \ell \mu(\cap_{m=1}^M r_{k,m}^{-1}\left([r_{k,m}(\bx)-\e_\ell,r_{k,m}(\bx)+\e_\ell]\right))
\right)}\\&\qquad\times{\1}_{\{\bx\in\bigcap_{m=1}^Mr_{k,m}^{-1}(I)\}}\mu(\textrm{d}\bx)+\mu\Big(\bigcup_{m=1}^Mr_{k,m}^{-1}(I^c)\Big)\\
&\quad\leq \max _\bu \bu e^{-\bu}\int_{\R^d} \frac{{\1}_{\{\bx\in\bigcap_{m=1}^Mr_{k,m}^{-1}(I)\}}}{\ell \mu(\cap_{m=1}^M r_{k,m}^{-1}\left([r_{k,m}(\bx)-\e_\ell,r_{k,m}(\bx)+\e_\ell]\right))}\mu(\textrm{d}\bx)\\&\qquad+\mu\Big(\bigcup_{m=1}^Mr_{k,m}^{-1}(I^c)\Big).
\end{align*}
Using the same arguments as in the proof of \autoref{prop:B}, the probability $\mathbb P\left(\sum_{i=1}^\ell W_{n,i}(\bX)\neq 1\right)$ is bounded by $\frac {e^{-1}}{\ell}\left\lceil{\frac {|I|}{\e_\ell}}\right\rceil^ M$. This bound vanishes as $n$ tends to infinity since, by assumption, $\lim_{\ell \to \infty}\ell \e_\ell^M=\infty$.
\end{proof}

\subsection{Proof of \autoref{pro:vitesse}}

Choose $\bx\in \R^d$. An easy
calculation yields that
\begin{align}
&\E[|T_n(\br_k(\bx))-T(\br_k(\bx))|^2\big|\br_k(\bX_1),\dots,\br_k(\bX_\ell),\mathcal{D}_k]\nonumber\\
&=\E\Bigg[\big|T_n(\br_k(\bx))-\E[T_n(\br_k(\bx))\big|\br_k(\bX_1),\dots,\br_k(\bX_\ell),\mathcal{D}_k]\big|^2
    \\ 
&\quad\Big|\br_k(\bX_1),\dots,\br_k(\bX_\ell),\mathcal{D}_k\Bigg]\nonumber%\\&\qquad
+\big|\E[T_n(\br_k(\bx))\big|\br_k(\bX_1),\dots,\br_k(\bX_\ell),\mathcal{D}_k]-T(\br_k(\bx))\big|^2\nonumber
\\ &:=E_1 + E_2.\label{eq:vitdec}
\end{align}

% Denoting the first term in the expression \eqref{eq:vitdec} by $E_1$ and the second one by $E_2$,
On the one hand, we have 
\begin{align*}E_1&=\E\Big[\big|T_n(\br_k(\bx))-\E[T_n(\br_k(\bx))|\br_k(\bX_1),
\dots,\br_k(\bX_\ell),\mathcal{D}_k]\big|^2\\ &\qquad\qquad\Big|\br_k(\bX_1),\dots,\br_k(\bX_\ell),\mathcal{D}_k\Big]
\\&=\E\left[\left|\sum_{i=1}^\ell W_{n,i}(\bx)(Y_i-\E[Y_i|\br_k(\bX_i)])\right|^2|\br_k(\bX_1),\dots,\br_k(\bX_\ell),\mathcal{D}_k\right].
\end{align*}
Developing the square and noticing that
$\E\big[Y_j|Y_i,\br_k(\bX_1),\dots,\br_k(\bX_\ell),\mathcal{D}_k\big]=\E[Y_j|\br_k(\bX_j)]$,
since  $Y_j$ is independent of $Y_i$ and of the $\bX_j$'s with $j\neq
i$, we have
\begin{align}
E_1&=\E\Bigg[\frac{\sum_{i=1}^\ell\1_{\bigcap_{m=1}^M\{|r_{k,m}(\bx)-r_{k,m}(\bX_i)|\leq \e_\ell \}}|Y_i-\E[Y_i|\br_k(\bX_i)]|^2}{\left|\sum_{i=1}^\ell\1_{\bigcap_{m=1}^M\{|r_{k,m}(\bx)-r_{k,m}(\bX_i)|\leq \e_\ell \}}\right|^2}\\&\qquad\qquad\Bigg|\br_k(\bX_1),\dots,\br_k(\bX_\ell),\mathcal{D}_k\Bigg]\nonumber\\
&=\sum_{i=1}^\ell \var(Y_i|\br_k(\bX_i))\frac{\1_{\bigcap_{m=1}^M\{|r_{k,m}(\bx)-r_{k,m}(\bX_i)|\leq \e_\ell \}}}{\left|\sum_{i=1}^\ell\1_{\bigcap_{m=1}^M\{|r_{k,m}(\bx)-r_{k,m}(\bX_i)|\leq \e_\ell \}}\right|^2}.\nonumber
\end{align}
Thus,
\begin{equation}
E_1\leq 4R^2\frac{\1_{\left\{\sum_{i=1}^\ell\1_{\bigcap_{m=1}^M\{|r_{k,m}(\bx)-r_{k,m}(\bX_i)|\leq \e_\ell \}}>0\right\}}}{\sum_{i=1}^\ell\1_{\bigcap_{m=1}^M\{|r_{k,m}(\bx)-r_{k,m}(\bX_i)|\leq \e_\ell \}}},\label{eq:vitE1}
\end{equation}
where $\var (Z)$ denotes the variance of a random variable $Z$.
On the other hand, recalling the notation $\Sigma$ introduced
in \autoref{section:cobra}, we obtain for the second term $E_2$:
\begin{align}
E_2&= \big|\E[T_n(\br_k(\bx))|\br_k(\bX_1),\dots,\br_k(\bX_\ell),\mathcal{D}_k]-T(\br_k(\bx))\big|^2\nonumber\\
&=\left|\sum_{i=1}^\ell W_{n,i}(\bx)\E[Y_i|\br_k(\bX_i)]-T(\br_k(\bx))\right|^2\1_{\{\Sigma >0\}}+T^2(\br_k(\bx))\1_{\{\Sigma =0\}}\nonumber\\
&\leq \frac{\sum_{i=1}^\ell\1_{\bigcap_{m=1}^M\{|r_{k,m}(\bx)-r_{k,m}(\bX_i)|\leq \e_\ell \}}
\left|\E[Y_i|\br_k(\bX_i)]-T(\br_k(\bx))\right|^2}{\sum_{j=1}^\ell\1_{\bigcap_{m=1}^M\{|r_{k,m}(\bx)-r_{k,m}(\bX_j)|\leq \e_\ell \}}}\1_{\{\Sigma >0\}}\\
&\qquad+T^2(\br_k(\bx))\1_{\{\Sigma =0\}}
\nonumber\\
&\qquad\textrm{(by Jensen's inequality)}
\nonumber\\
&= \frac{\sum_{i=1}^\ell\1_{\bigcap_{m=1}^M\{|r_{k,m}(\bx)-r_{k,m}(\bX_i)|\leq \e_\ell \}}
\left|T(\br_k(\bX_i))-T(\br_k(\bx))\right|^2}{\sum_{j=1}^\ell\1_{\bigcap_{m=1}^M\{|r_{k,m}(\bx)-r_{k,m}(\bX_j)|\leq \e_\ell \}}}\1_{\{\Sigma >0\}}
\\&\qquad+T^2(\br_k(\bx))\1_{\{\Sigma =0\}}\nonumber\\
&\leq L^2\e_\ell^2+T^2(\br_k(\bx))\1_{\{\Sigma =0\}}.\label{eq:vitE2}
\end{align}

Now,
$$\E|T_n(\br_k(\bX))-T(\br_k(\bX))|^2
\leq \int_{\R^d}\E|(T_n(\br_k(\bx))-T(\br_k(\bx))|^2\mu (\mathrm{d}\bx).$$
Then, using the decomposition \eqref{eq:vitdec} and the upper bounds \eqref{eq:vitE1} and \eqref{eq:vitE2}, 
\begin{align*}&\E|T_n(\br_k(\bX))-T(\br_k(\bX))|^2\\&\quad\leq 
\int_{\R^d}\E\left[\frac{4R^2\1_{\{\Sigma >0\}}}{B}\right]\mu(\mathrm{d} \bx)
+L^2\e_\ell ^2+\int_{\R^d}\E\left[T^2(\br_k(\bx))\1_{\{\Sigma =0\}}\right]\mu(\mathrm{d}\bx)
\\&\quad\leq 
\int_{\R^d}\E\left\{\E\left[\frac{4R^2\1_{\{\Sigma >0\}}}{B}\Big|\mathcal{D}_k\right]\right\}\mu(\mathrm{d} \bx)
+L^2\e_\ell ^2\\&\qquad+\int_{\R^d}\E\left\{\E\left[T^2(\br_k(\bx))\1_{\{\Sigma =0\}}|\mathcal{D}_k\right]\right\}\mu(\mathrm{d}\bx).
\end{align*}
Thus, thanks to \autoref{lem:binom}, 
\begin{align*}
&\E|T_n(\br_k(\bX))-T(\br_k(\bX))|^2\\&\quad\leq\frac{8R^2}{(\ell+1)}\int_{\R^d}
\frac{1}{\mu(\bigcap_{m=1}^M\{|r_{k,m}(\bx)-r_{k,m}(\bX)|\leq \e_\ell \})}\mu(\mathrm{d} \bx)+L^2\e_\ell ^2\\&\qquad+\int_{\R^d} T^2(\br_k(\bx)) \left(1-\mu(\bigcap_{m=1}^M\{|r_{k,m}(\bx)-r_{k,m}(\bX)|\leq \e_\ell \})\right)^\ell\mu(\mathrm{d}\bx).
\end{align*}
Consequently,
\begin{align*}
&\E|T_n(\br_k(\bX))-T(\br_k(\bX))|^2\\&\quad\leq\frac{8R^2}{(\ell+1)}\int_{\R^d}\frac{1}{\mu(\bigcap_{m=1}^M\{|r_{k,m}(\bx)-r_{k,m}(\bX)|\leq \e_\ell \})}\mu(\mathrm{d} \bx)+L^2\e_\ell ^2\\&\qquad+\int_{\R^d} T^2(\br_k(\bx)) \exp\left(-\ell\mu(\bigcap_{m=1}^M\{|r_{k,m}(\bx)-r_{k,m}(\bX)|\leq \e_\ell \})\right)\mu(\mathrm{d}\bx)
\\
&\quad\leq\frac{8R^2}{(\ell+1)}\int_{\R^d}\frac{1}{\mu(\bigcap_{m=1}^M\{|r_{k,m}(\bx)-r_{k,m}(\bX)|\leq \e_\ell \})}\mu(\mathrm{d} \bx)+L^2\e_\ell ^2\\&\qquad+\Bigg(\sup_{\bx\in \R^d}T^2(\br_k(\bx)) \max_{\mathbf u\in\R^+}\mathbf u e^{-\mathbf u}\\&\qquad\qquad\times\int_{\R^d} \frac{1}{\ell\mu(\bigcap_{m=1}^M\{|r_{k,m}(\bx)-r_{k,m}(\bX)|\leq \e_\ell \})}\mu(\mathrm{d}\bx)\Bigg).
\end{align*}

Introducing a bounded interval $I$ as in the proof of
\autoref{prop:cons}, we observe that the boundedness of the $\br_k$ yields
that $$\mu\left(\bigcup_{m=1}^Mr_{k,m}^{-1}(I^c)\right) = 0,$$
 as soon as $I$ is sufficiently large, independently of $k$.
Then, proceeding as in the proof of \autoref{prop:cons}, we obtain
\begin{align*}
&\E|T_n(\br_k(\bX))-T(\br_k(\bX))|^2\\&\quad\leq
8R^2\left\lceil\frac{|I|}{\e_\ell }\right\rceil^M\frac{1}{\ell+1}+L^2\e_\ell ^2+R^2\max_{\mathbf u\in\R^+}\mathbf ue^{-\mathbf u}\left\lceil\frac{|I|}{\e_\ell }\right\rceil^M\frac{1}{\ell}
\\&\quad\leq C_1\frac{R^2}{\ell\e_\ell ^M}+L^2\e_\ell ^2,
\end{align*}
for some positive constant $C_1$, independent of $k$.
Hence, for the choice $\e_\ell\propto \ell^{-\frac 1{M+2}}$,
we obtain $$\E|T_n(\br_k(\bX))-T(\br_k(\bX))|^2\leq C \ell^{-\frac
  2{M+2}},$$ for some positive constant $C$ depending on $L$, $R$
and independent of $k$, as desired.

\section{Numerical results}

\begin{table}[h]
  \caption{Quadratic errors of the implemented machines and \cobra in
    high-dimensional situations. Means and standard deviations over
    $200$ independent replications.}
  \label{tab-HD}
  \begin{center}
    \begin{tabular}{cccccccc}
      & & \texttt{lars} & \texttt{ridge} & \texttt{fnn} &
      \texttt{tree} & \texttt{rf} & \cobra \\ \hline\hline
      \multirow{2}{*}{\autoref{mHD1}} & m. & 1.5698 & 2.9752 & 3.9285 &
      1.8646 & 1.5001 & \textbf{0.9996} \\
      & sd. & 0.2357 & 0.4171 & 0.5356 & 0.3751 & 0.2491 & 0.1733 \\
      \multirow{2}{*}{\autoref{mHD2}} & m. & 5.2356 & 5.1748 & 6.1395 &
      6.1585 & 4.8667 & \textbf{2.7076} \\ & sd. & 0.6885 & 0.7139 & 0.9192 &
      0.9298 & 0.6634 & 0.3810 \\
      \multirow{2}{*}{\autoref{mHD3}} & m. & 0.1584 & 0.1055 & 0.1363 &
      0.0058 & 0.0327 & \textbf{0.0049} \\ & sd. & 0.0199 & 0.0119 & 0.0176 &
      0.0010 & 0.0052 & 0.0009 % \\
      % \multirow{2}{*}{Model 4} & m. & & & & & & \\ & sd. & & & & & &
    \end{tabular}
  \end{center}
\end{table}

\begin{table}[h]
  \caption{Quadratic errors of exponentially weighted aggregate
    (EWA) and \cobra. $200$ independent replications.}
  \label{tab-EWA}
  \begin{center}
    \begin{tabular}{cccc}
      & & \texttt{EWA} & \cobra \\ \hline\hline
      \multirow{2}{*}{\autoref{mHD1}} & m. & 1.1712 & \textbf{1.1360} \\
      & sd. &  0.2090 & 0.2468  \\
      \multirow{2}{*}{\autoref{mHD2}} & m. & \textbf{9.4789} & 12.4353 \\
      & sd. & 5.6275 & 9.1267  \\
      \multirow{2}{*}{\autoref{mHD3}} & m. & 0.0244 & \textbf{0.0128}  \\
      & sd. & 0.0042 & 0.0237  \\
      \multirow{2}{*}{\autoref{mSTAB}} & m. & 0.4175 &
      \textbf{0.3124} \\
      & sd. & 0.0513  & 0.0884 
    \end{tabular}
  \end{center}
\end{table}

%\FloatBarrier

\begin{figure}[t]
  \caption{Examples of calibration of parameters $\e_\ell$ and $\alpha$. The bold point is the minimum.}
  \label{epscal}
  \subfloat[\autoref{m5}, uncorrelated design.]{\includegraphics[width = .49\textwidth]{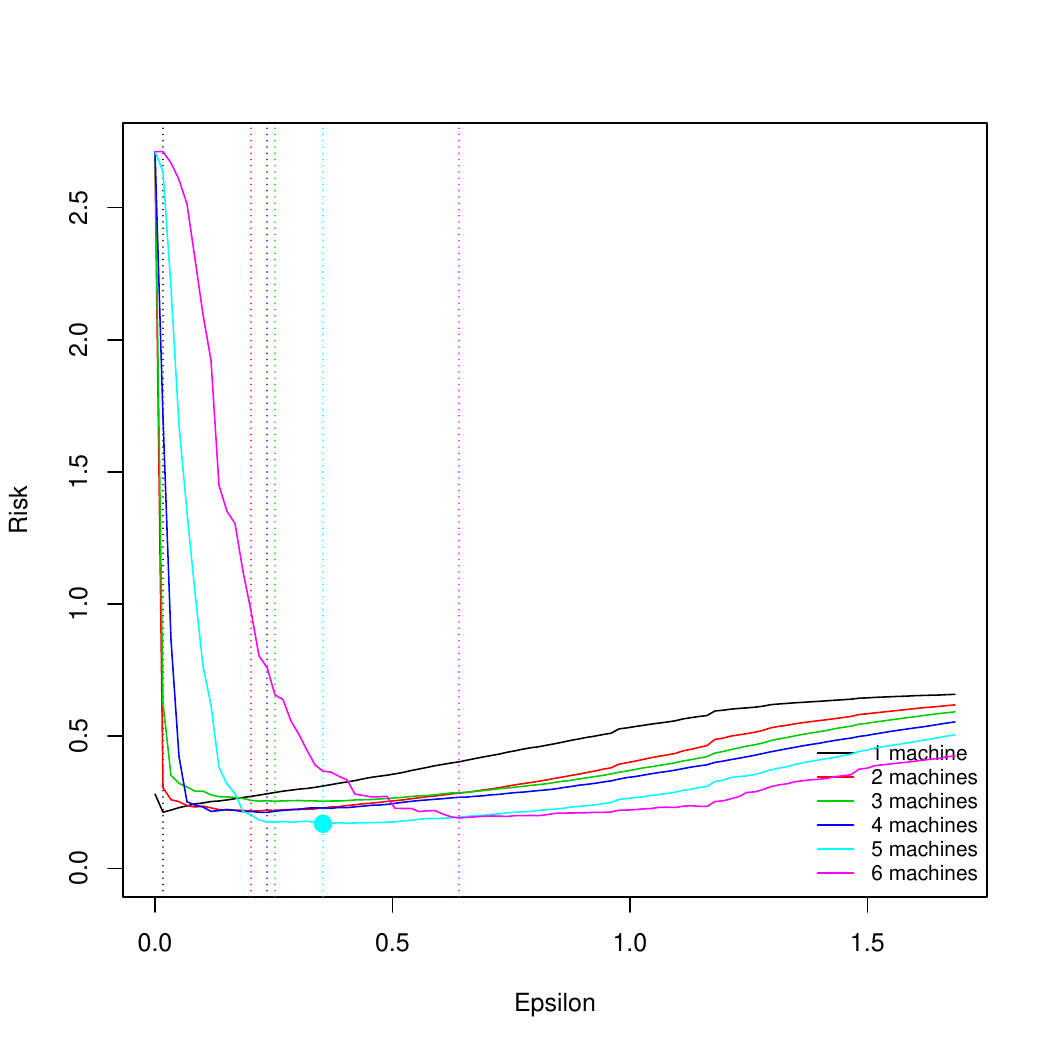}}
  \subfloat[\autoref{m5}, correlated design.]{\includegraphics[width =
    .49\textwidth]{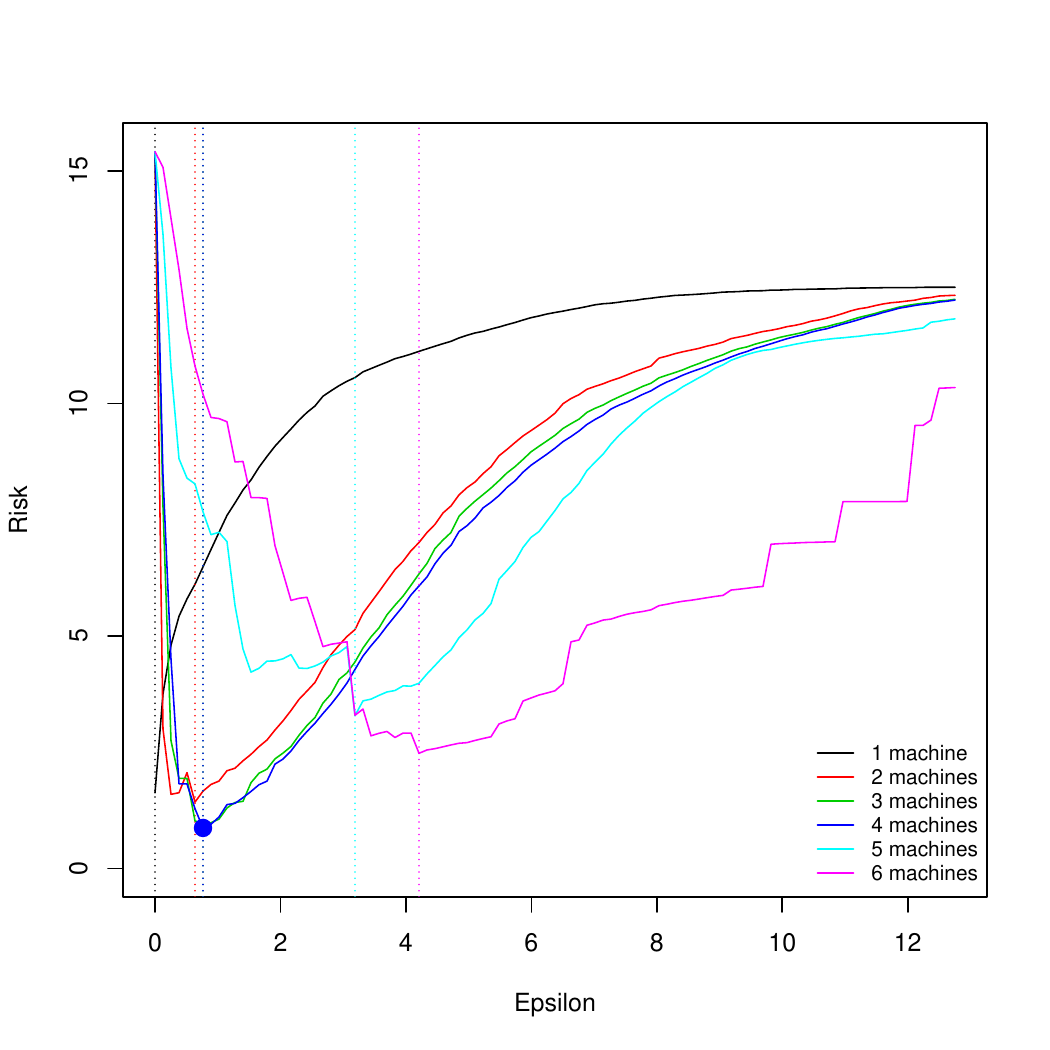}}
  \\
  \subfloat[\autoref{mHD1}.]{\includegraphics[width =
    .49\textwidth]{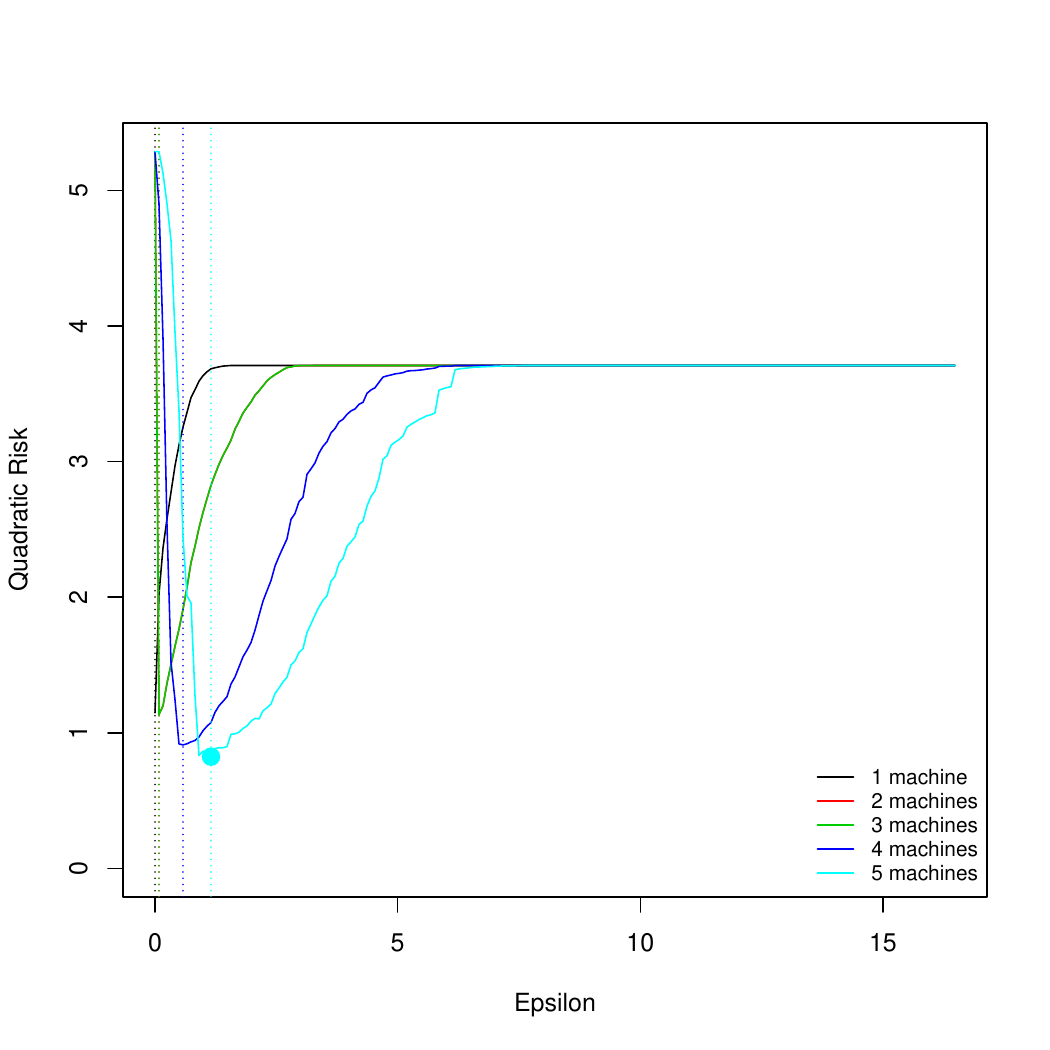}}
  \subfloat[\autoref{mSTAB}.]{\includegraphics[width = .49\textwidth]{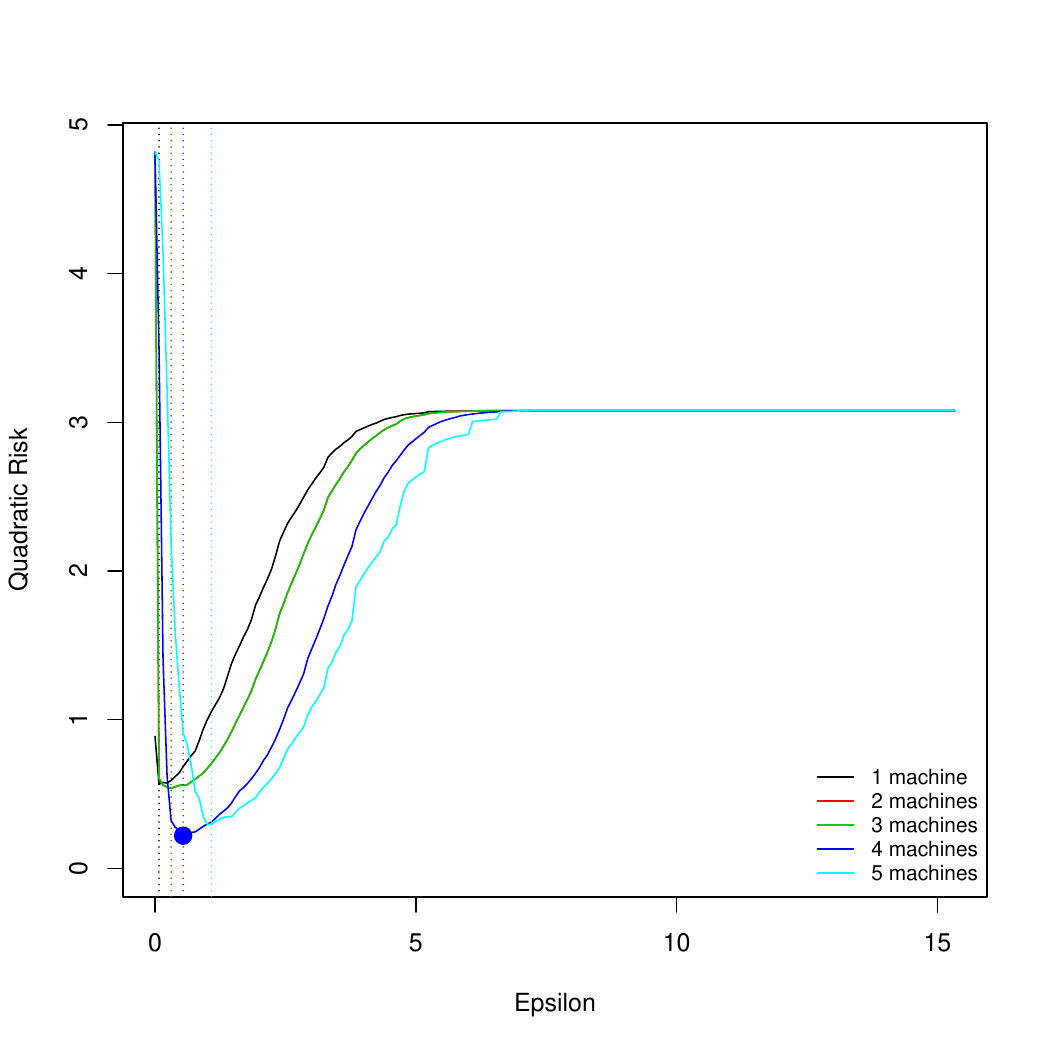}}
\end{figure}

\begin{figure}[t]
  \caption{Boxplots of quadratic errors, uncorrelated design. From
    left to right: \texttt{lars}, \texttt{ridge}, \texttt{fnn},
    \texttt{tree}, \texttt{randomForest}, \cobra.}
  \label{boxplot-U}
  \subfloat[\autoref{m1}.]{\includegraphics[width = .24\textwidth]{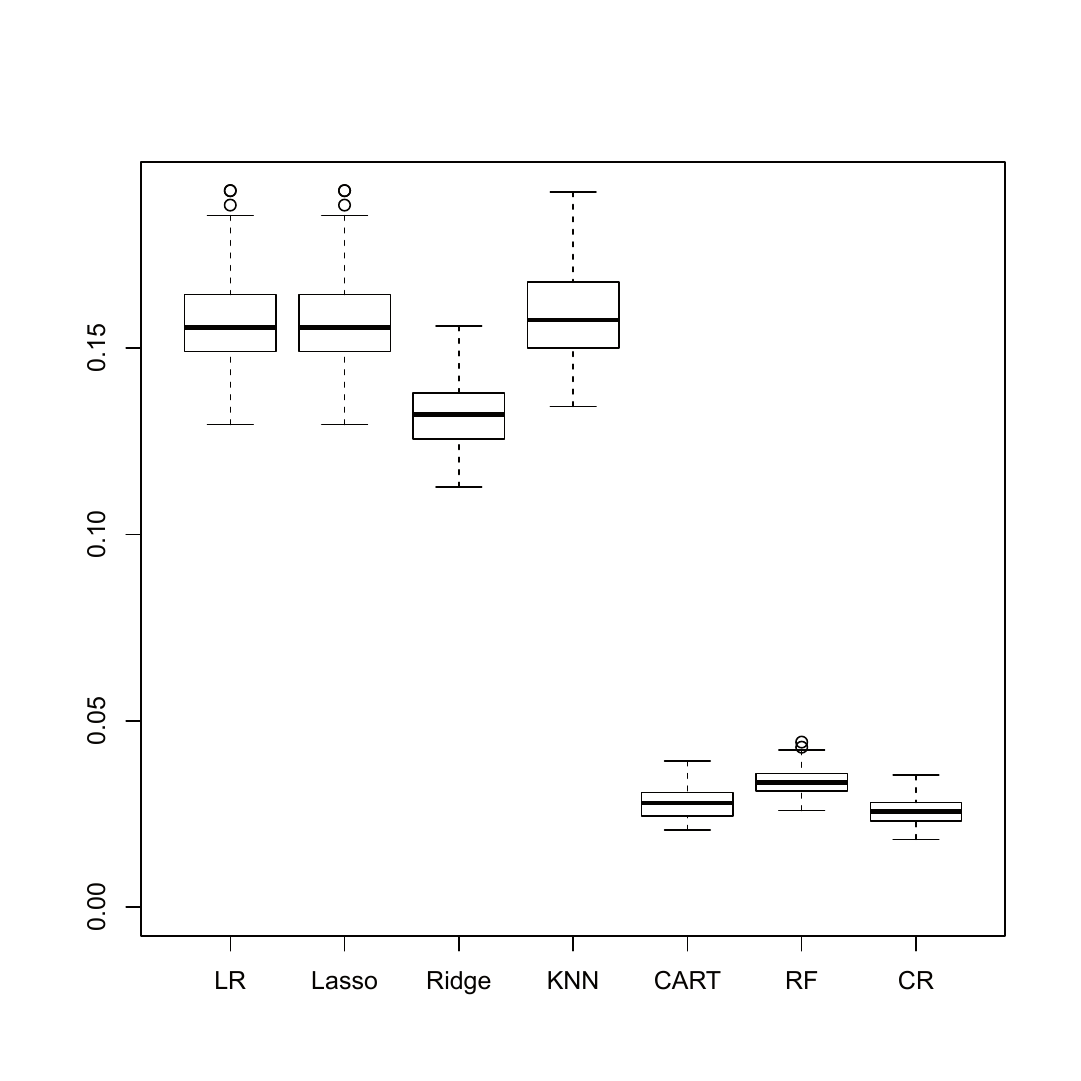}}
  \subfloat[\autoref{m2}.]{\includegraphics[width = .24\textwidth]{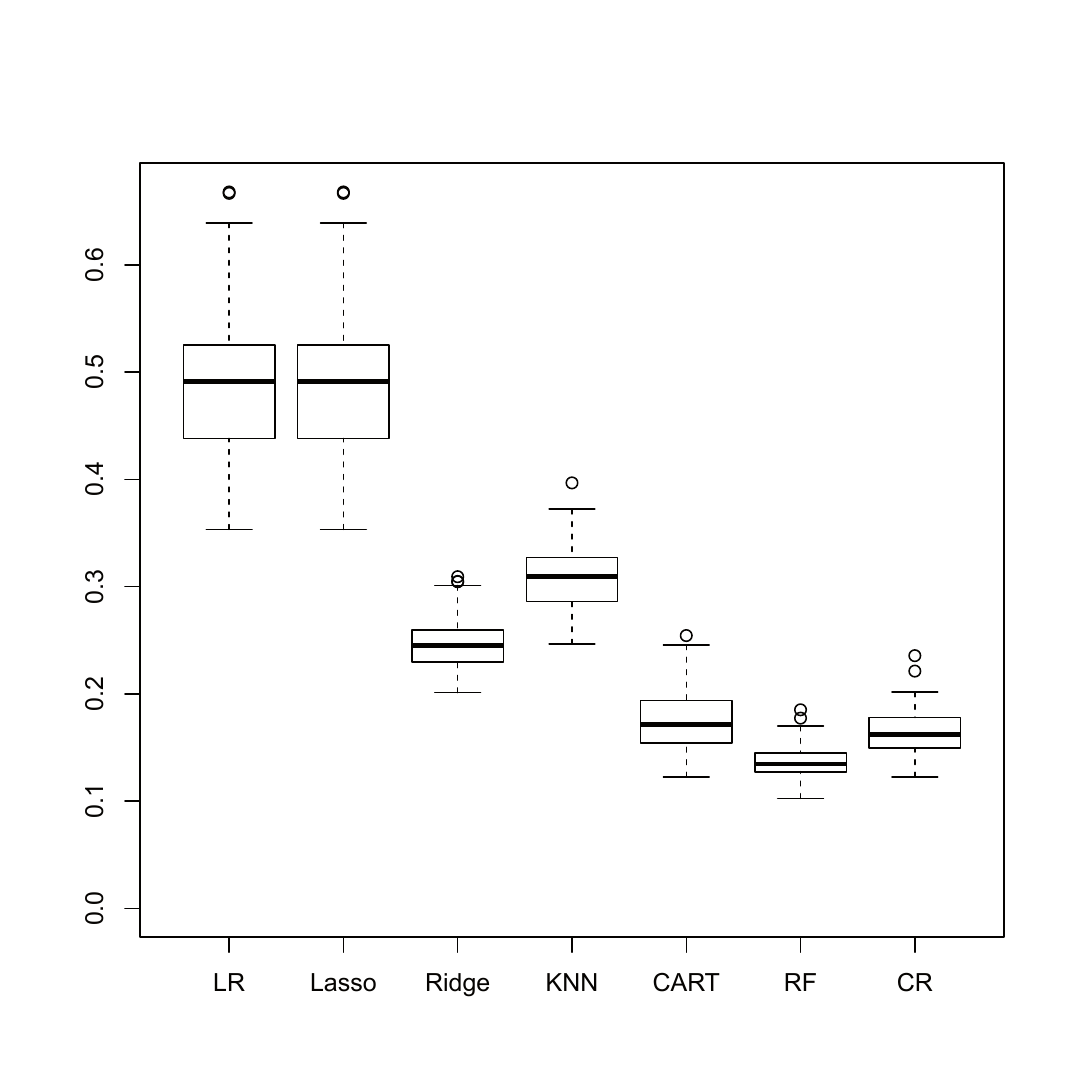}}
  \subfloat[\autoref{m3}.]{\includegraphics[width = .24\textwidth]{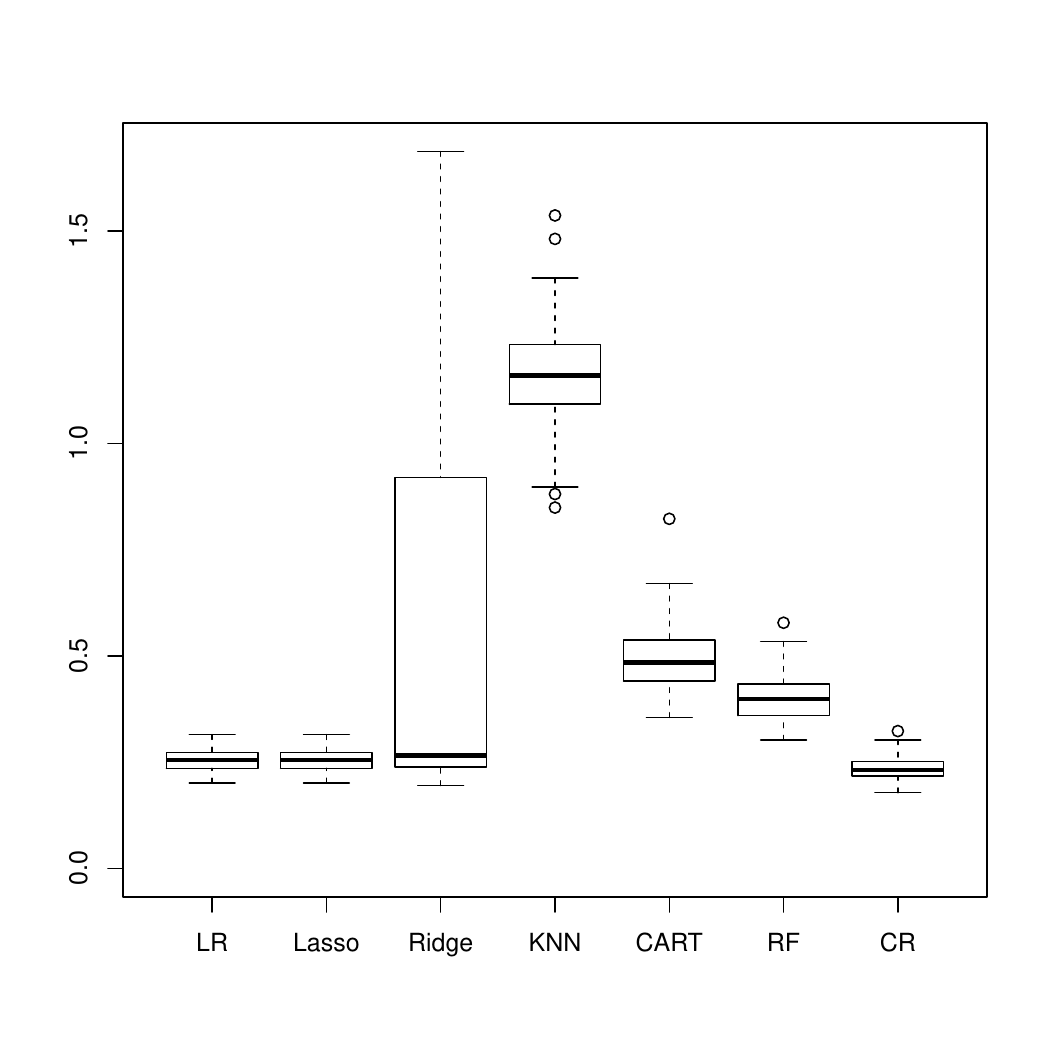}}
  \subfloat[\autoref{m4}.]{\includegraphics[width = .24\textwidth]{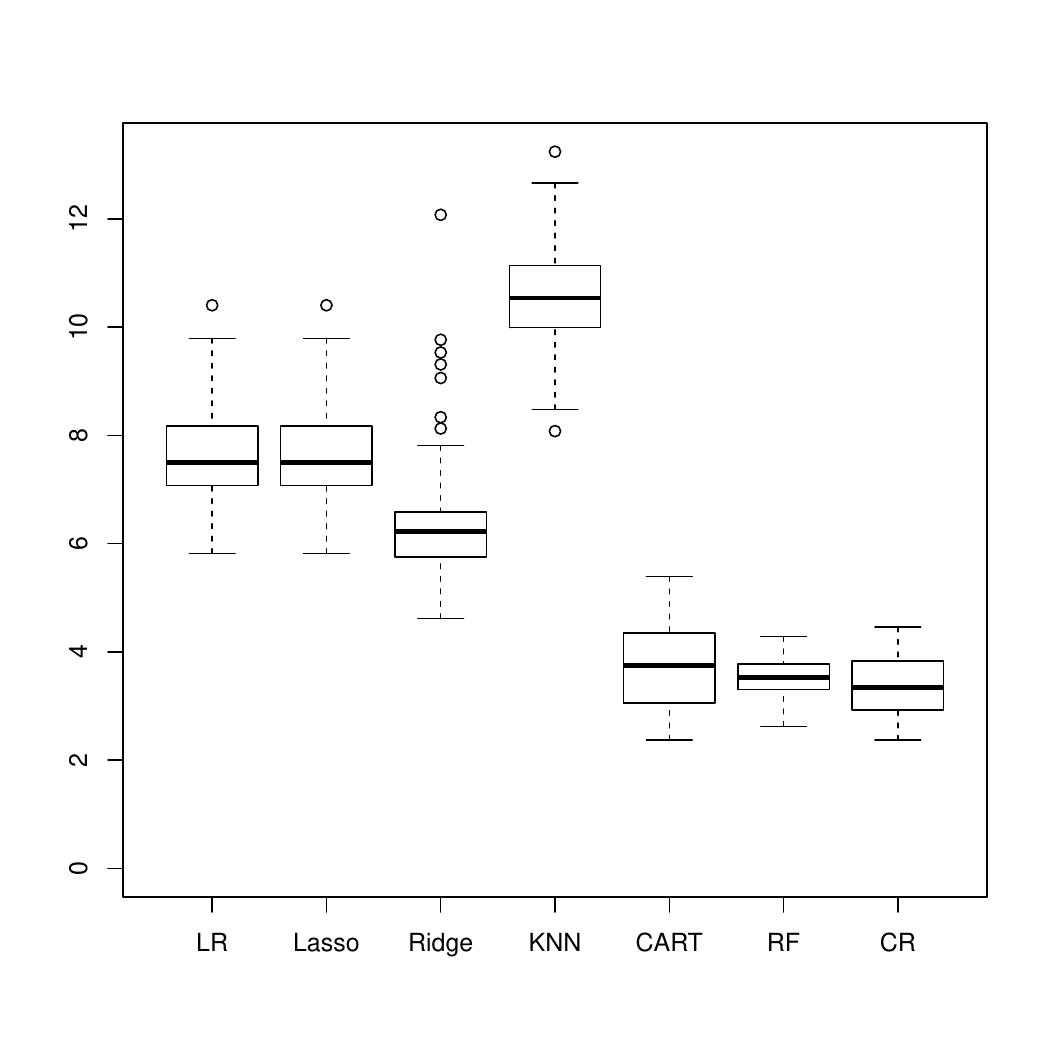}}
  \\
  \subfloat[\autoref{m5}.]{\includegraphics[width = .24\textwidth]{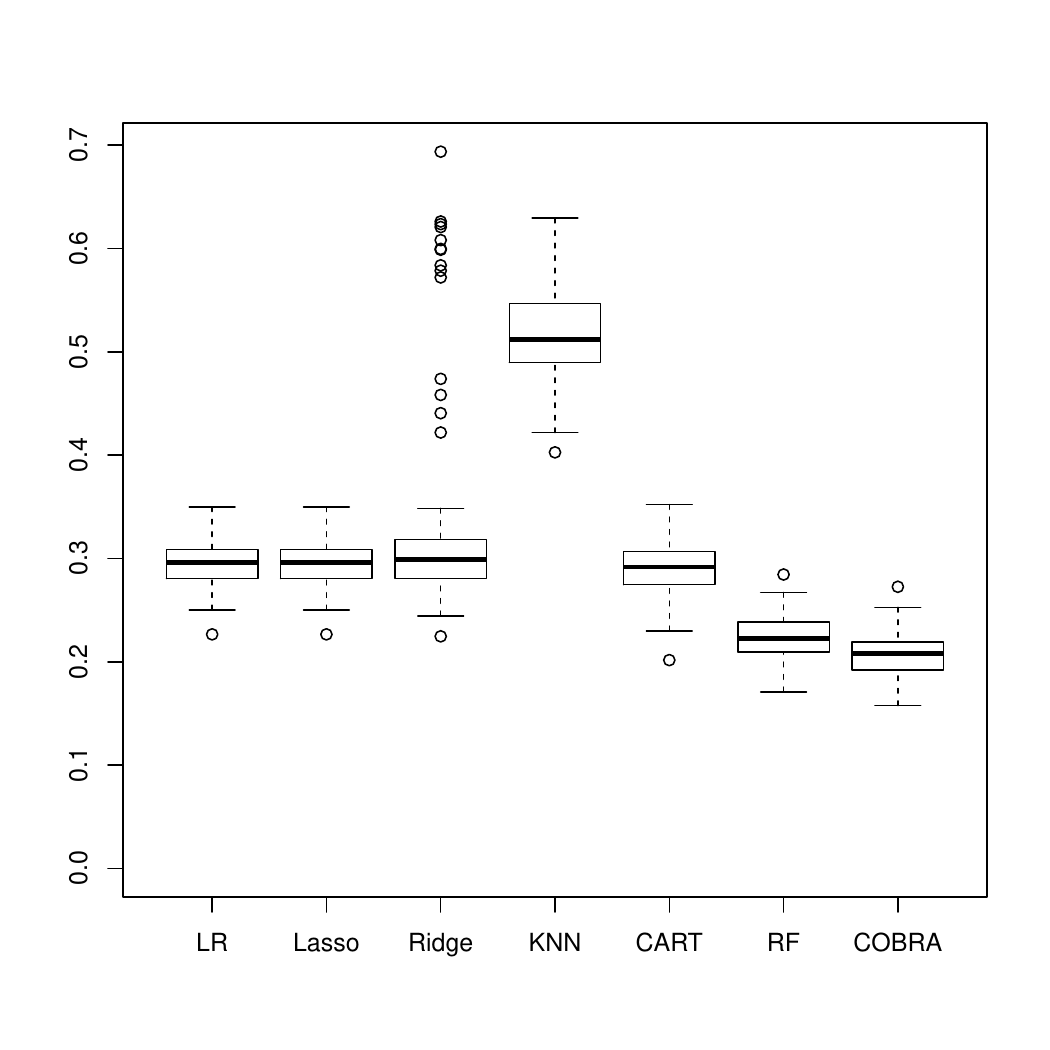}}
  \subfloat[\autoref{m6}.]{\includegraphics[width = .24\textwidth]{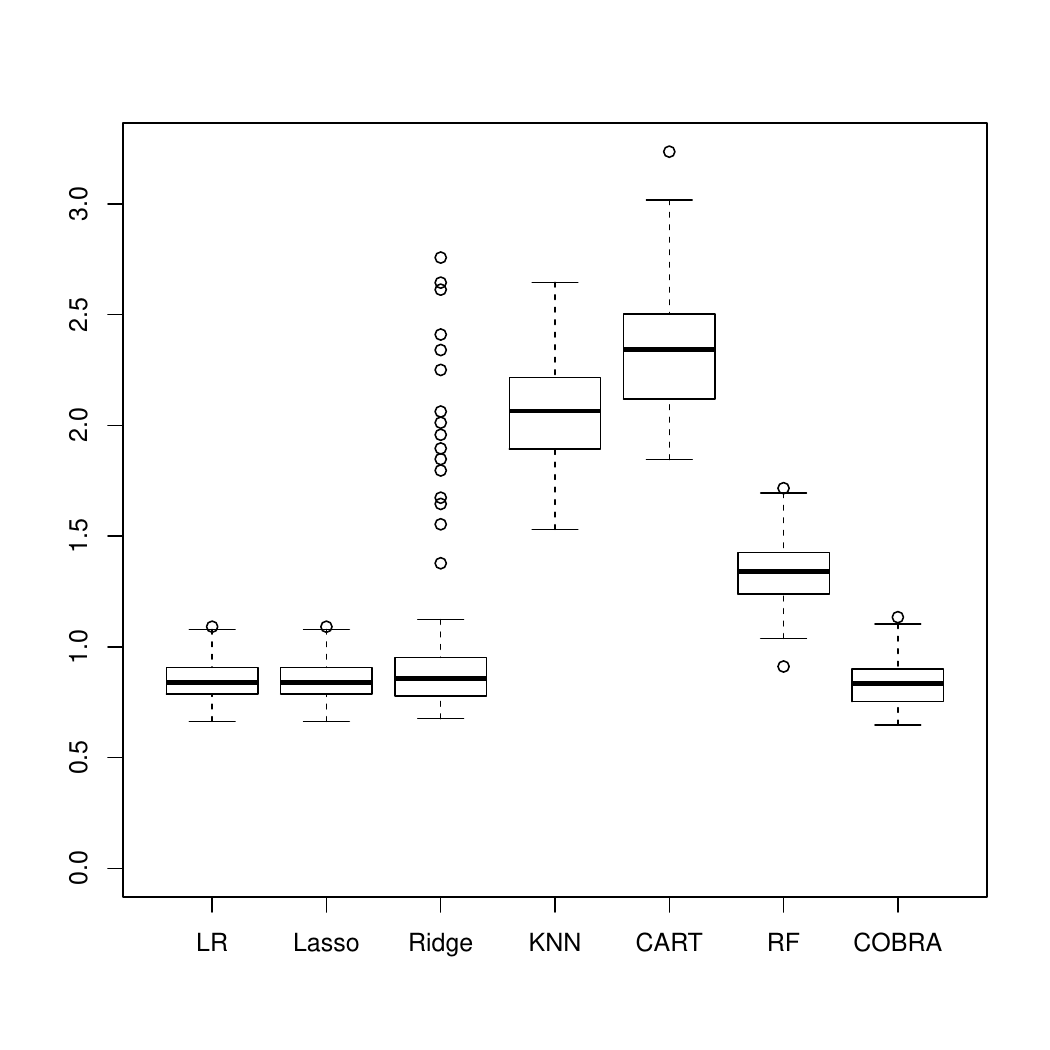}}
  \subfloat[\autoref{m7}.]{\includegraphics[width = .24\textwidth]{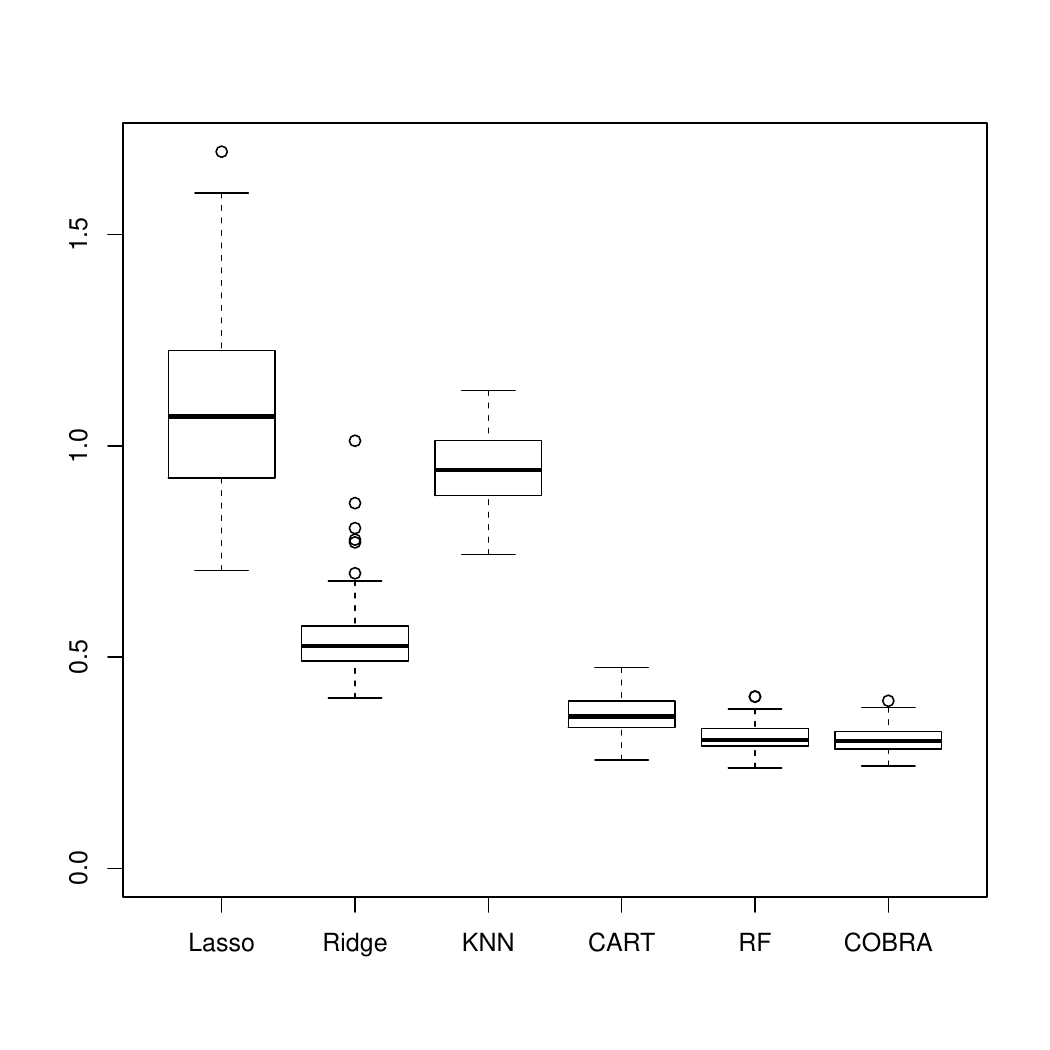}}
  \subfloat[\autoref{m8}.]{\includegraphics[width = .24\textwidth]{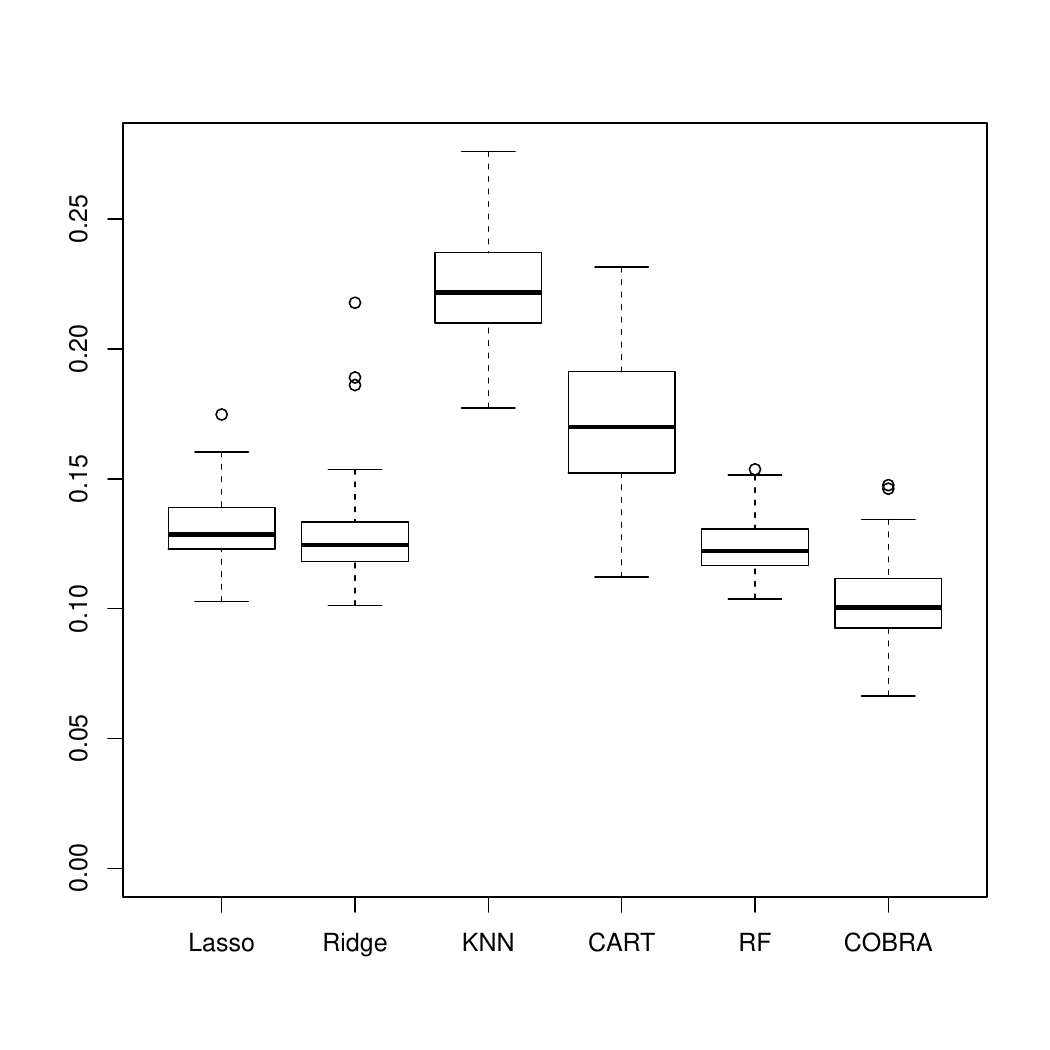}}
\end{figure}
\begin{figure}[t]
  \caption{Boxplots of quadratic errors, correlated design. From
    left to right: \texttt{lars}, \texttt{ridge}, \texttt{fnn},
    \texttt{tree}, \texttt{randomForest}, \cobra.}
  \label{boxplot-C}
  \subfloat[\autoref{m1}.]{\includegraphics[width = .24\textwidth]{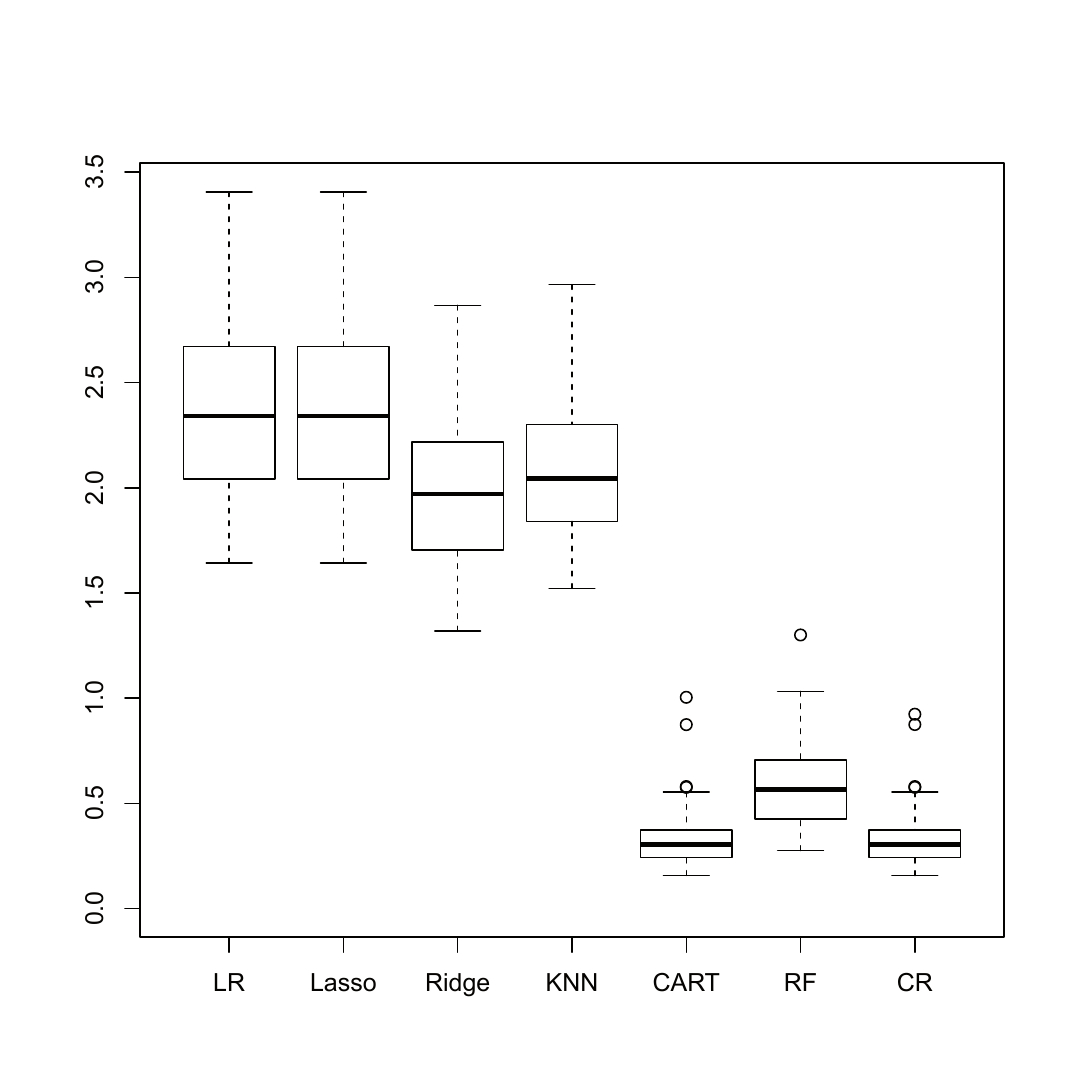}}
  \subfloat[\autoref{m2}.]{\includegraphics[width = .24\textwidth]{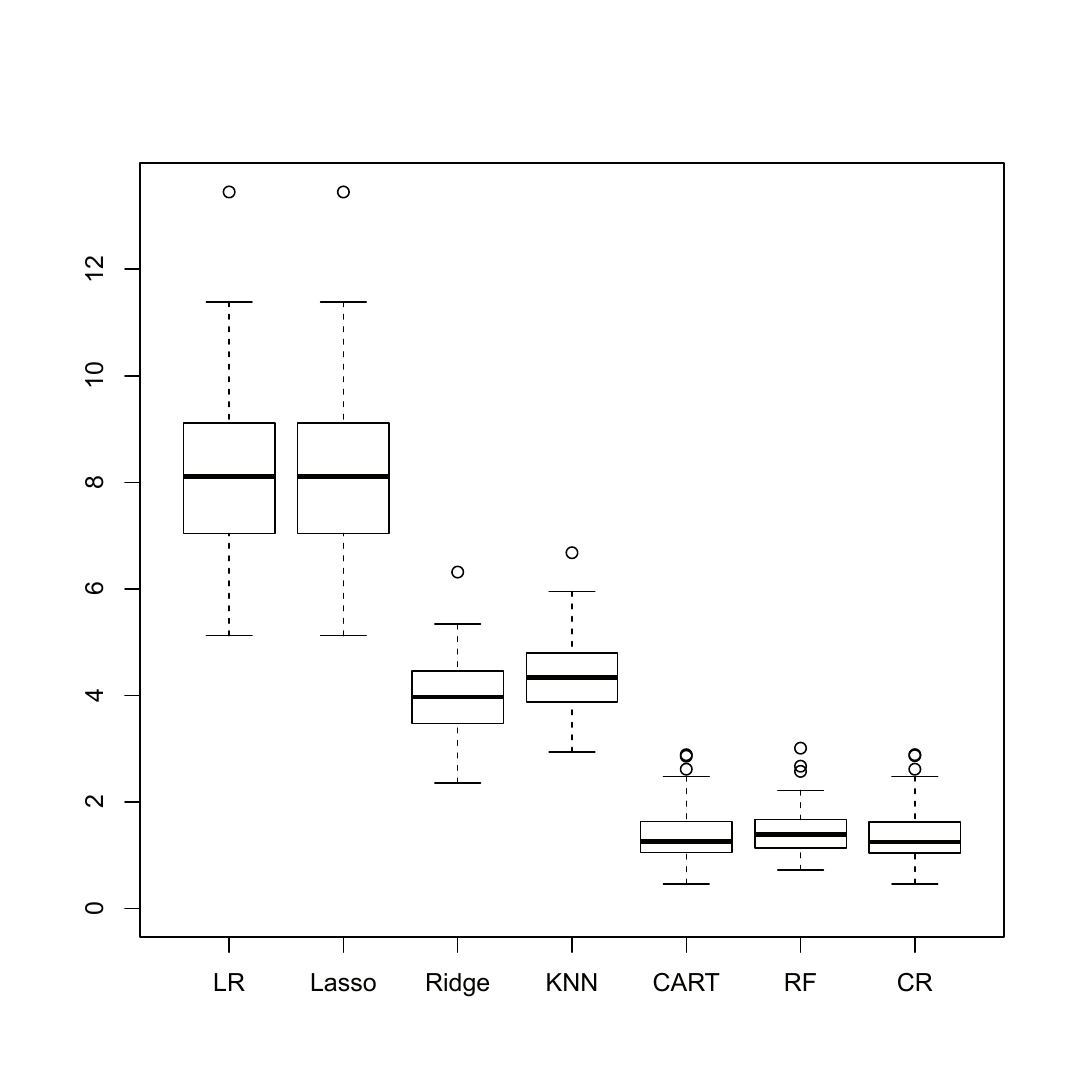}}
  \subfloat[\autoref{m3}.]{\includegraphics[width = .24\textwidth]{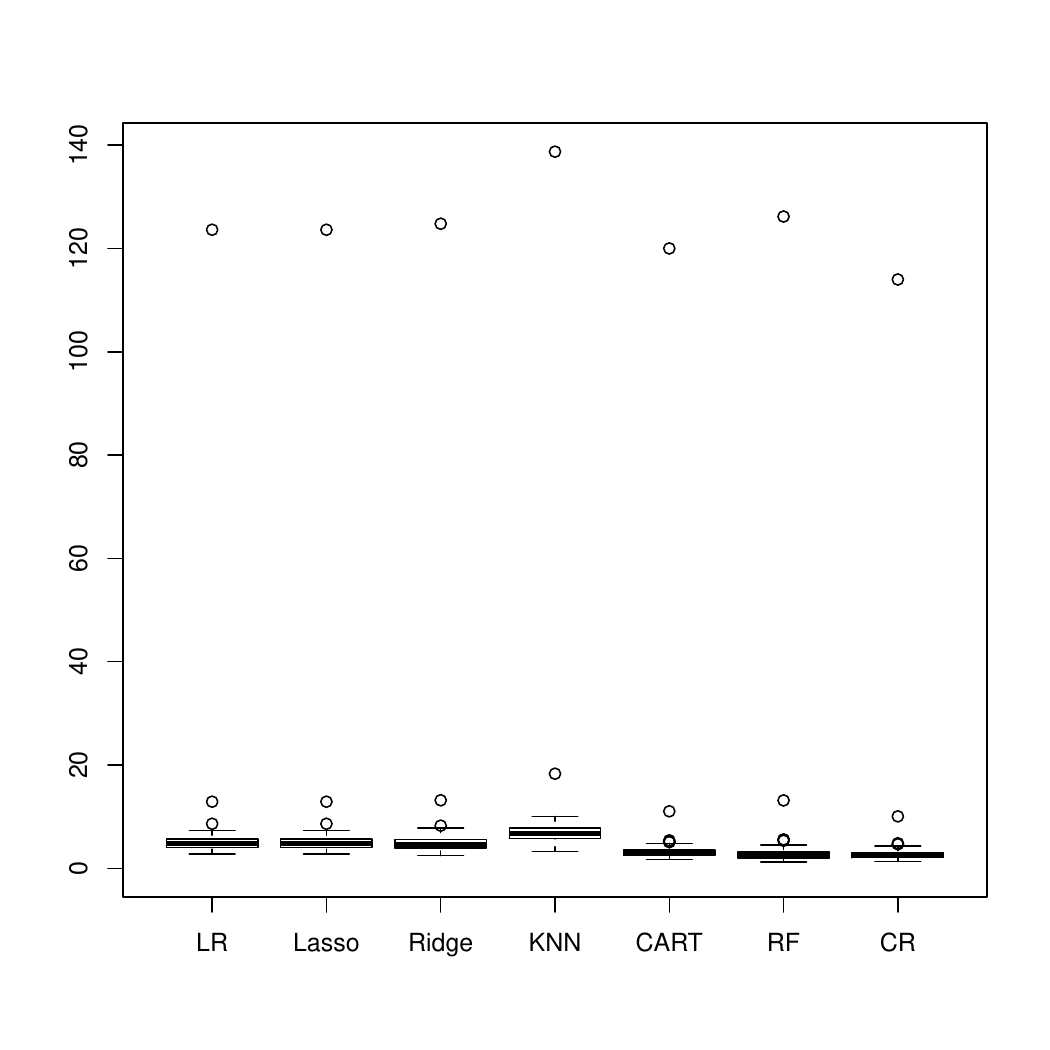}}
  \subfloat[\autoref{m4}.]{\includegraphics[width = .24\textwidth]{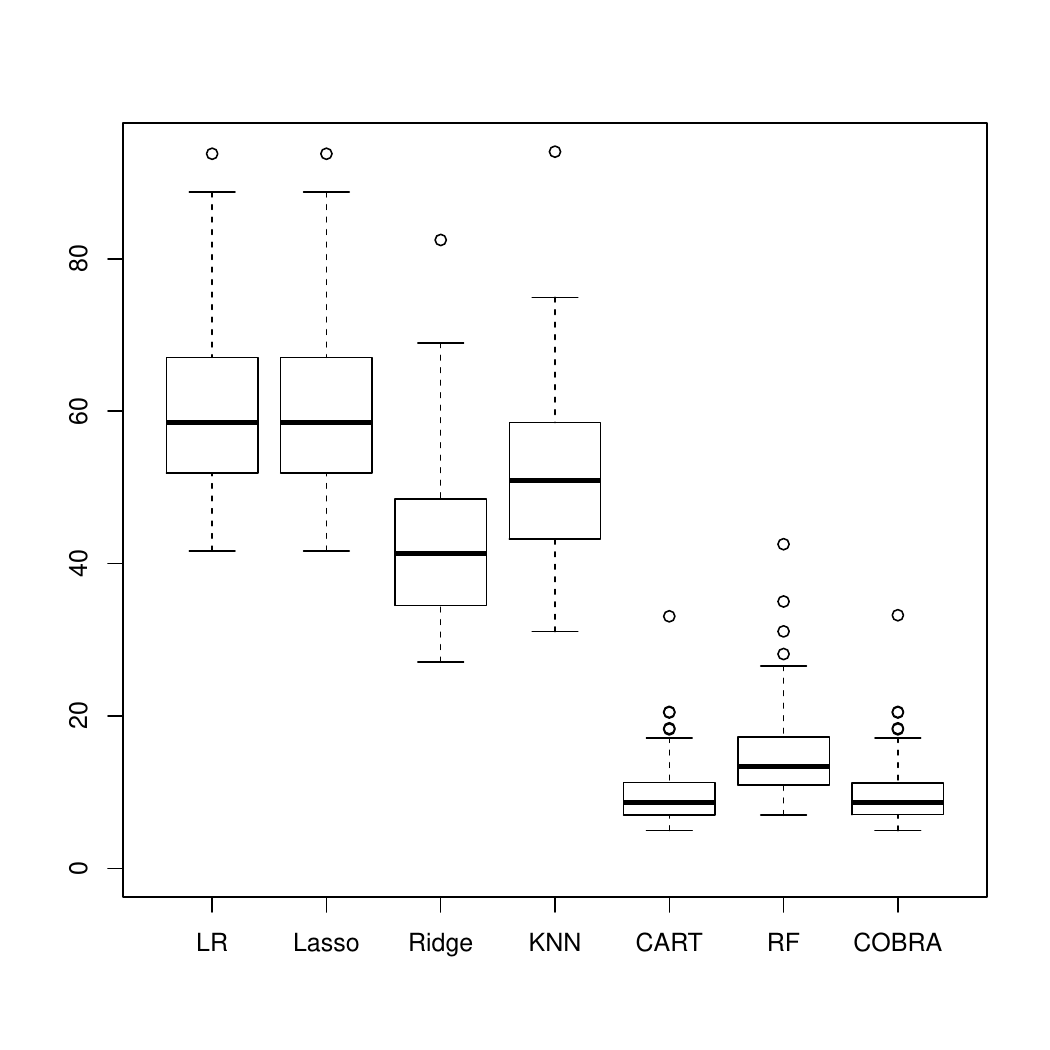}}
  \\
  \subfloat[\autoref{m5}.]{\includegraphics[width = .24\textwidth]{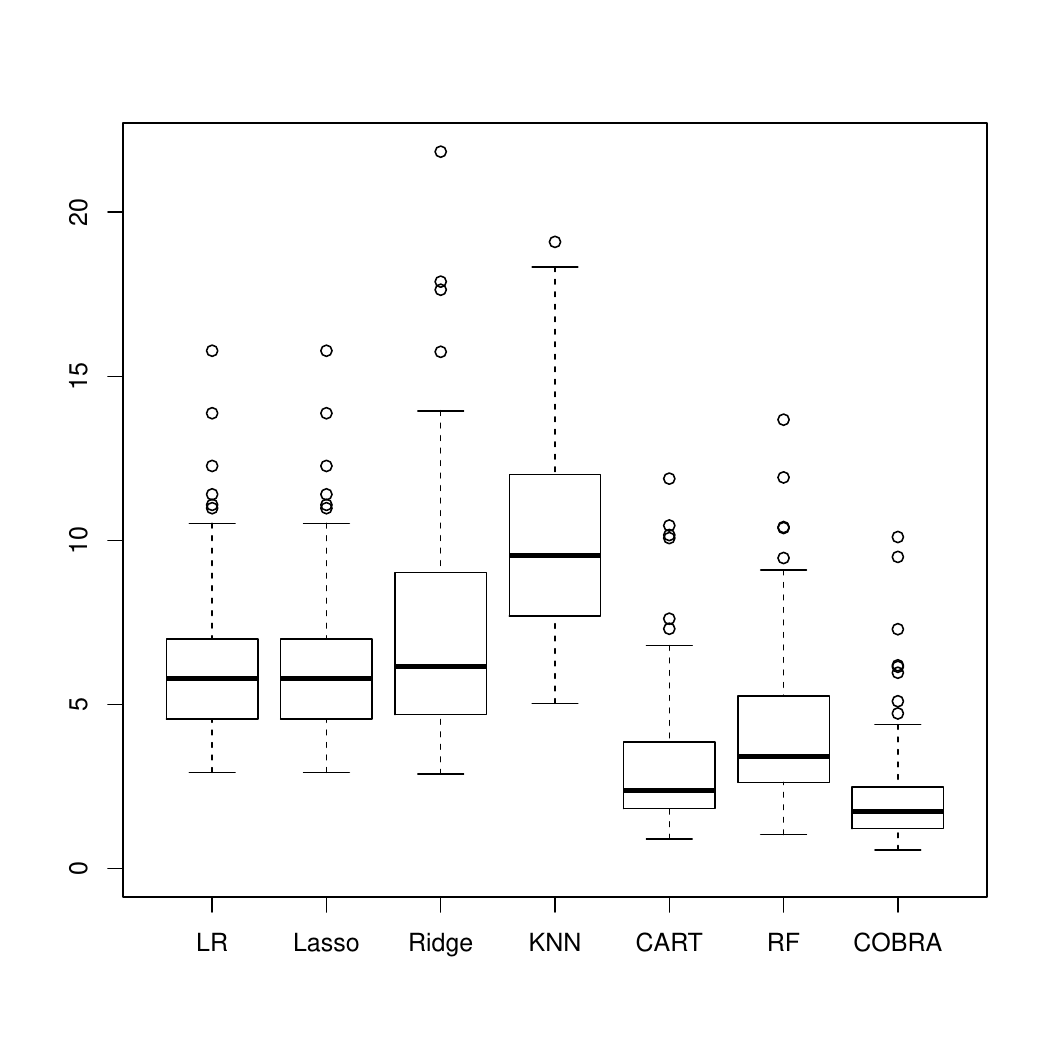}}
  \subfloat[\autoref{m6}.]{\includegraphics[width = .24\textwidth]{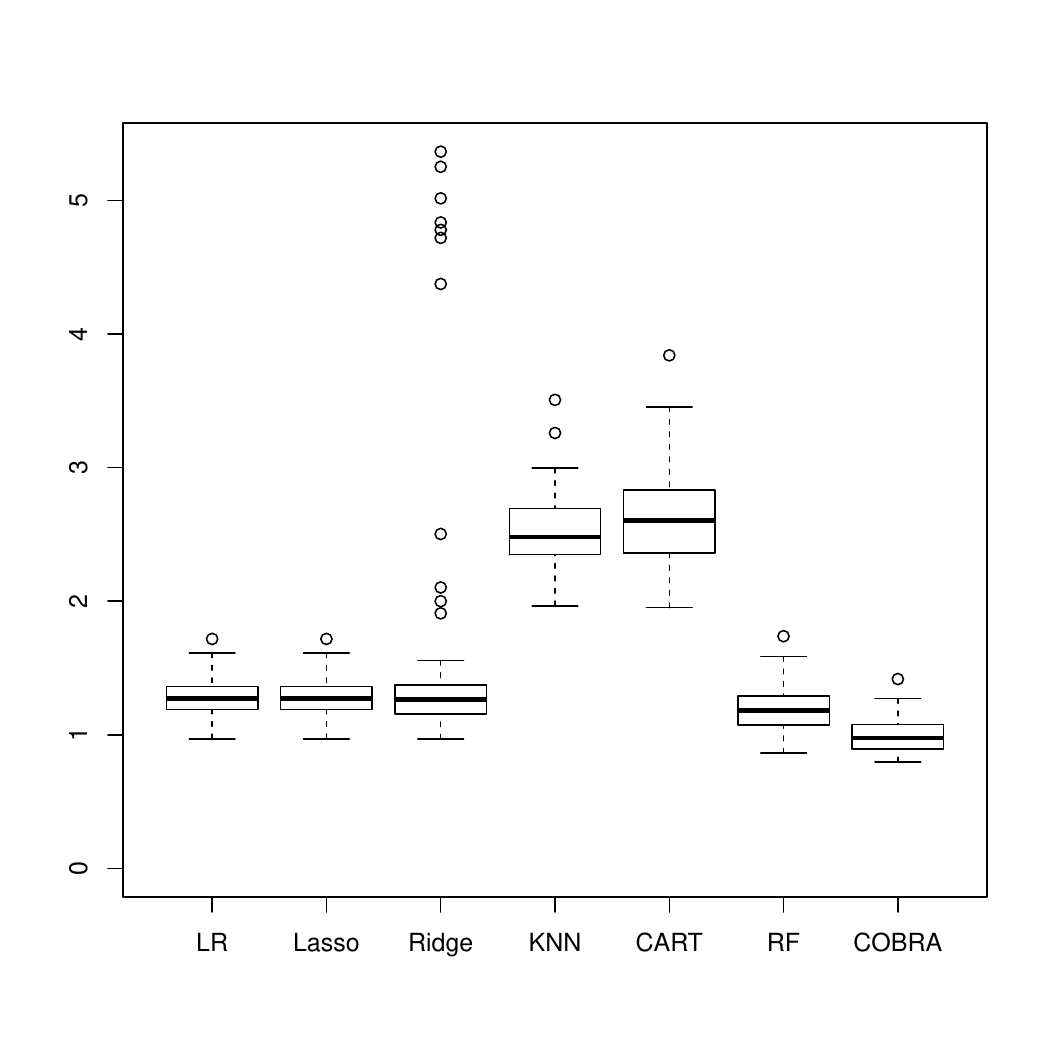}}
  \subfloat[\autoref{m7}.]{\includegraphics[width = .24\textwidth]{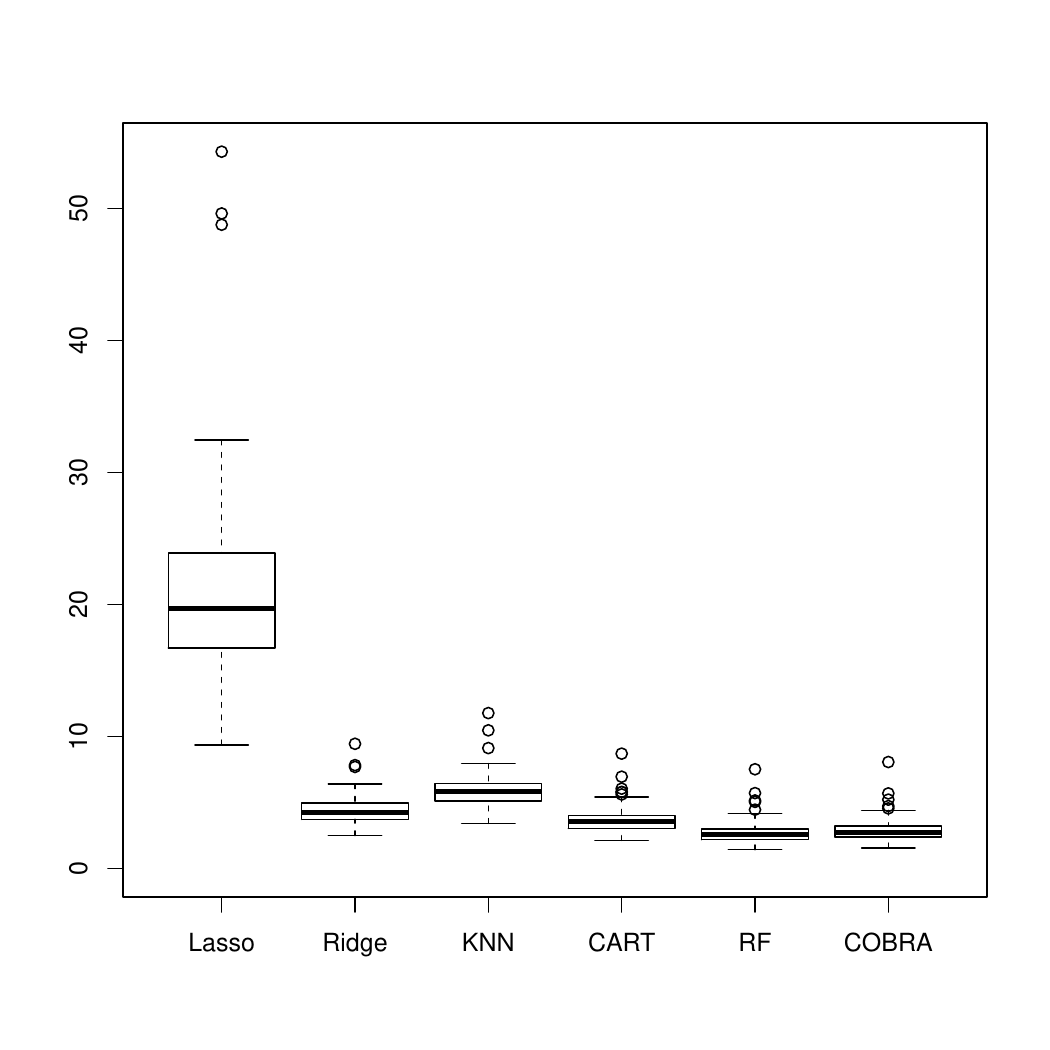}}
  \subfloat[\autoref{m8}.]{\includegraphics[width = .24\textwidth]{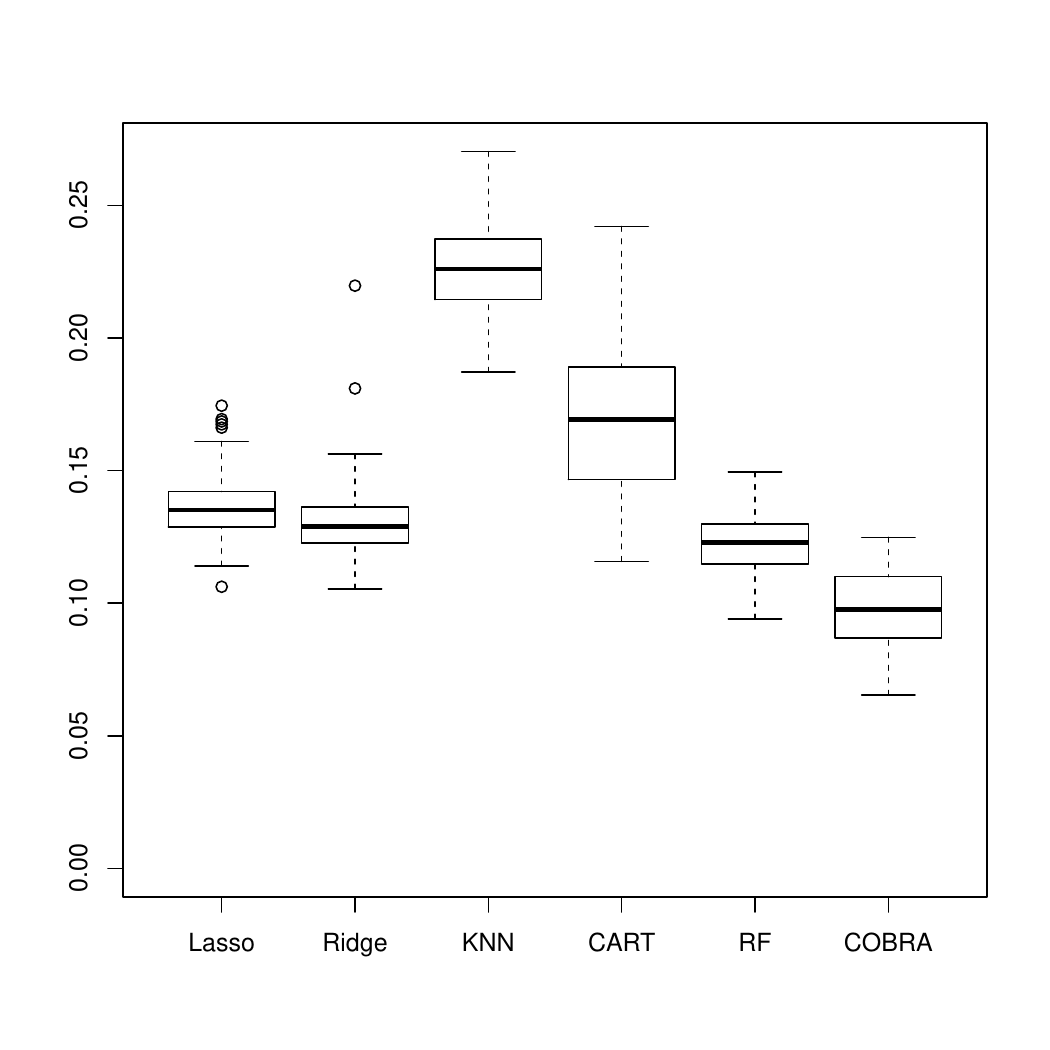}}
\end{figure}

\begin{figure}[t]
  \caption{Prediction over the testing set, uncorrelated design. The more points on the
    first bissectrix, the better the prediction.}
  \label{pred-U}
  \subfloat[\autoref{m1}.]{\includegraphics[width = .24\textwidth]{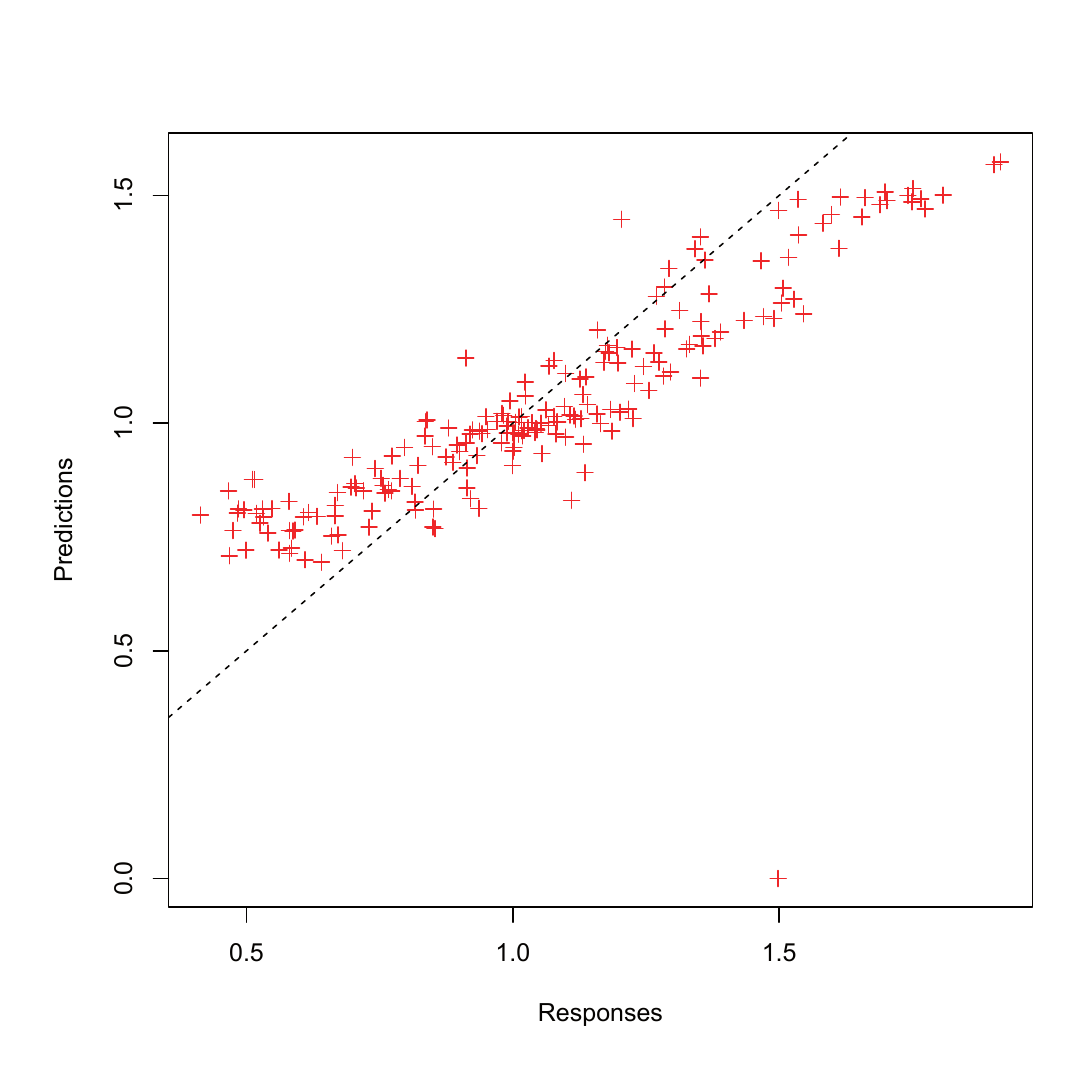}}
  \subfloat[\autoref{m2}.]{\includegraphics[width = .24\textwidth]{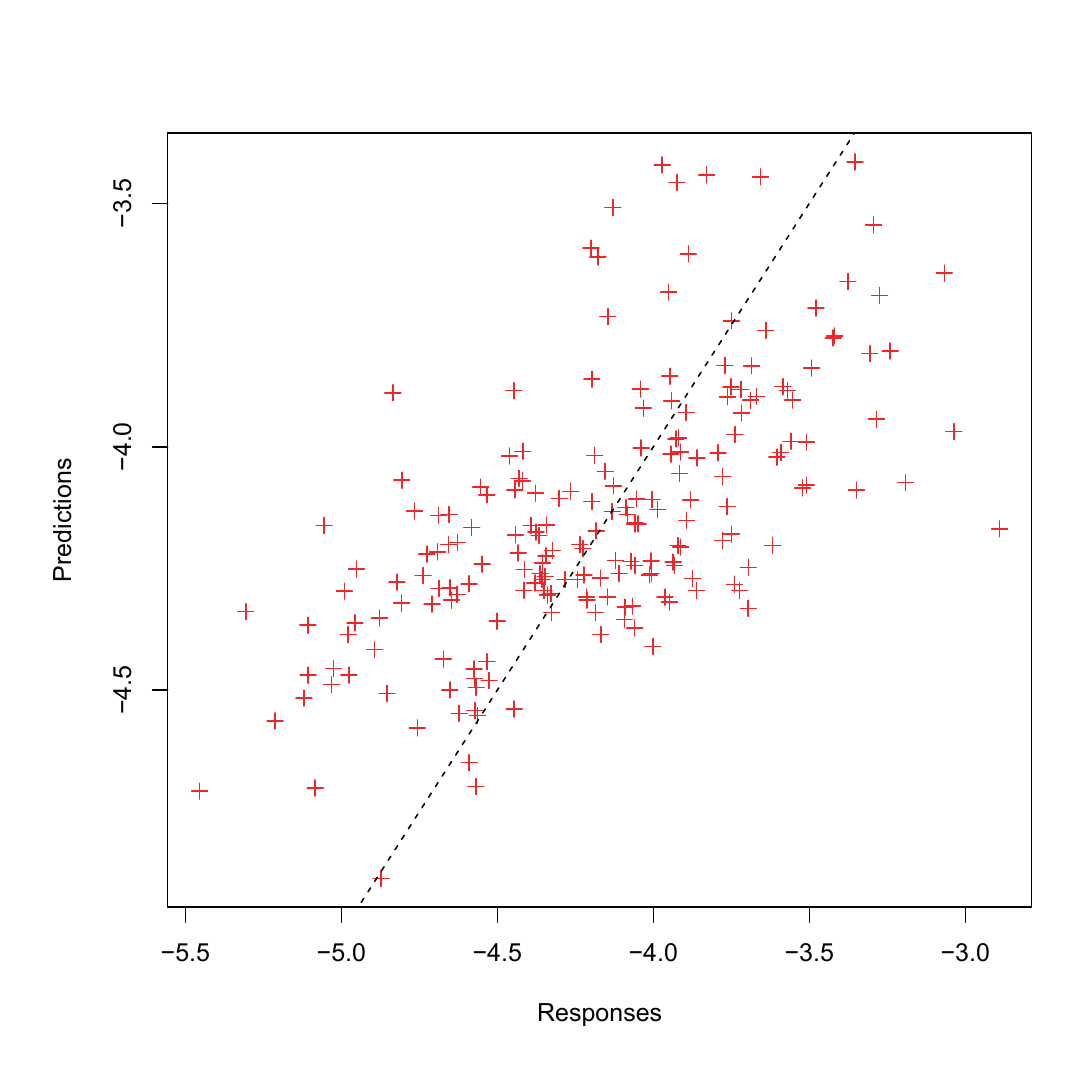}}
  \subfloat[\autoref{m3}.]{\includegraphics[width = .24\textwidth]{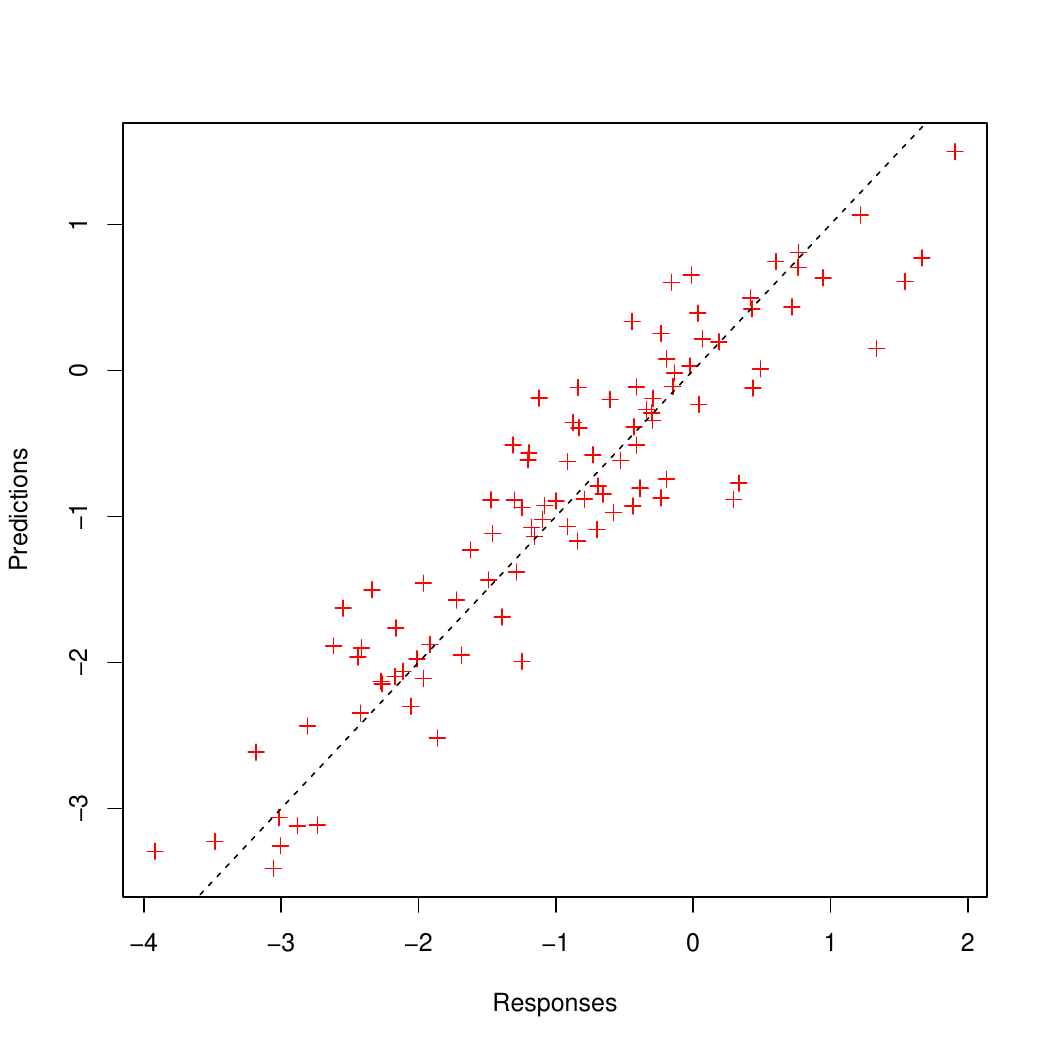}}
  \subfloat[\autoref{m4}.]{\includegraphics[width = .24\textwidth]{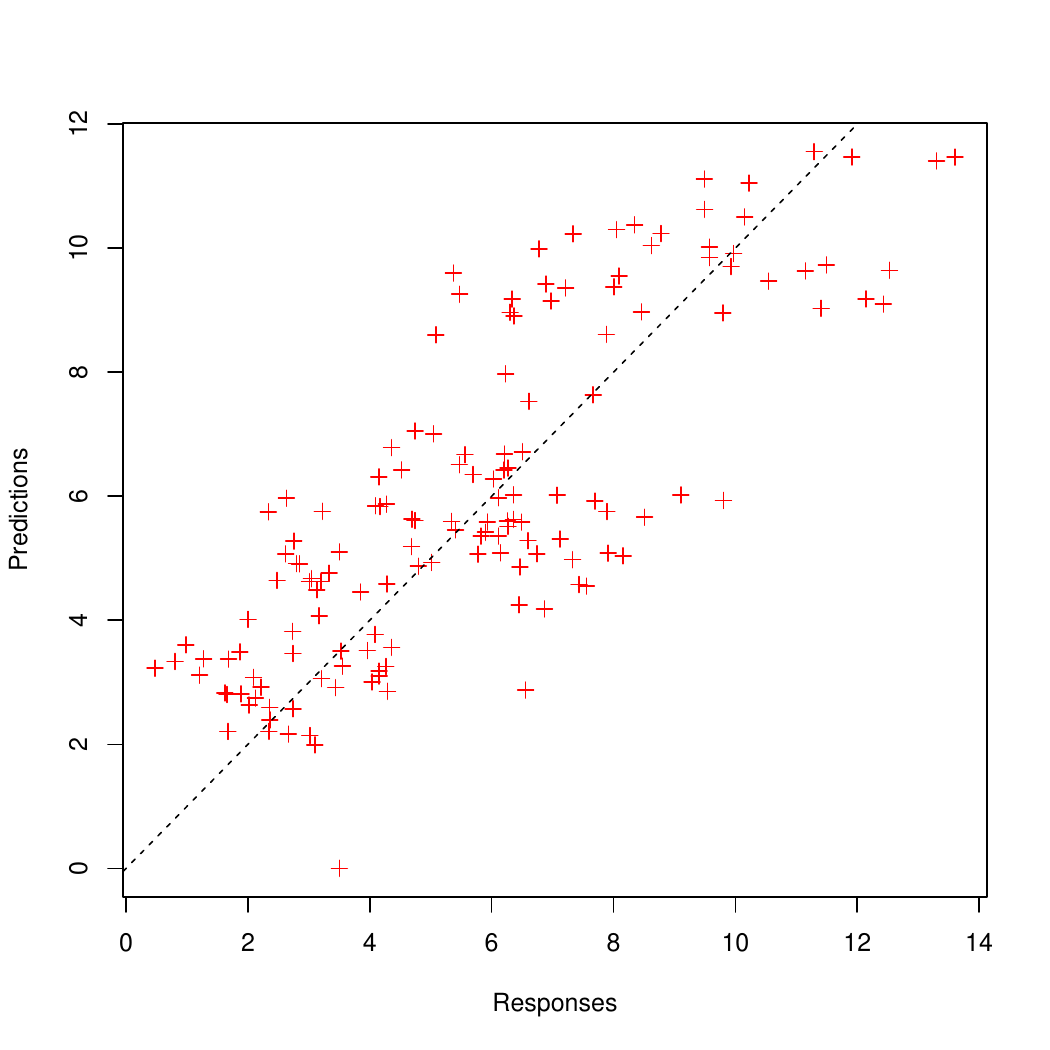}}
  \\
  \subfloat[\autoref{m5}.]{\includegraphics[width = .24\textwidth]{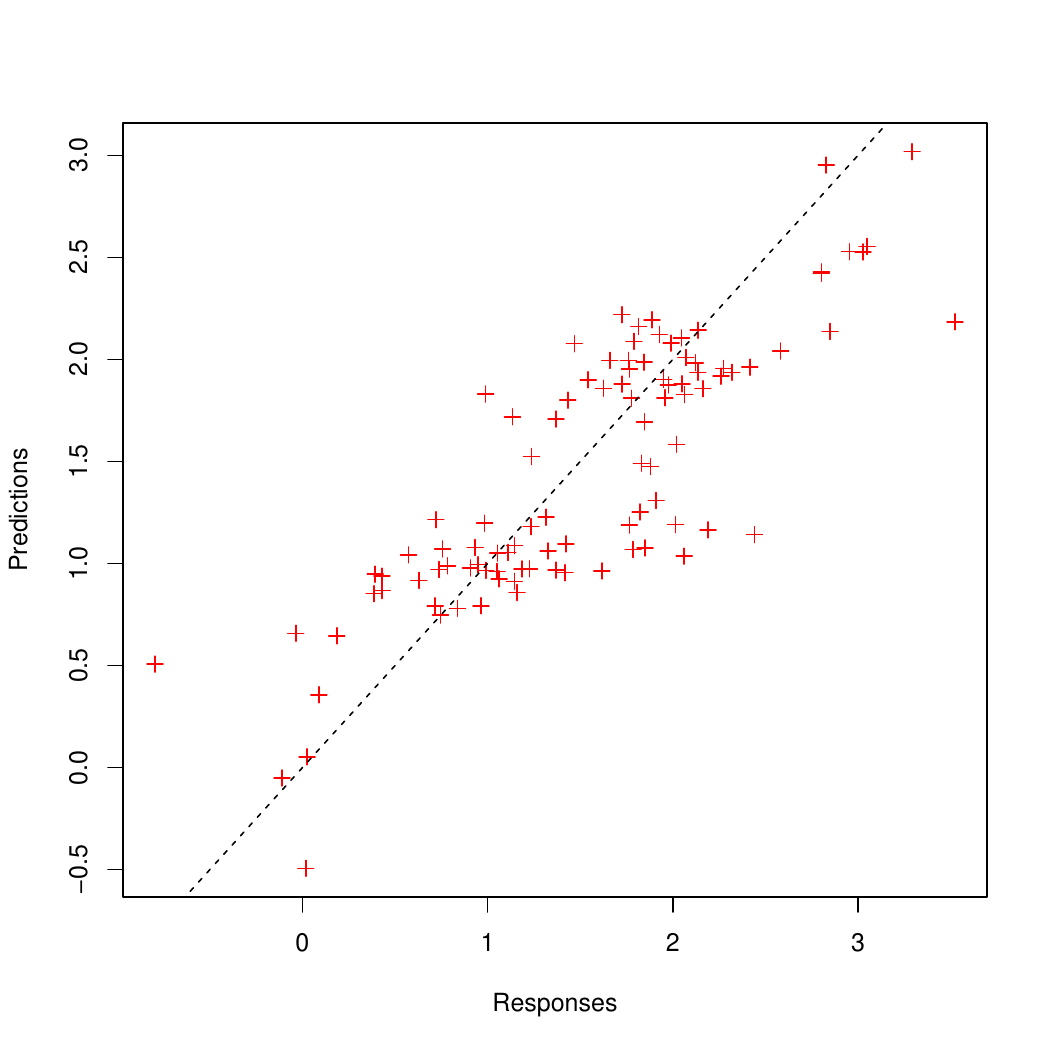}}
  \subfloat[\autoref{m6}.]{\includegraphics[width = .24\textwidth]{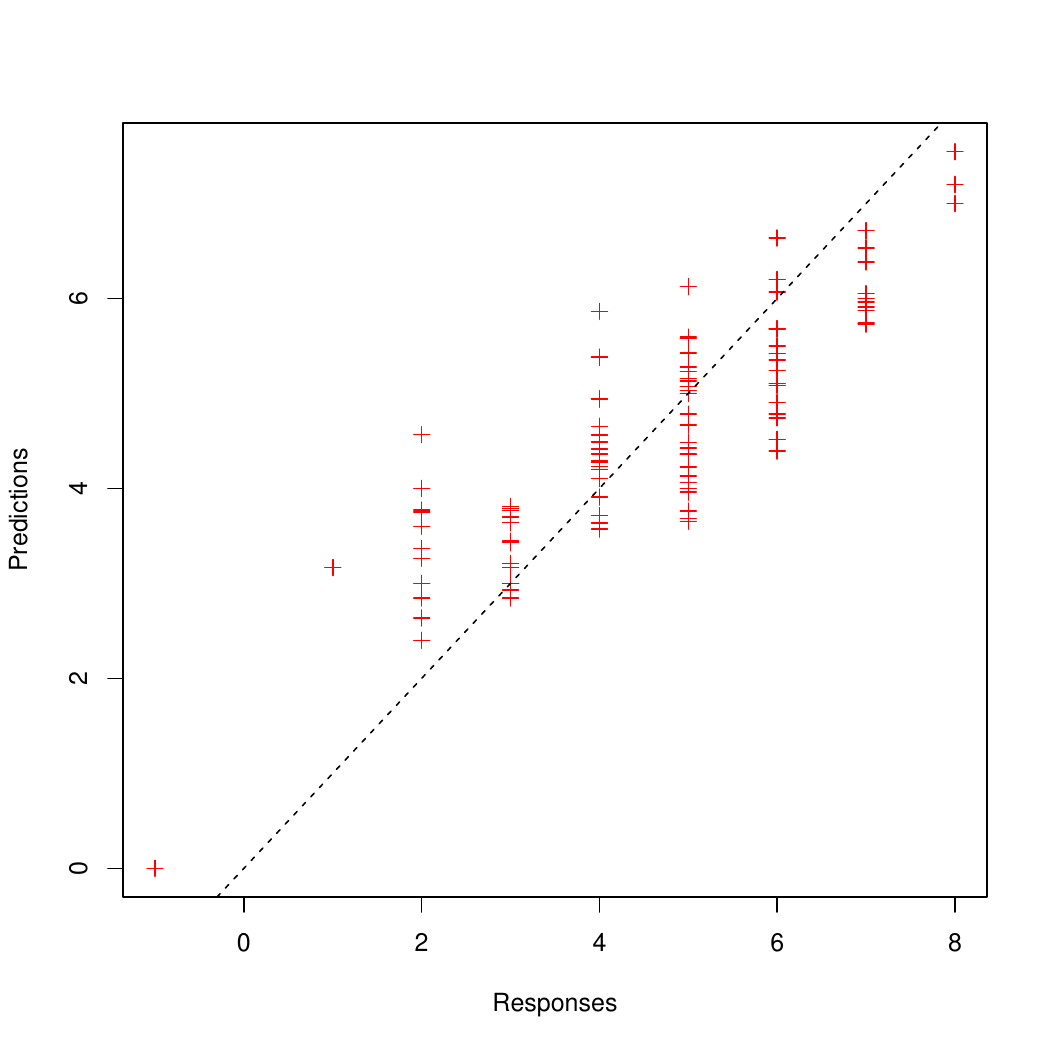}}
  \subfloat[\autoref{m7}.]{\includegraphics[width = .24\textwidth]{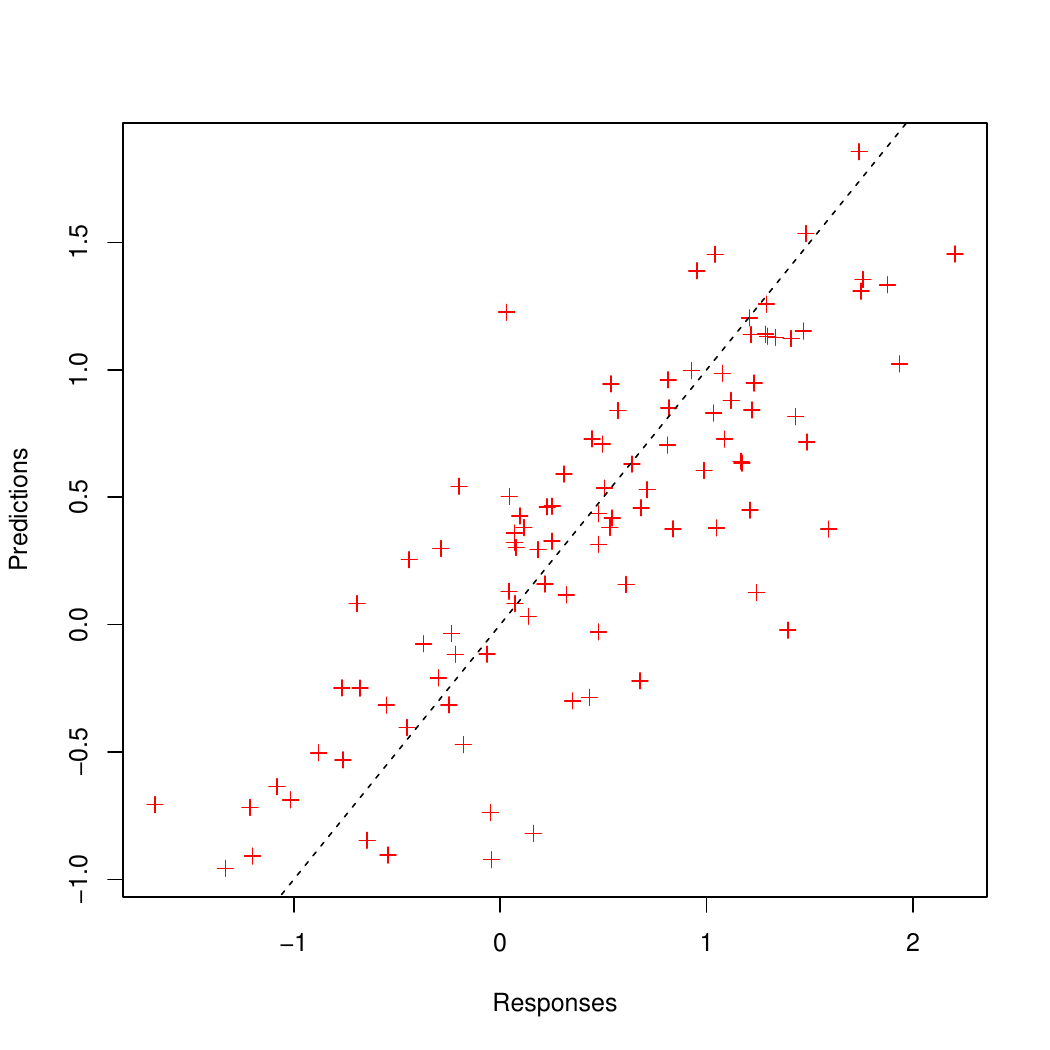}}
  \subfloat[\autoref{m8}.]{\includegraphics[width = .24\textwidth]{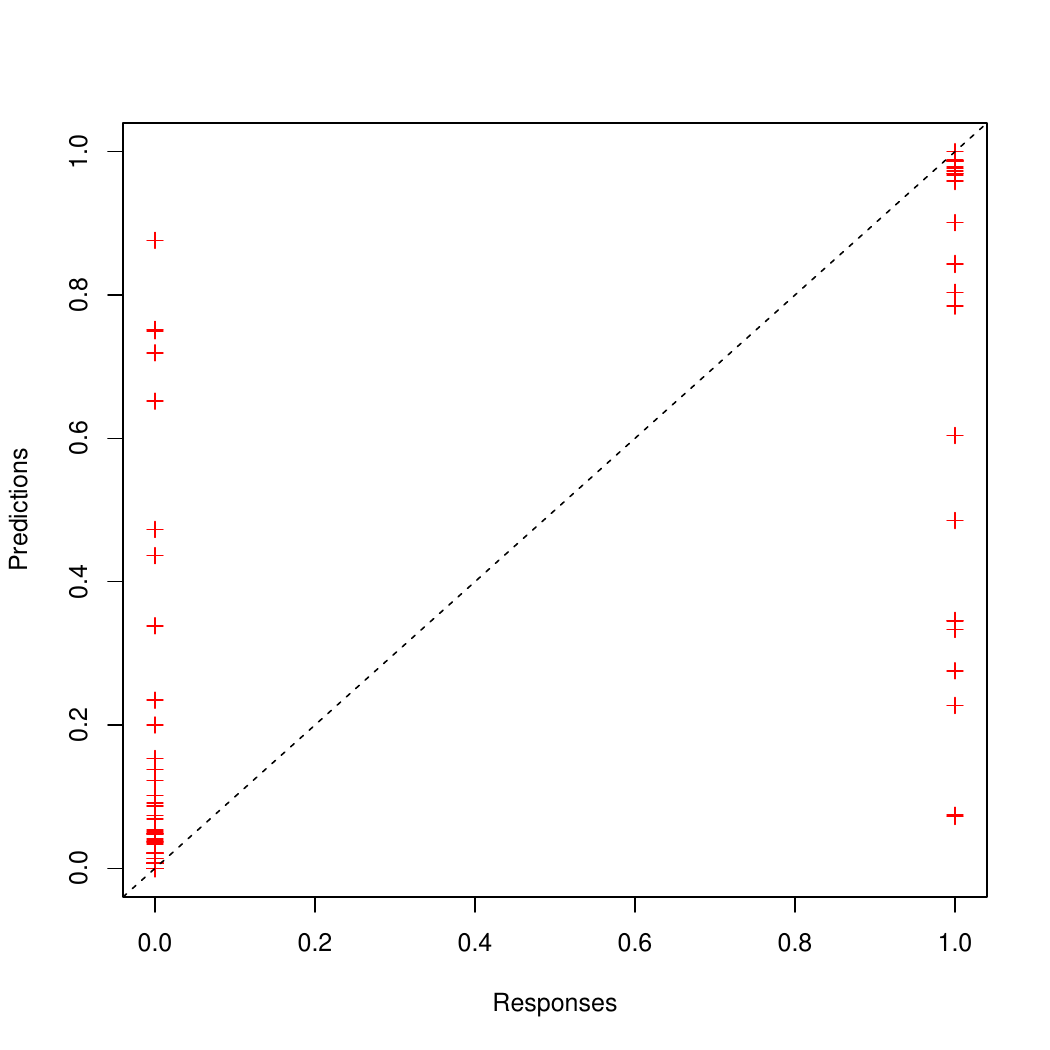}}
\end{figure}
\begin{figure}[t]
  \caption{Prediction over the testing set, correlated design. The more points on the
    first bissectrix, the better the prediction.}
  \label{pred-C}
  \subfloat[\autoref{m1}.]{\includegraphics[width = .24\textwidth]{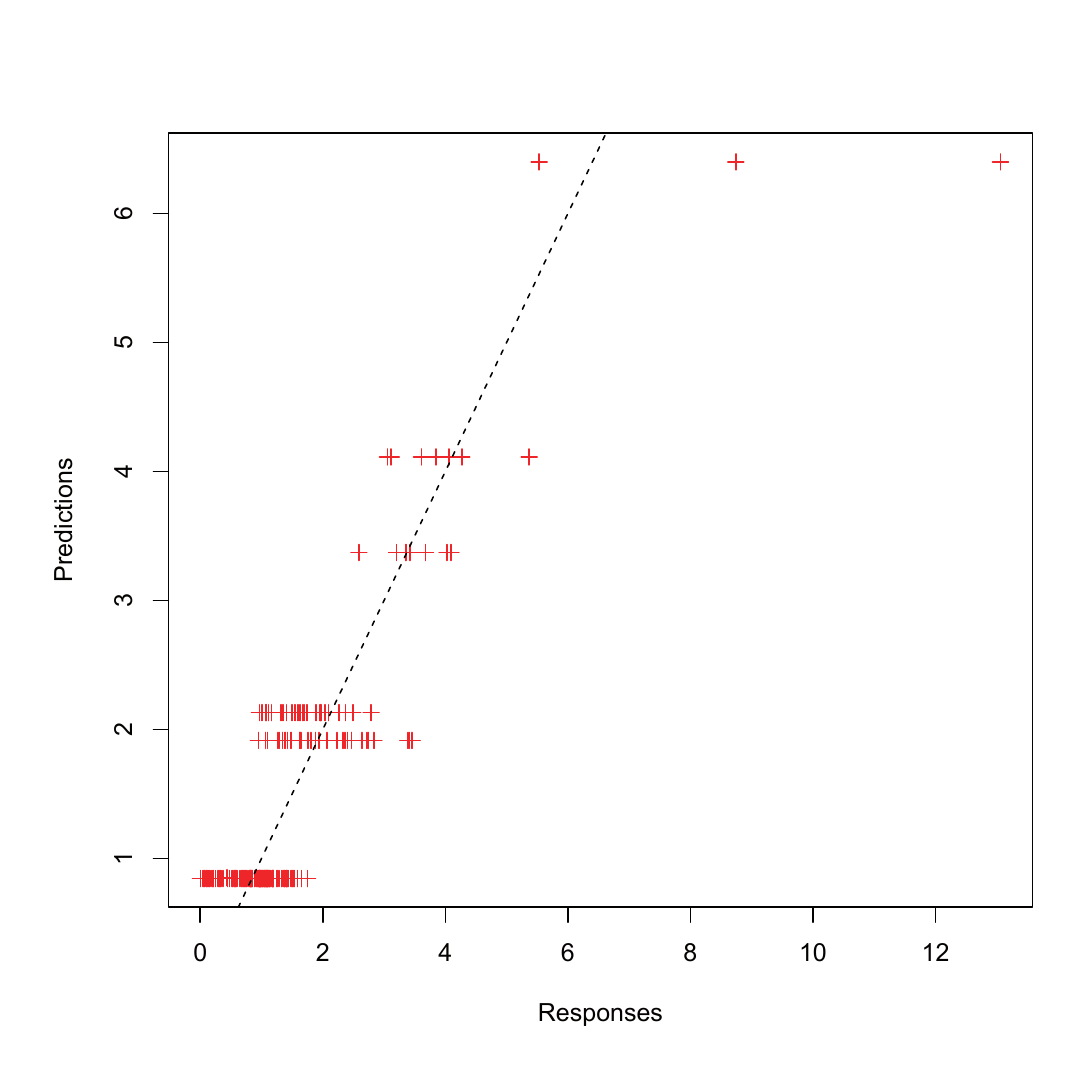}}
  \subfloat[\autoref{m2}.]{\includegraphics[width = .24\textwidth]{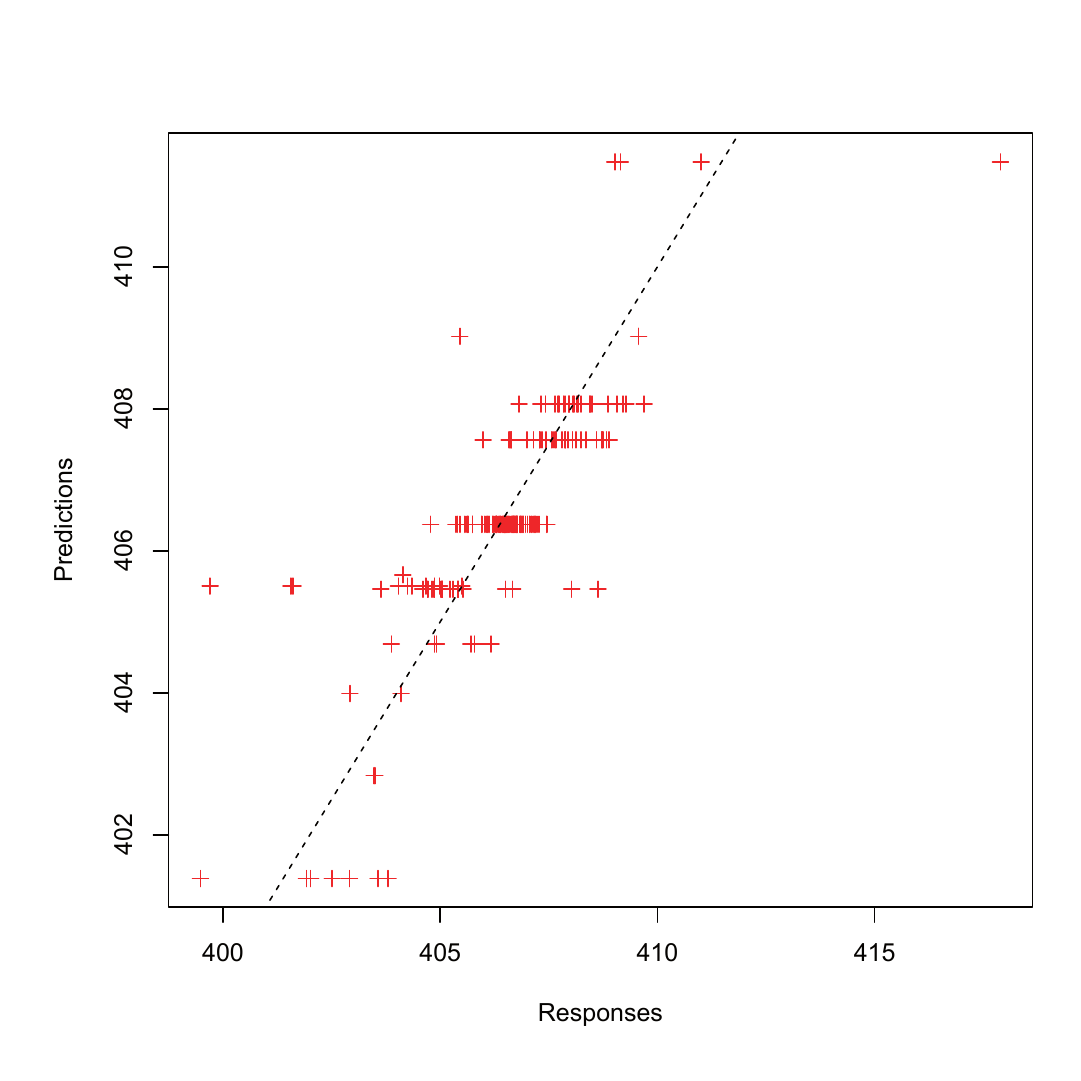}}
  \subfloat[\autoref{m3}.]{\includegraphics[width = .24\textwidth]{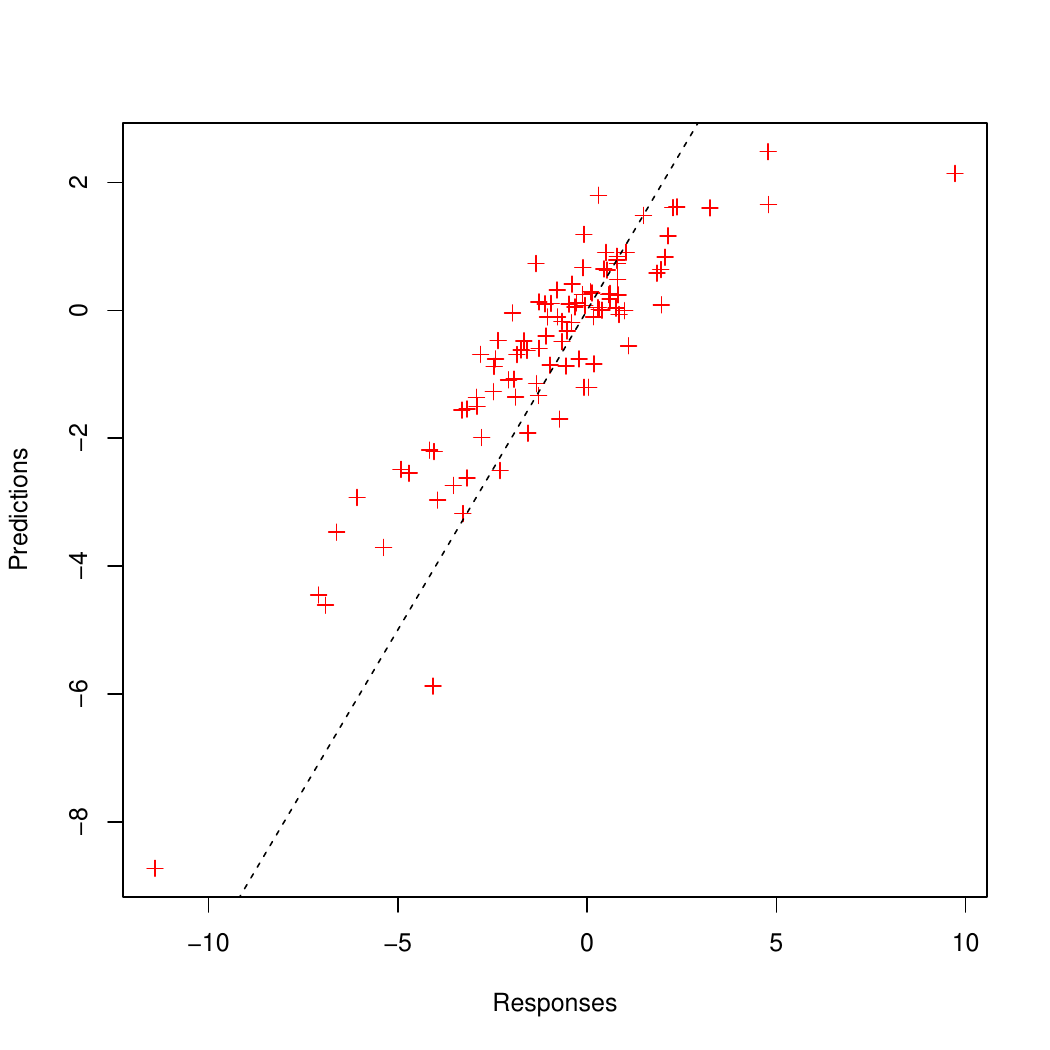}}
  \subfloat[\autoref{m4}.]{\includegraphics[width = .24\textwidth]{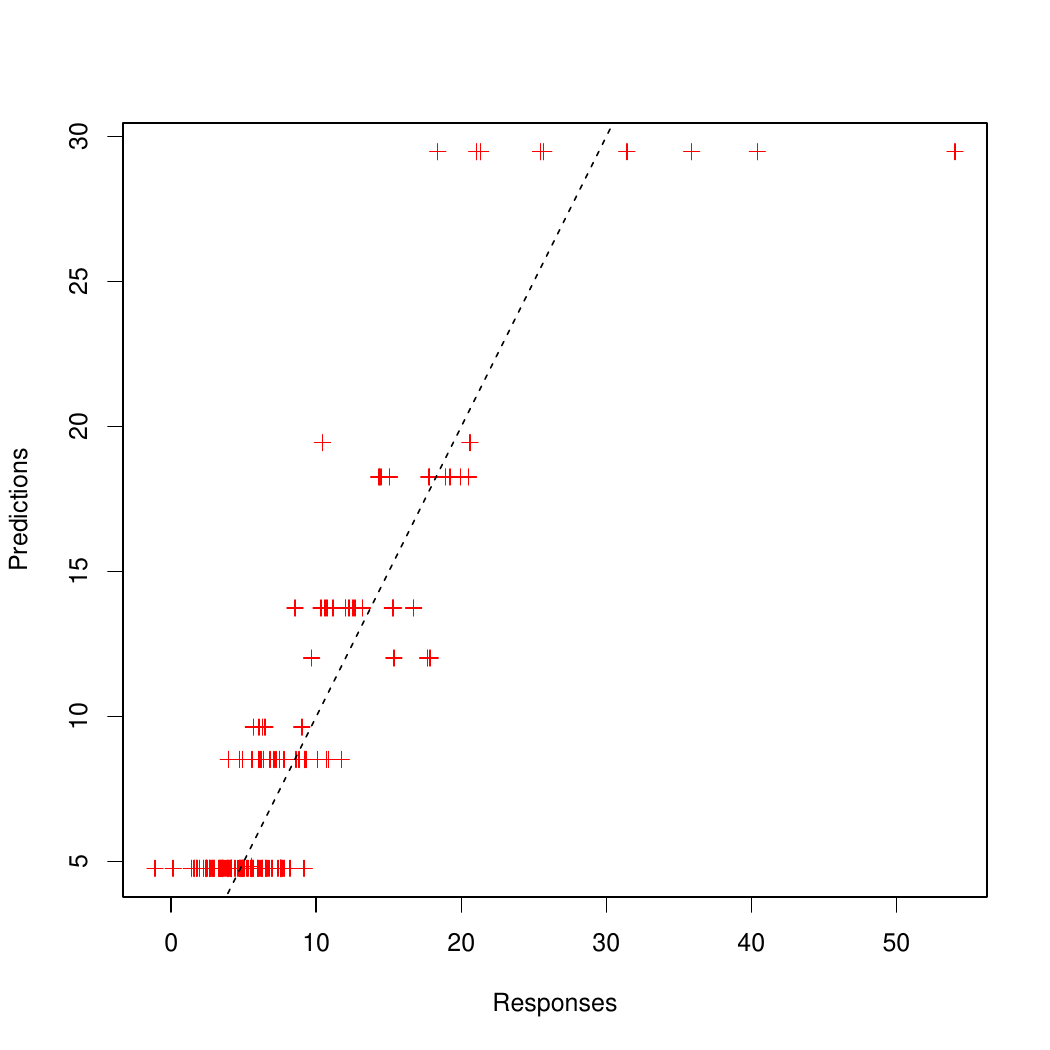}}
  \\
  \subfloat[\autoref{m5}.]{\includegraphics[width = .24\textwidth]{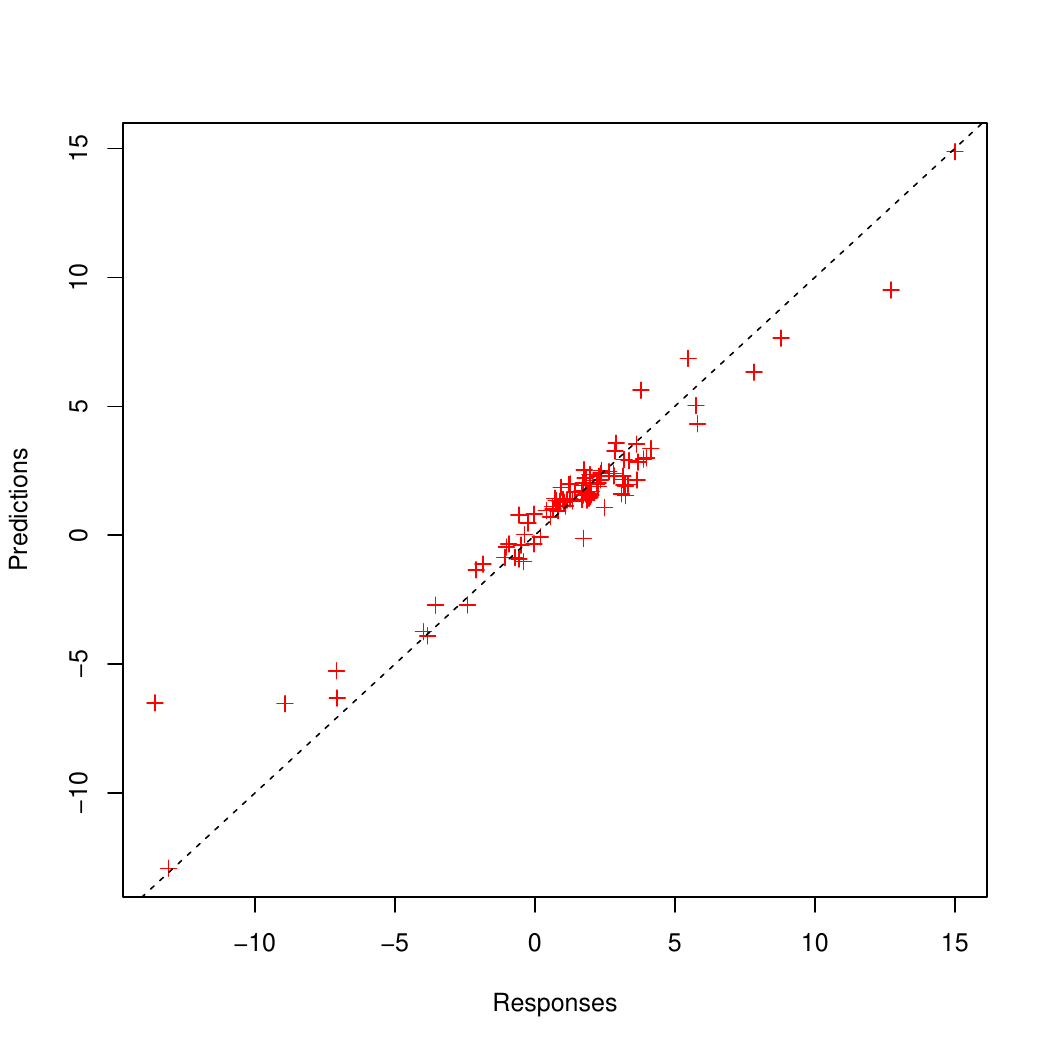}}
  \subfloat[\autoref{m6}.]{\includegraphics[width = .24\textwidth]{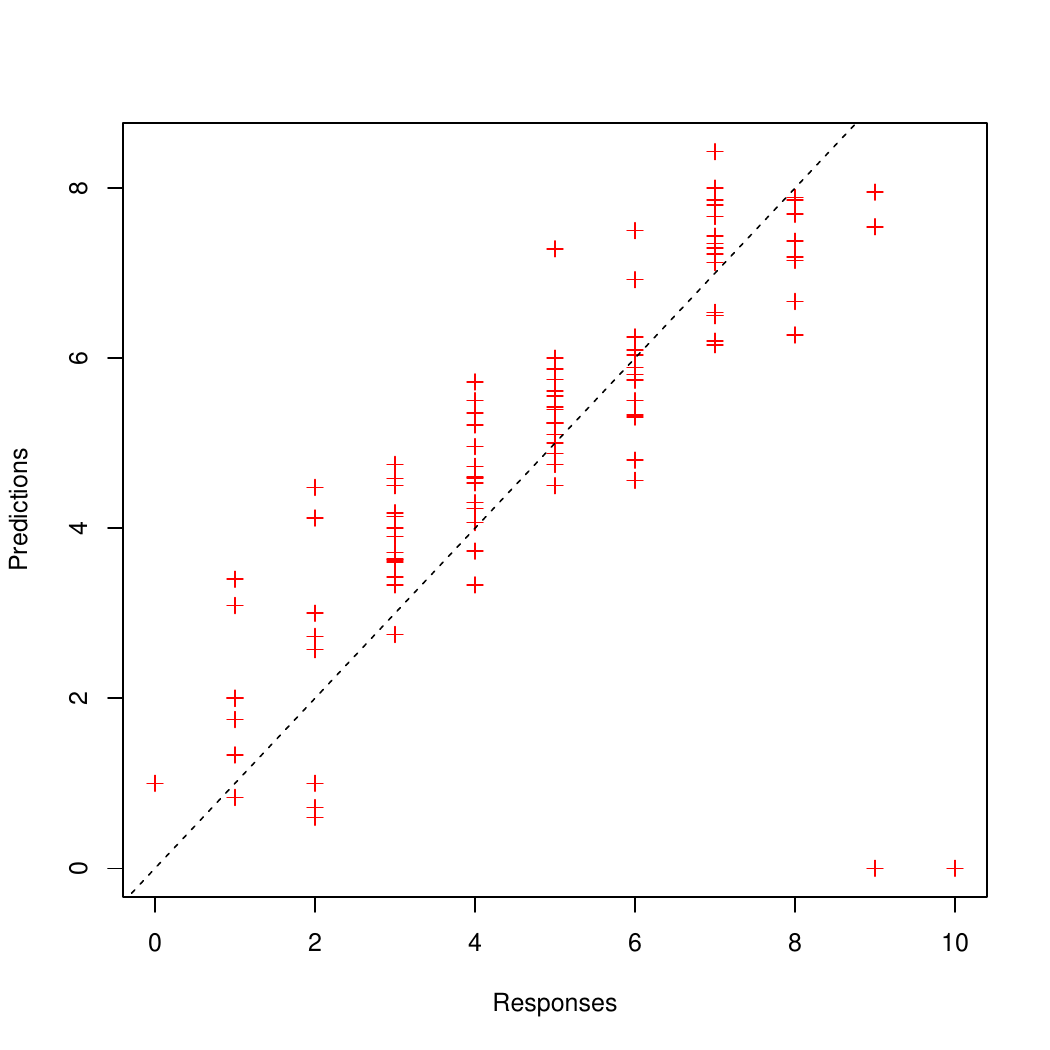}}
  \subfloat[\autoref{m7}.]{\includegraphics[width = .24\textwidth]{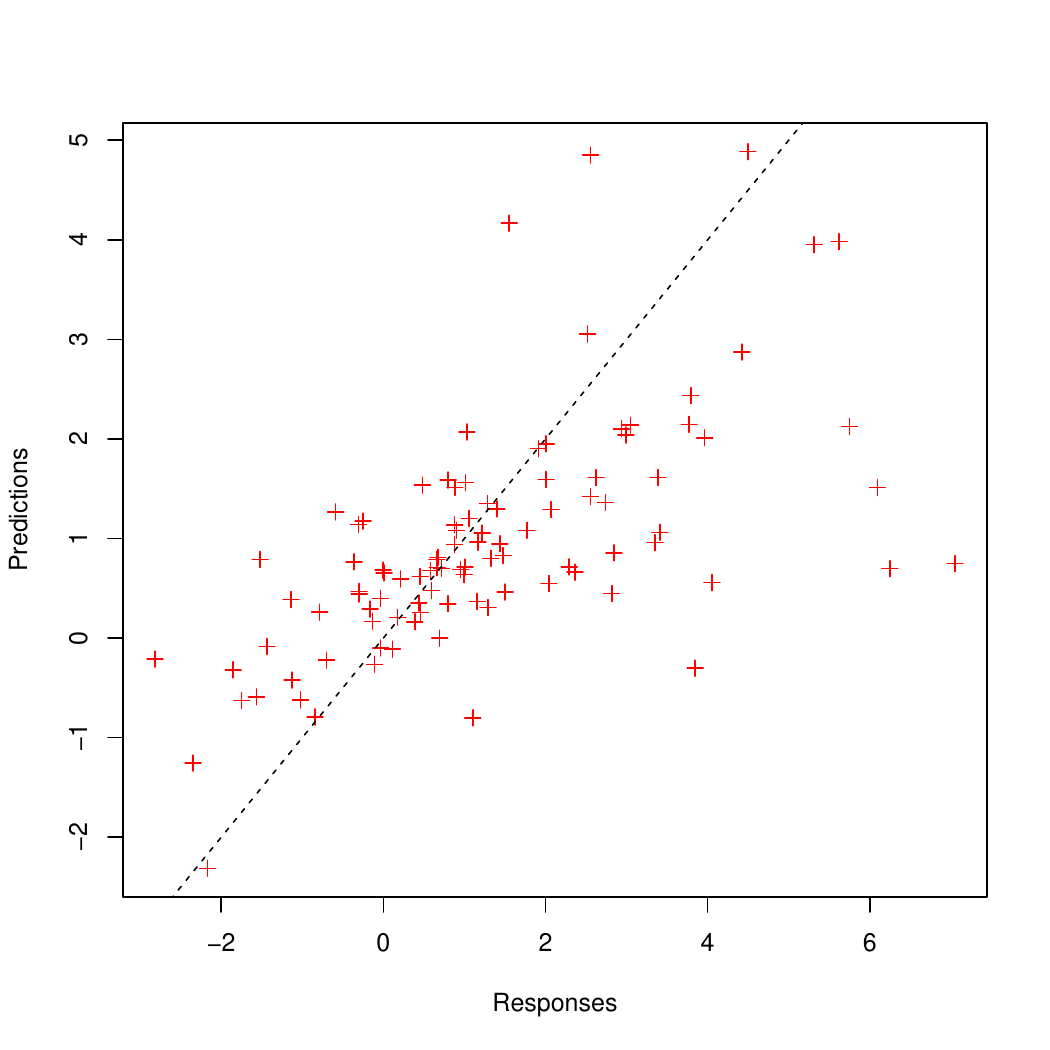}}
  \subfloat[\autoref{m8}.]{\includegraphics[width = .24\textwidth]{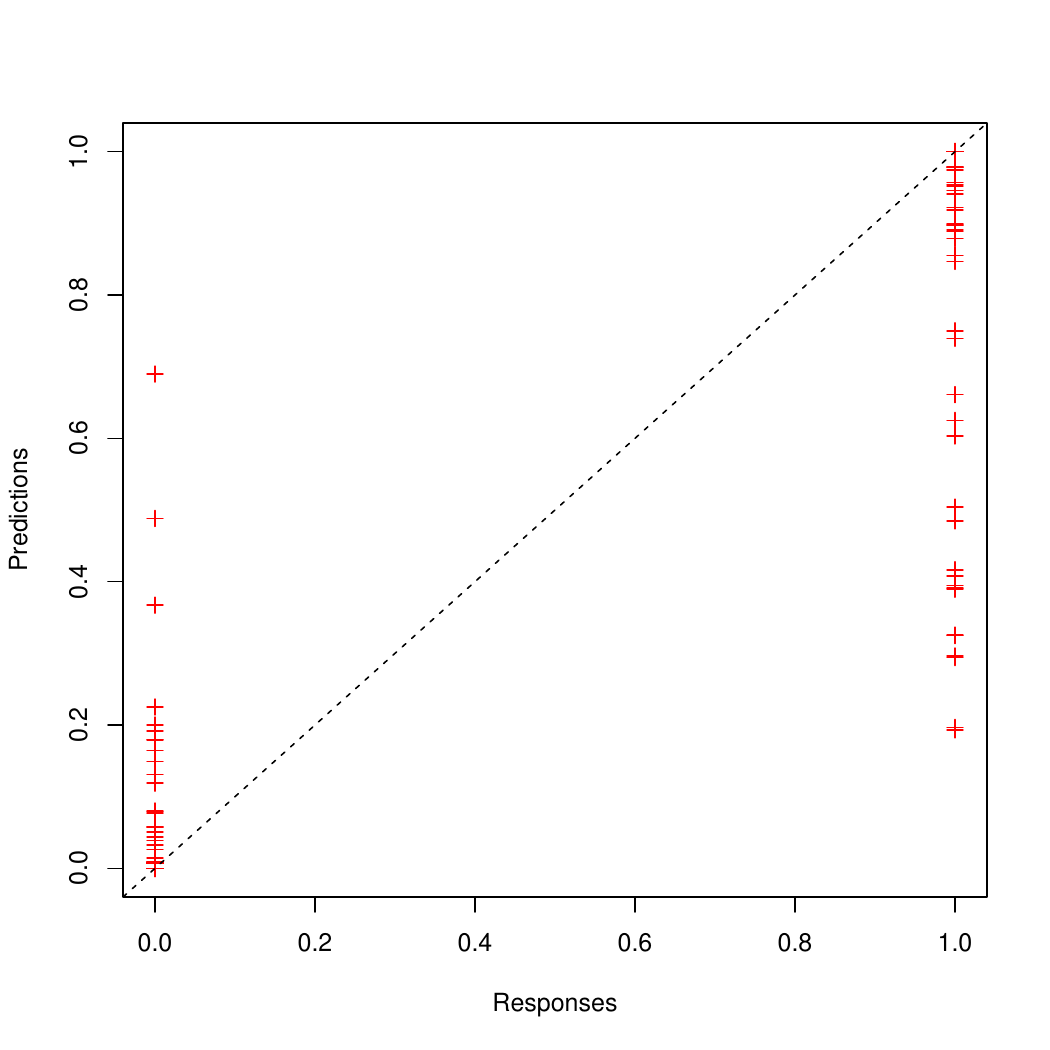}}
\end{figure}

\begin{figure}[t]
  \caption{Examples of reconstruction of the functional dependencies, for covariates
    $1$ to $4$.}
  \label{func}
  \subfloat[\autoref{m1}, uncorrelated design.]{\includegraphics[width = .49\textwidth]{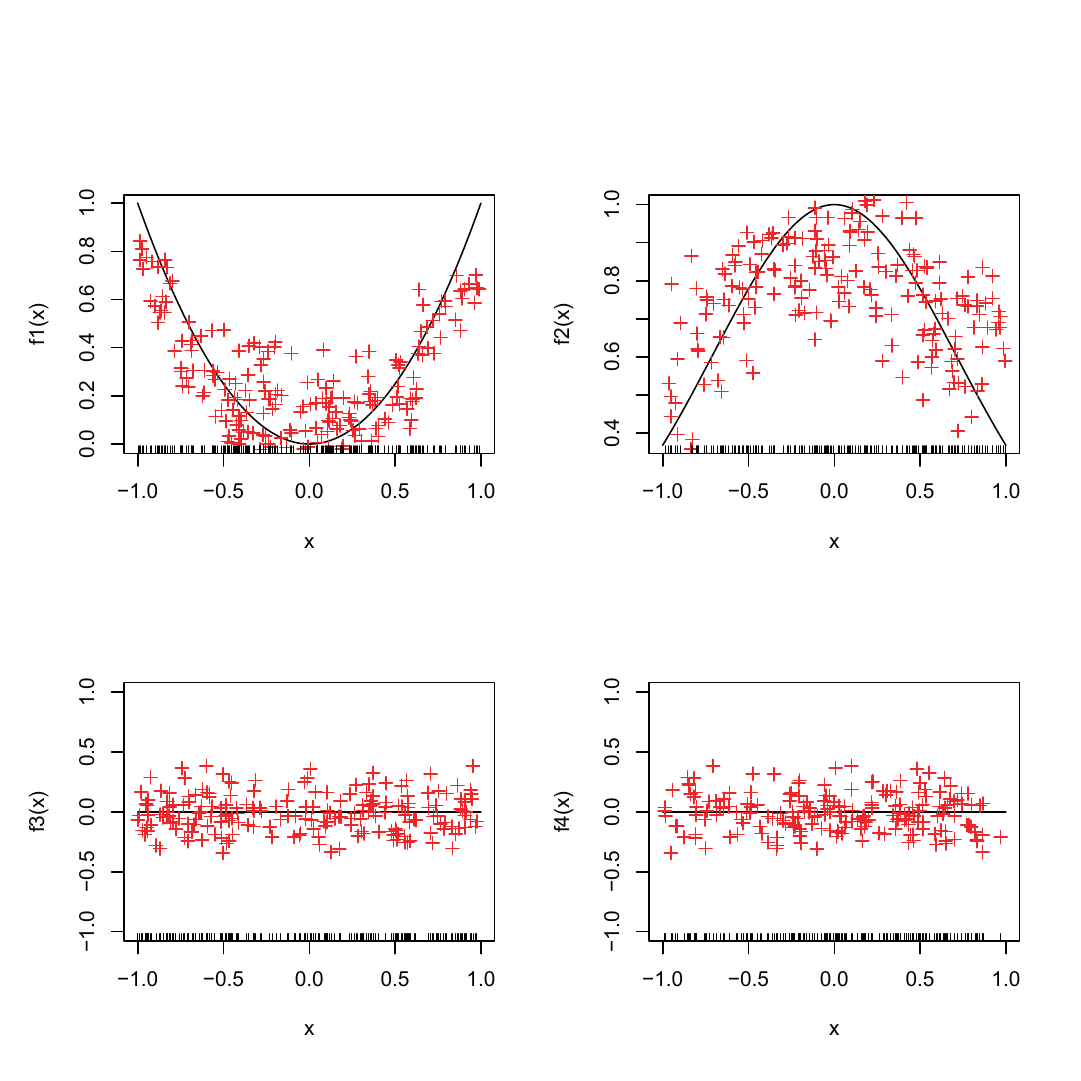}}
  \subfloat[\autoref{m1}, correlated design.]{\includegraphics[width =
    .49\textwidth]{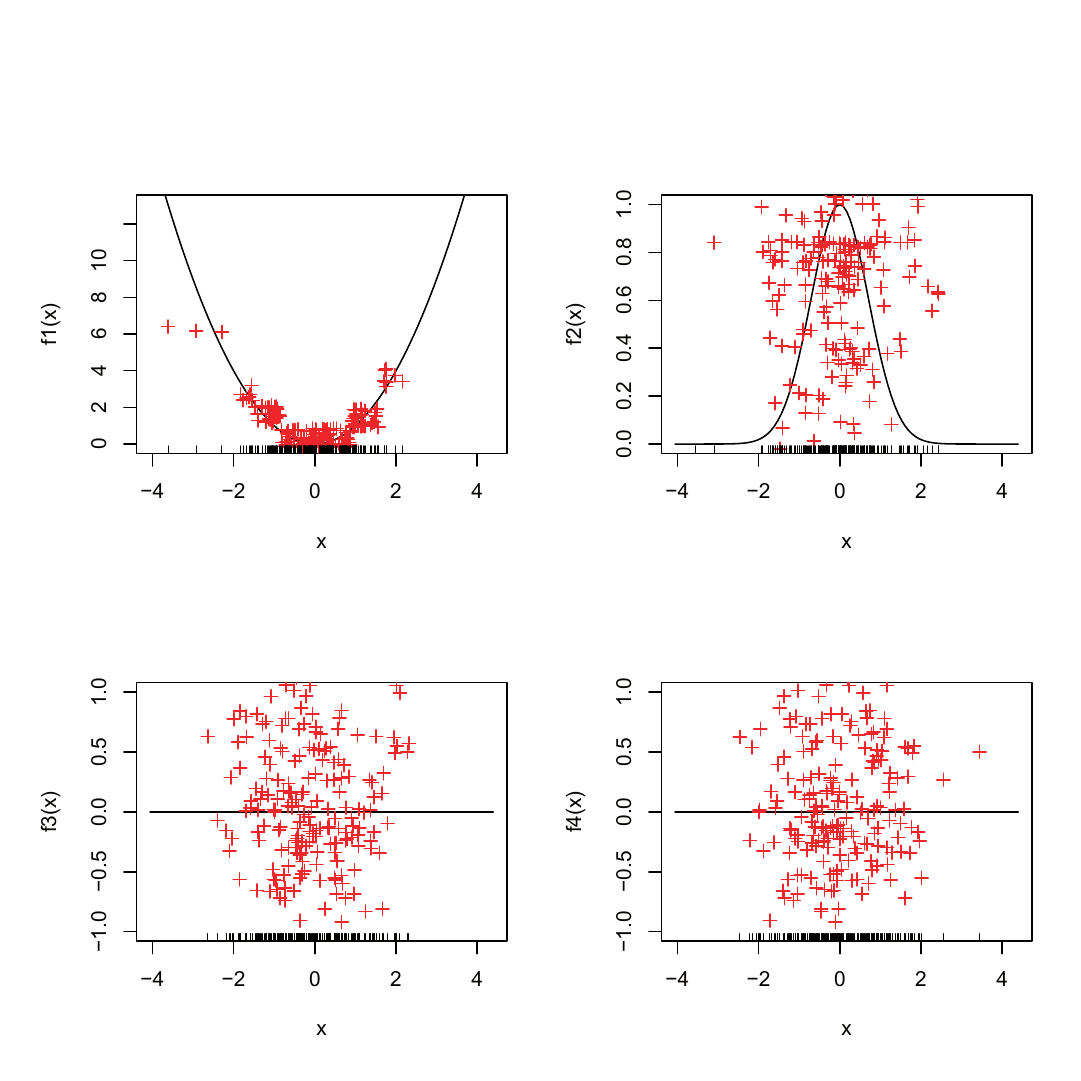}}
  \\
  \subfloat[\autoref{m3}, uncorrelated design.]{\includegraphics[width = .49\textwidth]{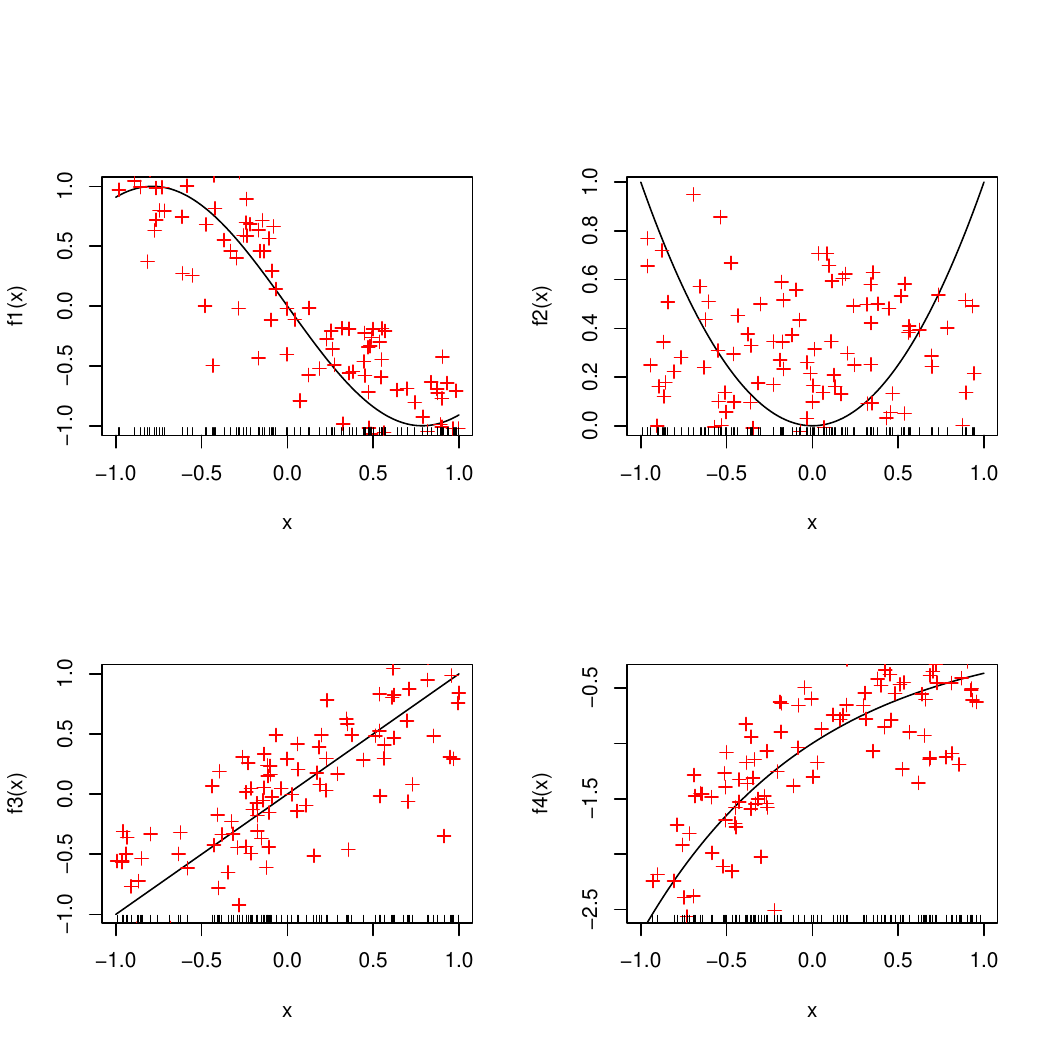}}
  \subfloat[\autoref{m3}, correlated design.]{\includegraphics[width = .49\textwidth]{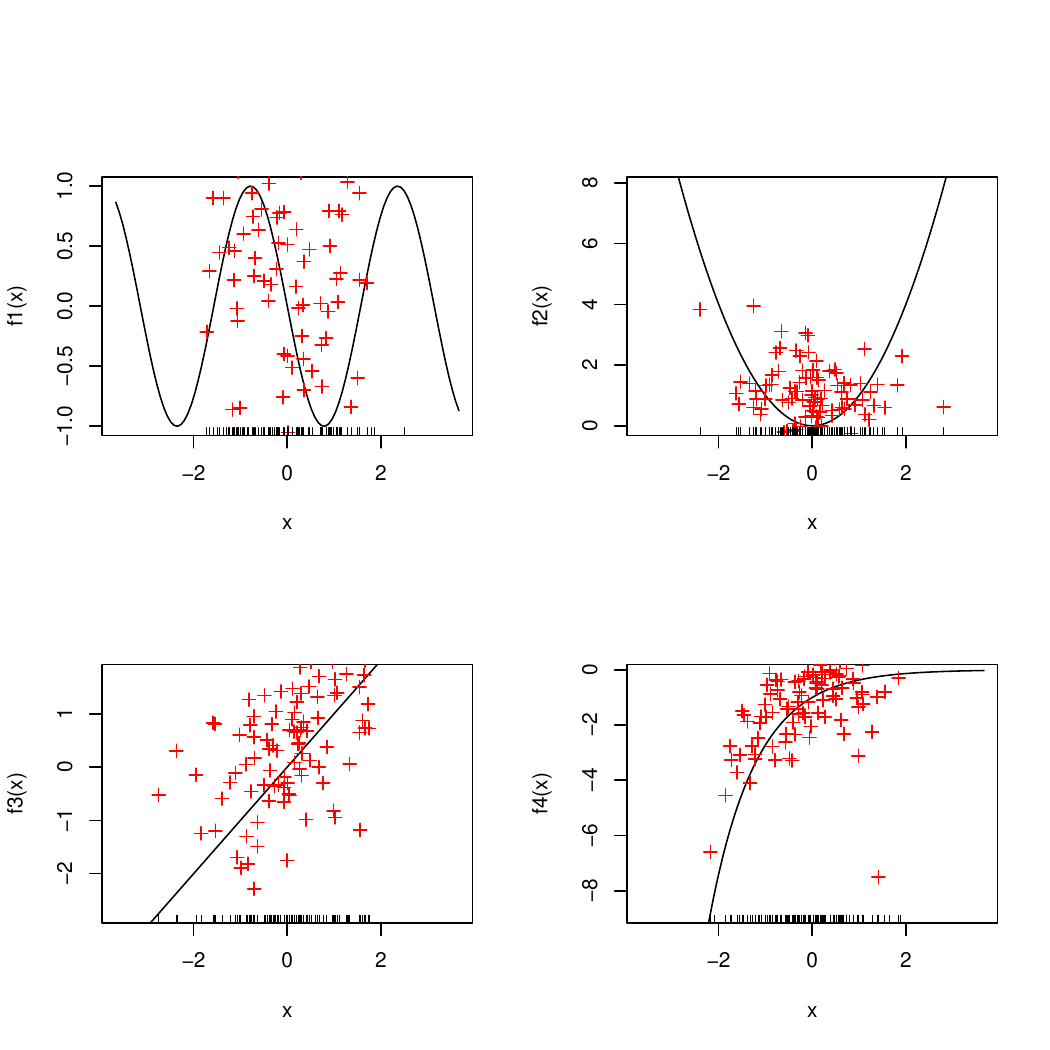}}
  % \\
\end{figure}

\begin{figure}[t]
  \caption{Boxplot of errors, high-dimensional models.}
  \label{box-HD}
  \subfloat[\autoref{mHD1}]{\includegraphics[width=.32\textwidth]{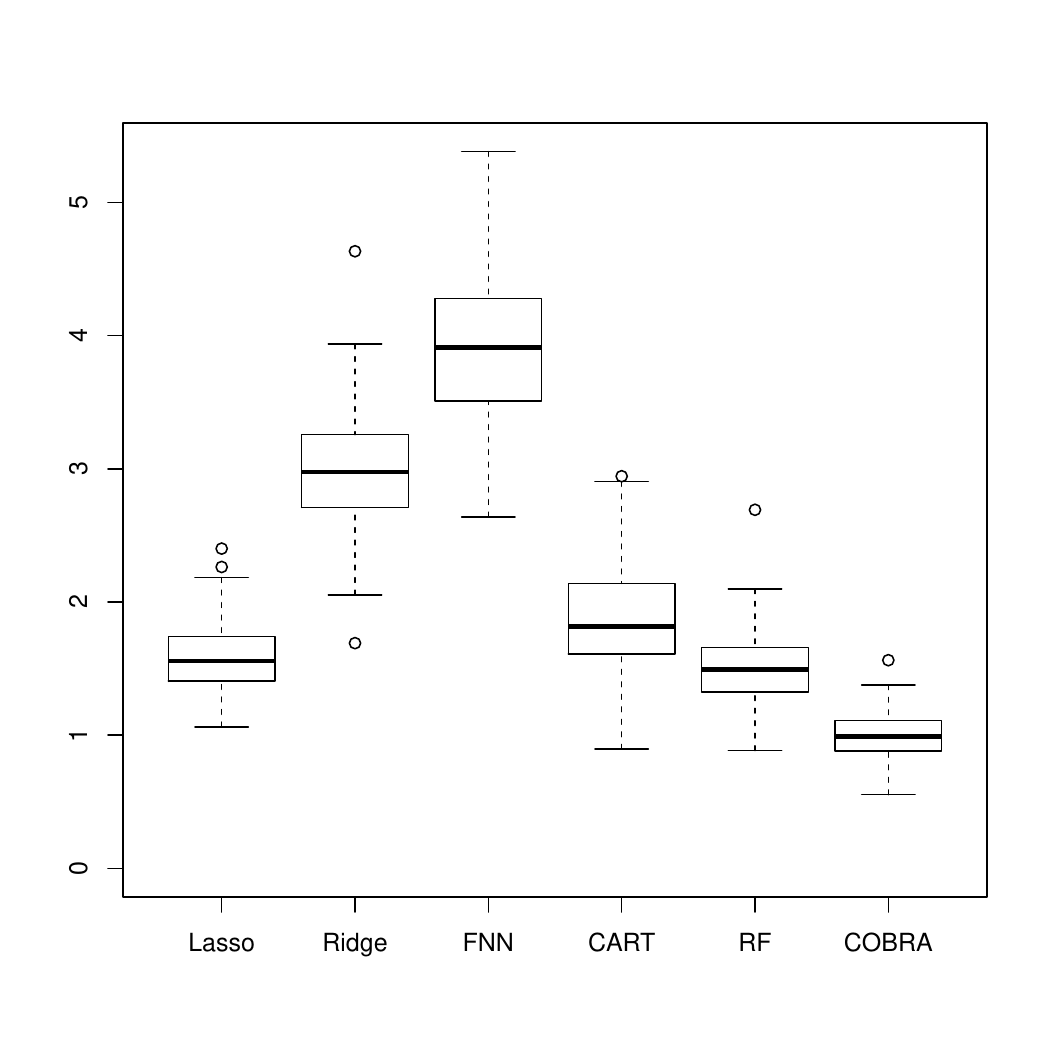}}
  \hfill
  \subfloat[\autoref{mHD2}]{\includegraphics[width=.32\textwidth]{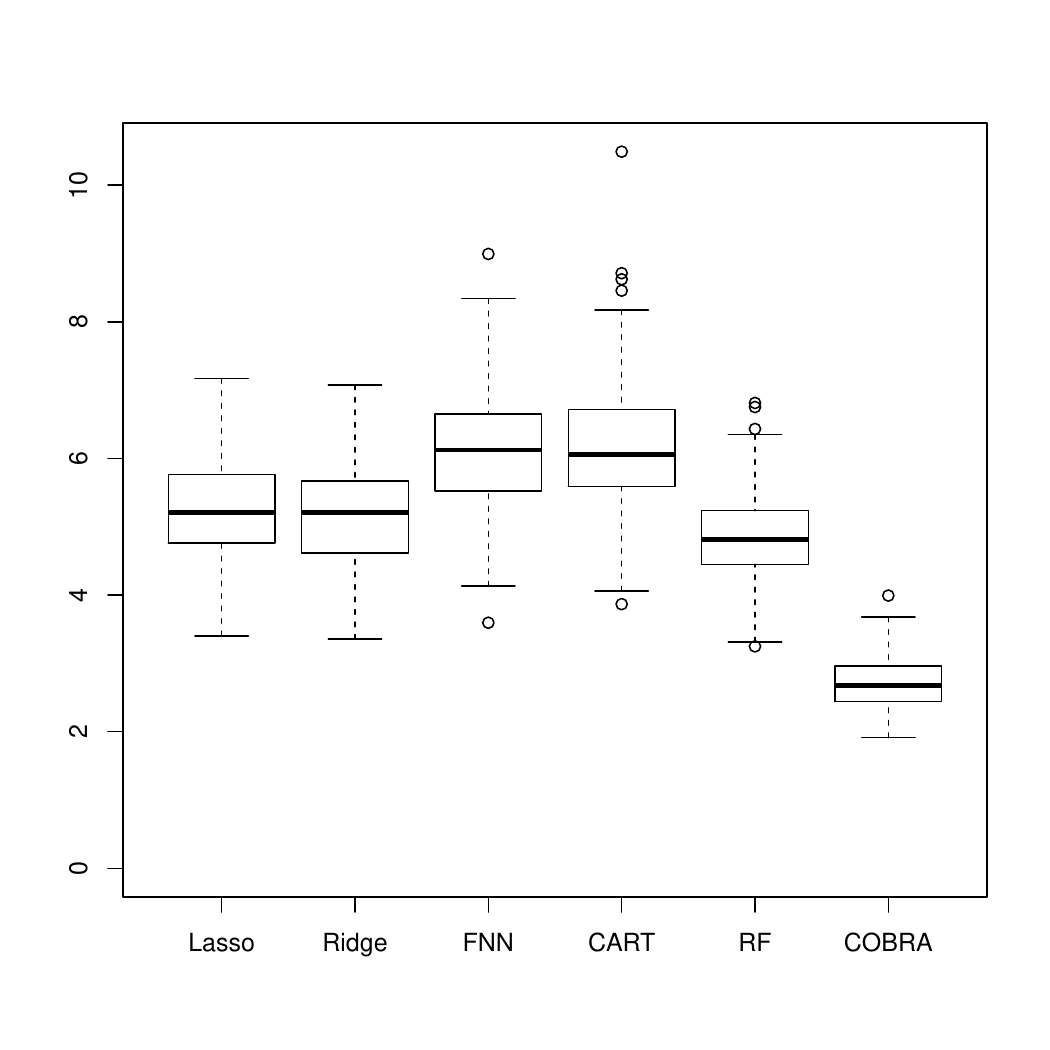}} \hfill
  \subfloat[\autoref{mHD3}]{\includegraphics[width=.32\textwidth]{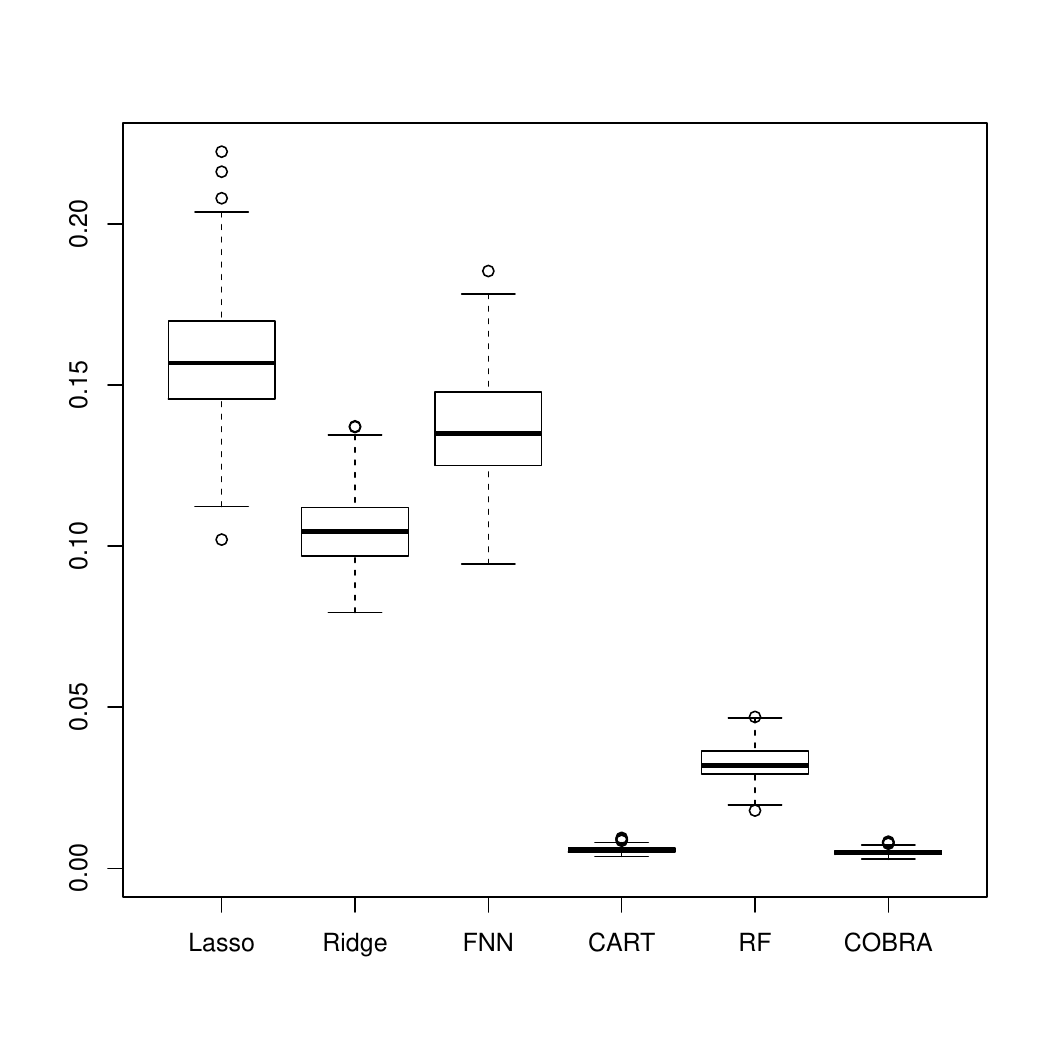}}
\end{figure}

\begin{figure}[h]
  \caption{How stable is \cobra?}
  \label{box-Stab}
  \subfloat[Boxplot of errors: Initial sample is randomly cut ($1000$ replications of \autoref{mSTAB}).]{\includegraphics[width=.49\textwidth]{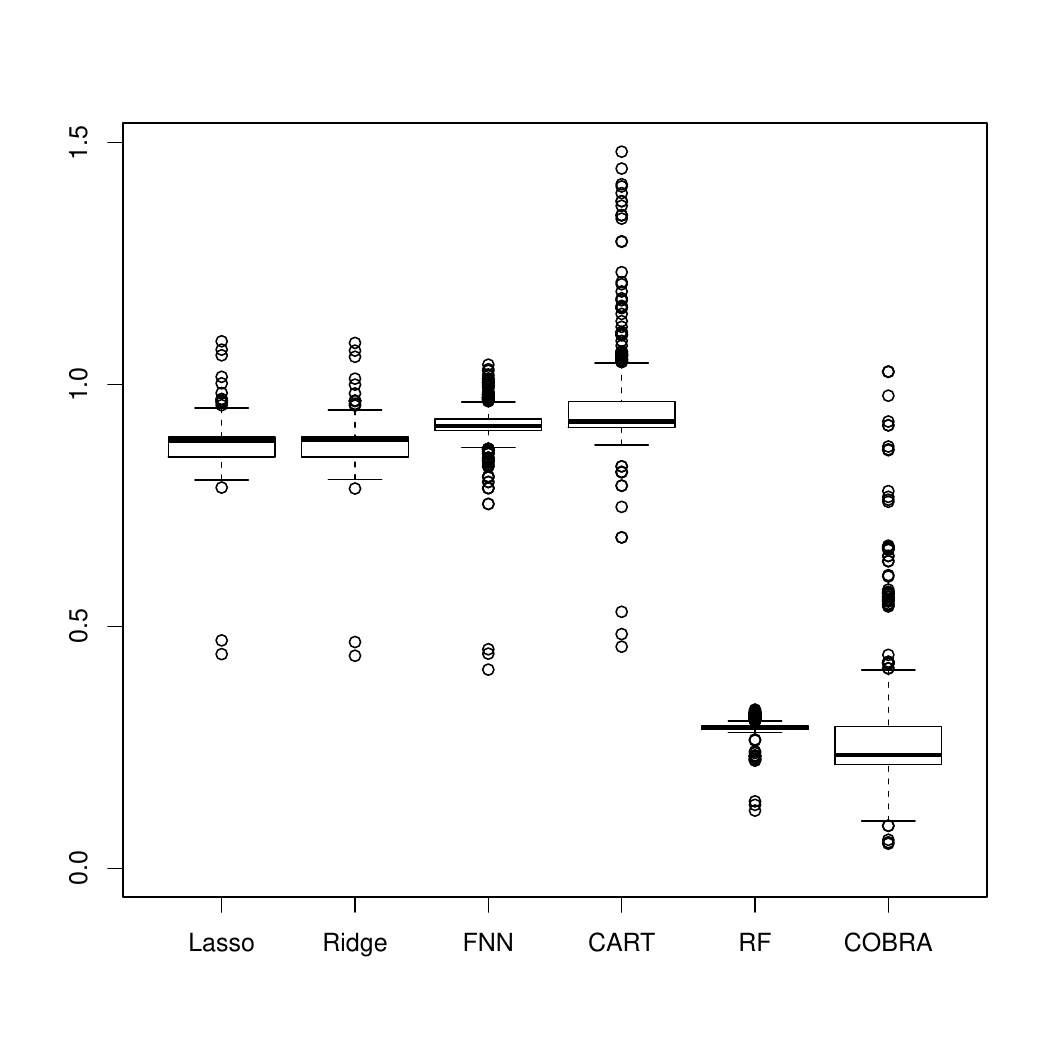}}
  \hfill
  \subfloat[Empirical risk with respect to the size of subsample
    $\mathcal{D}_k$, in \autoref{mSTAB}.]{\includegraphics[width=.49\textwidth]{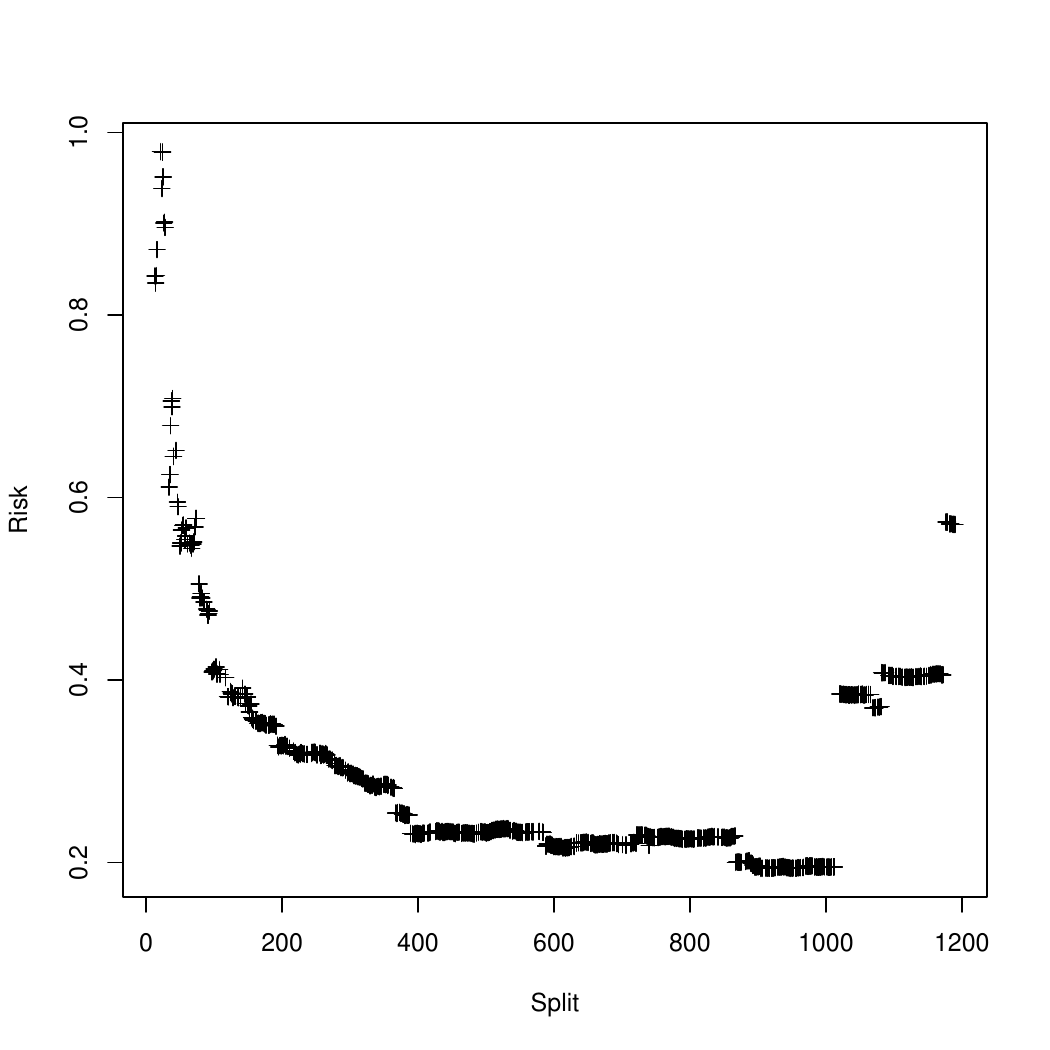}}
\end{figure}

\begin{figure}[h]
  \caption{Boxplot of errors: \texttt{EWA} vs \cobra}
  \label{box-EWA}
  % \begin{center}
  \subfloat[\autoref{mHD1}.]{\includegraphics[width=.24\textwidth]{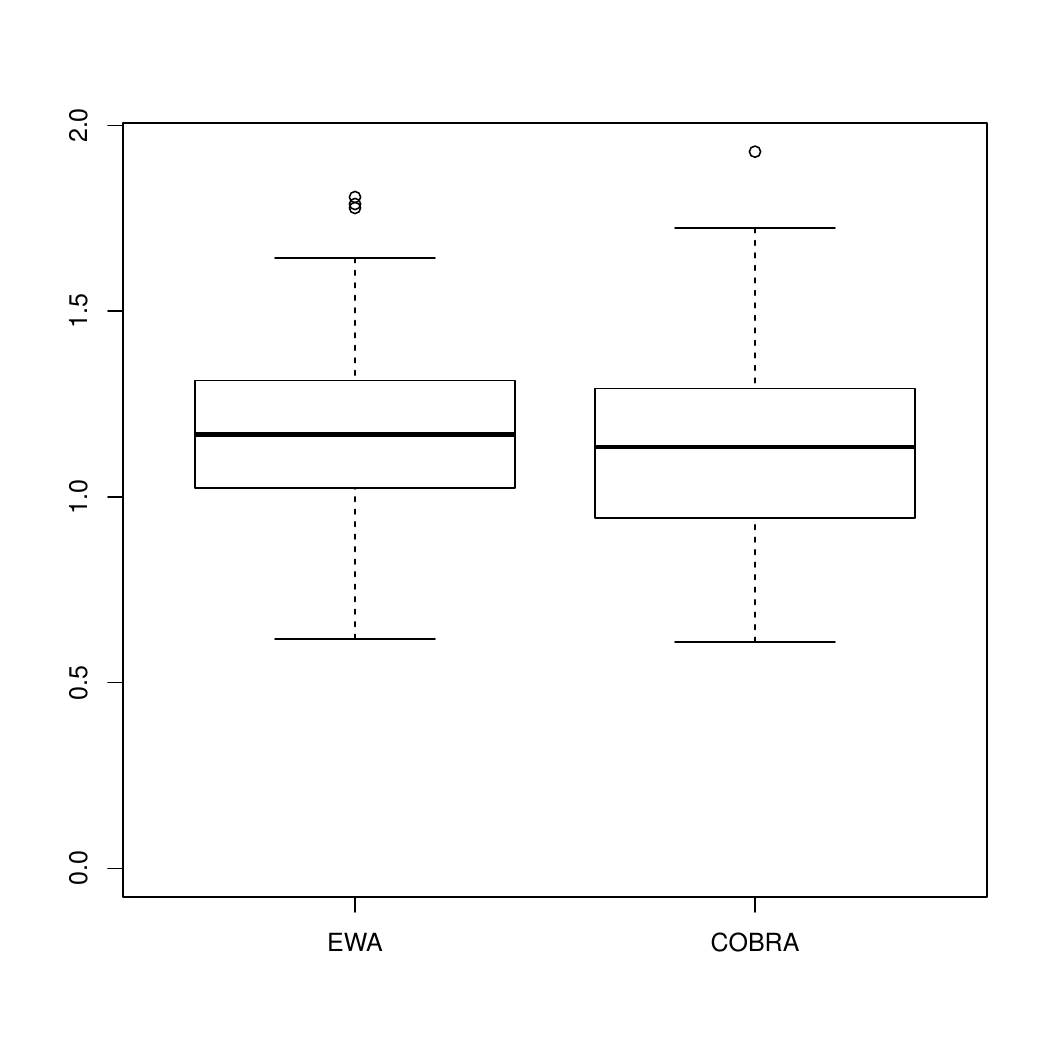}} \hfill
  \subfloat[\autoref{mHD2}.]{\includegraphics[width=.24\textwidth]{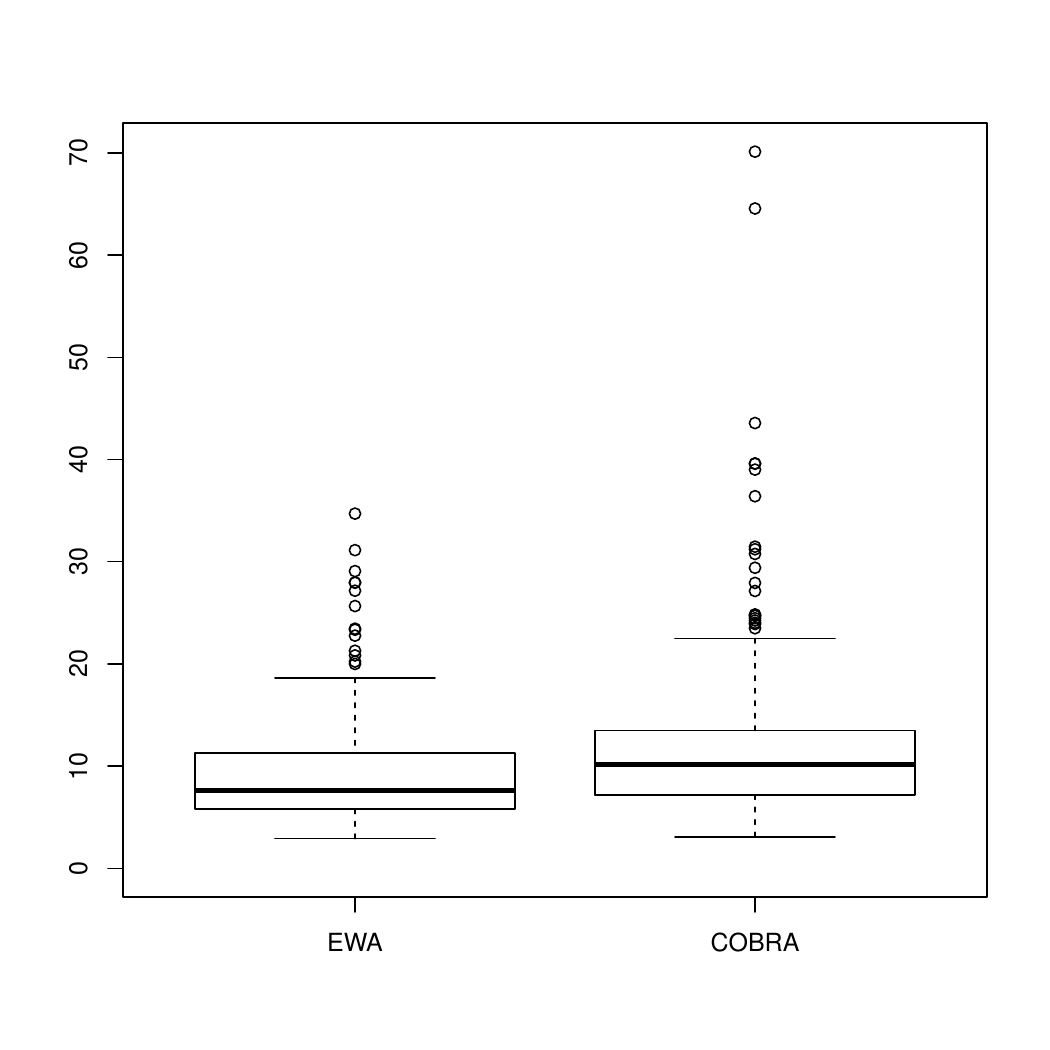}} \hfill
  \subfloat[\autoref{mHD3}.]{\includegraphics[width=.24\textwidth]{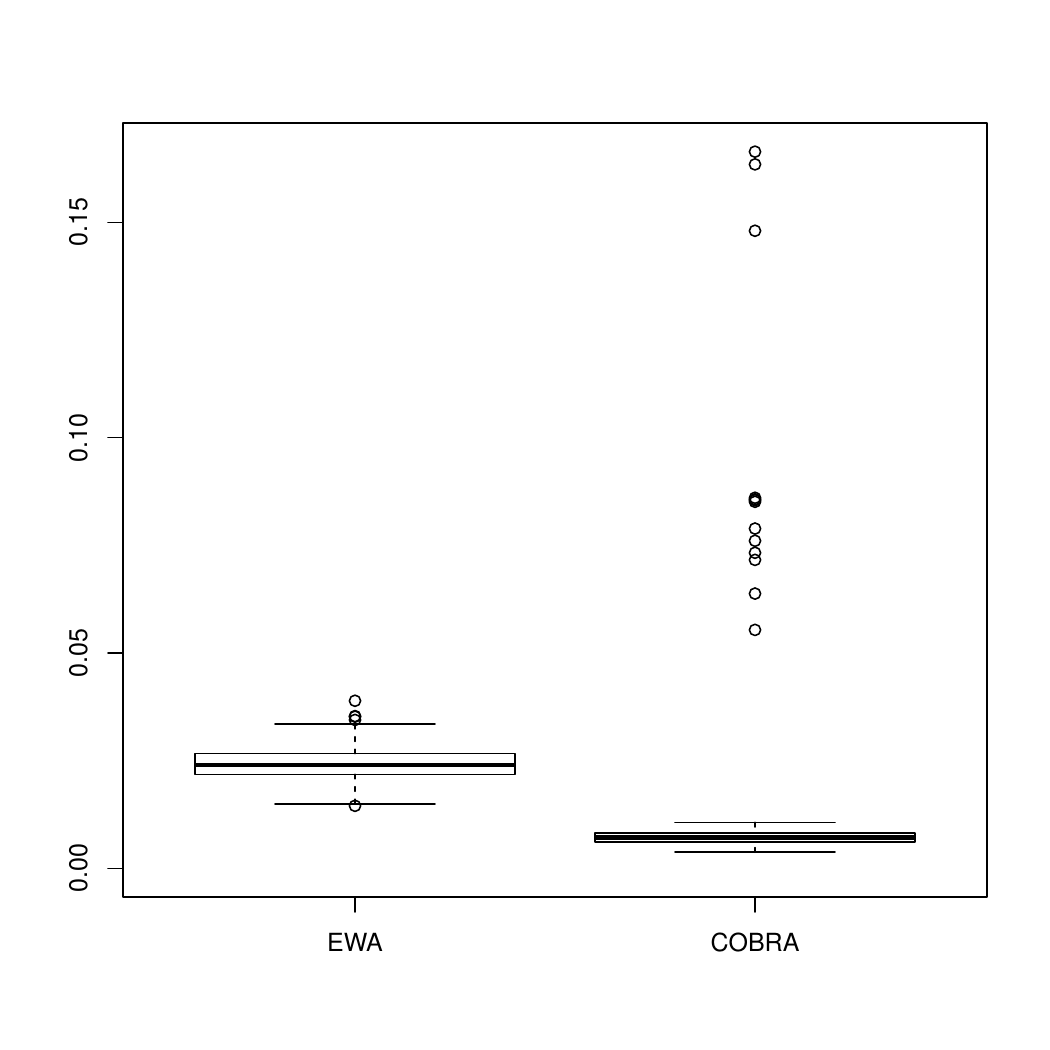}} \hfill
  \subfloat[\autoref{mSTAB}.]{\includegraphics[width=.24\textwidth]{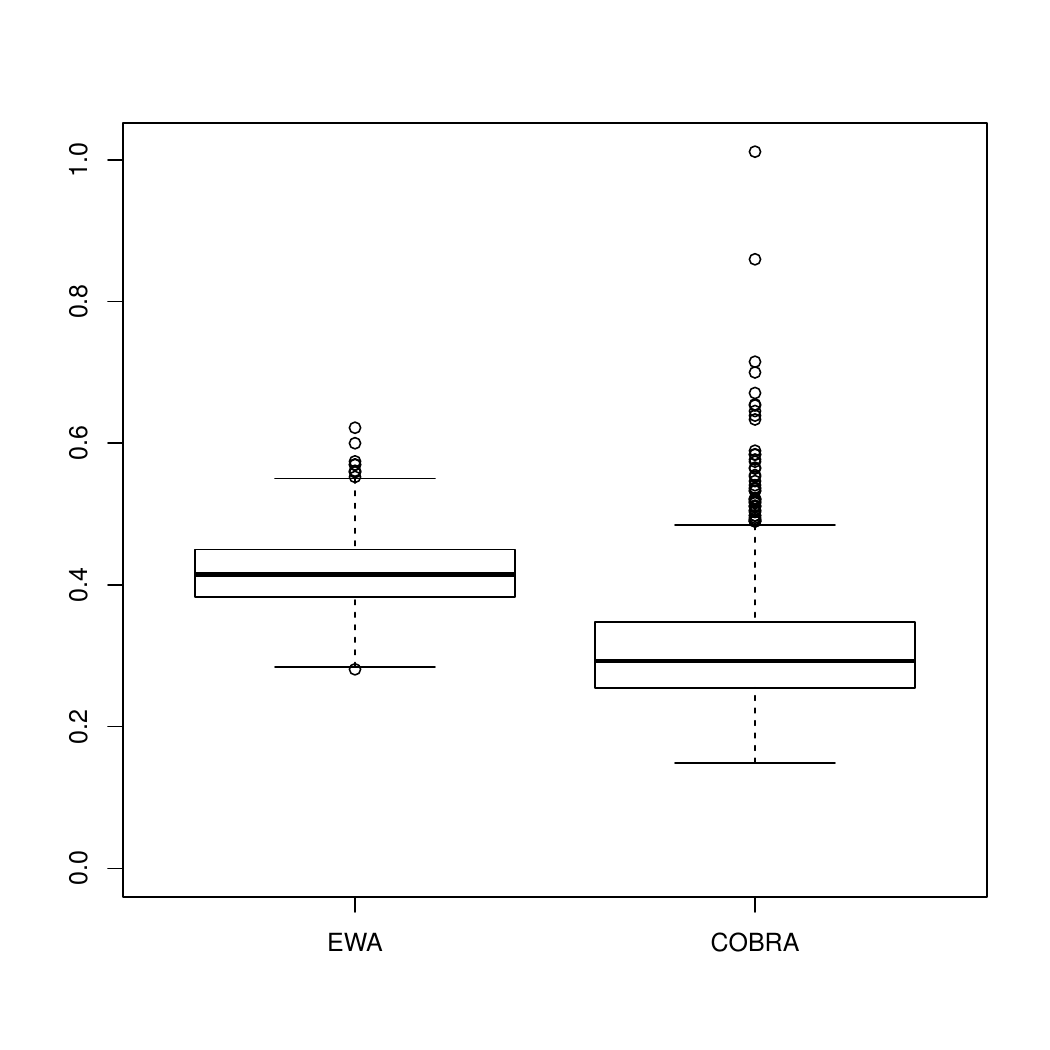}}
  % \end{center}
\end{figure}

\begin{figure}[h]
  \caption{Prediction over the testing set, real-life data sets.}
  \label{real}
  \subfloat[Concrete Slump Test.]{\includegraphics[width = .49\textwidth]{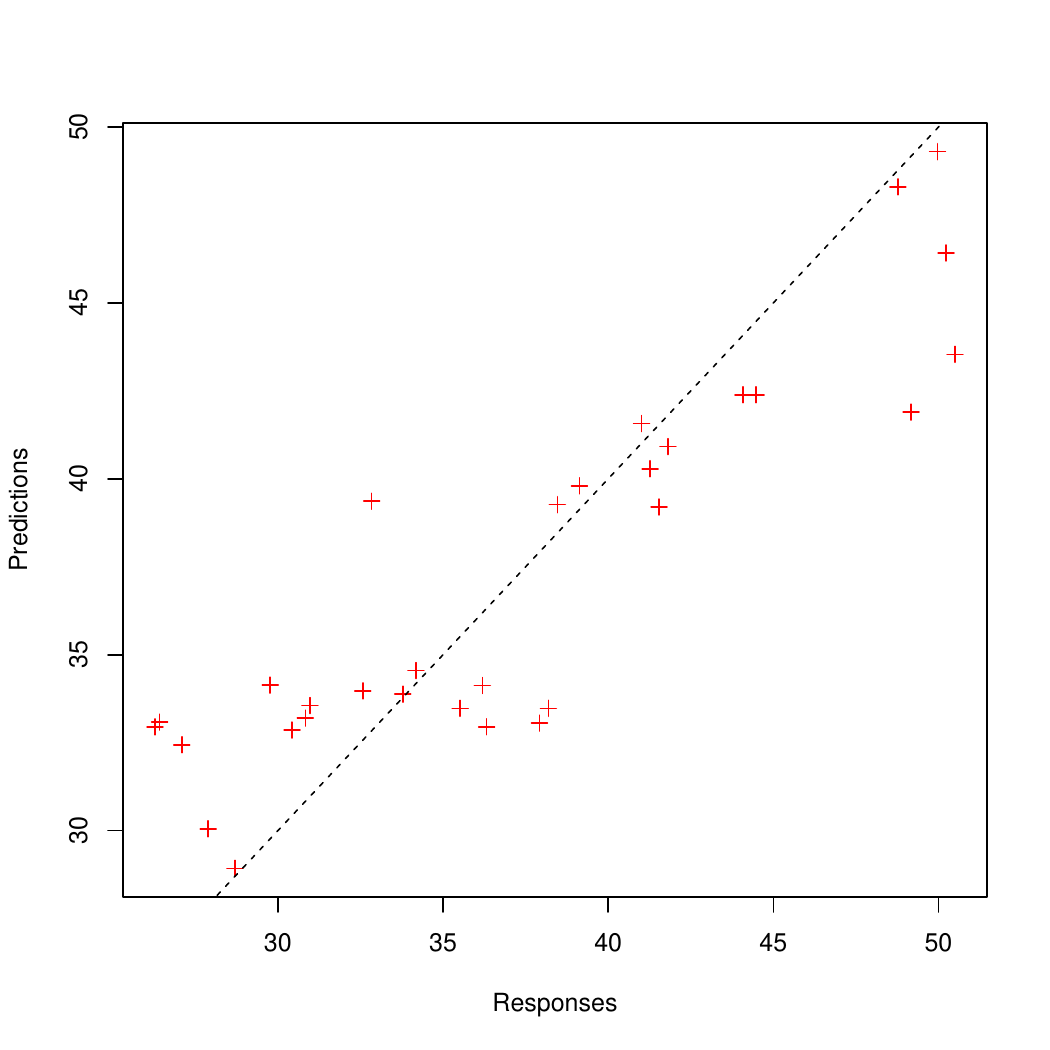}}
  \subfloat[Concrete Compressive Strength.]{\includegraphics[width =
    .49\textwidth]{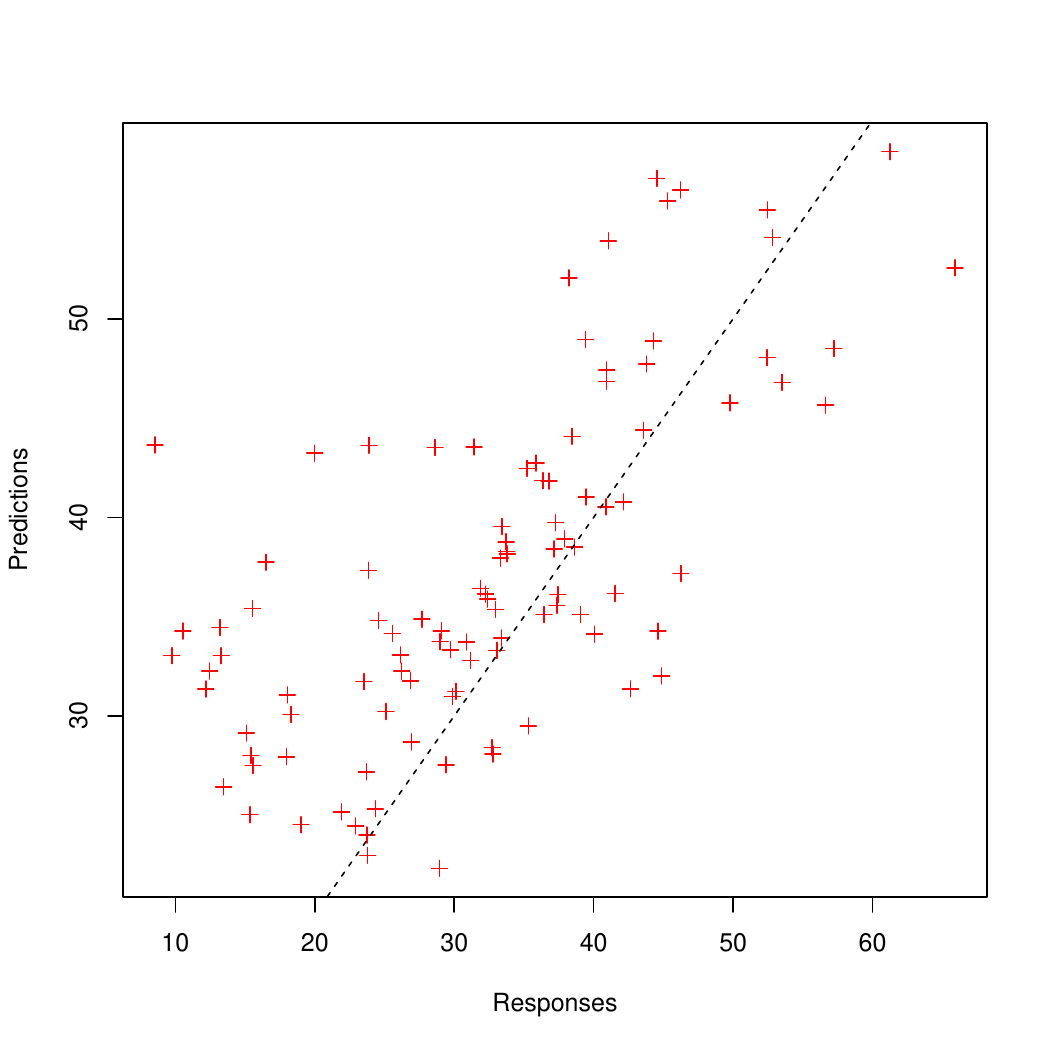}} \\
  \subfloat[Wine Quality, red wine.]{\includegraphics[width =
    .49\textwidth]{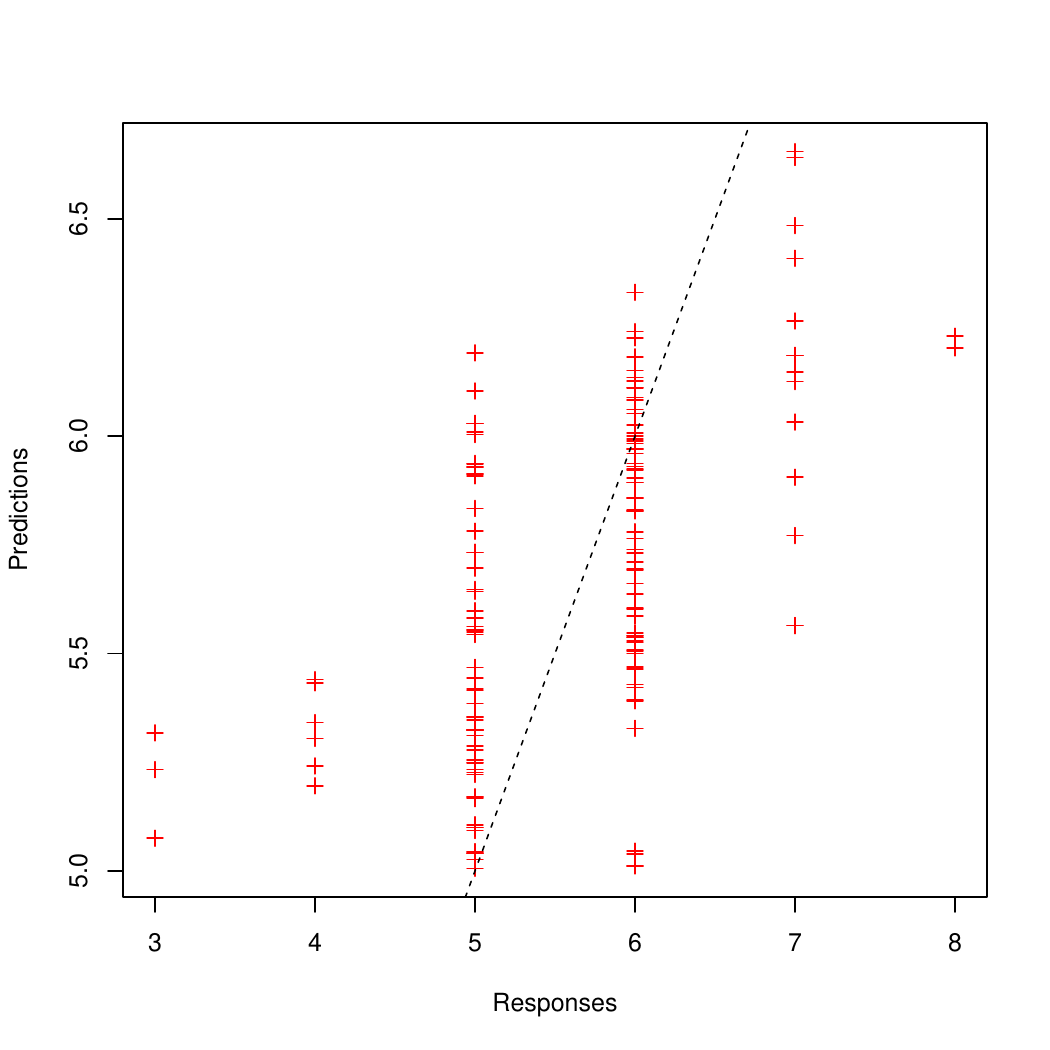}}
  \subfloat[Wine Quality, white wine.]{\includegraphics[width =
    .49\textwidth]{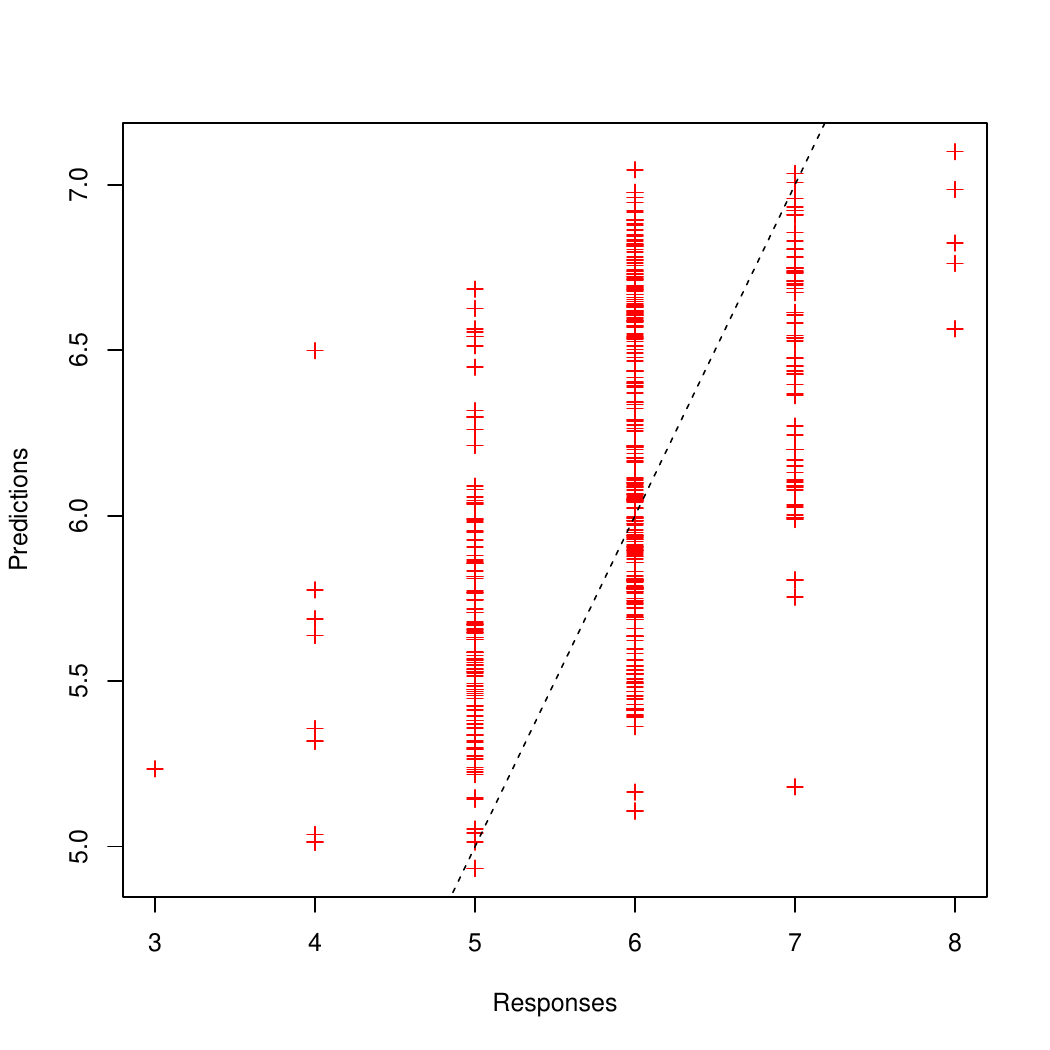}}
\end{figure}
\begin{figure}[h]
  \caption{Boxplot of quadratic errors, real-life data sets.}
  \label{box-real}
  \subfloat[Concrete Slump Test.]{\includegraphics[width = .24\textwidth]{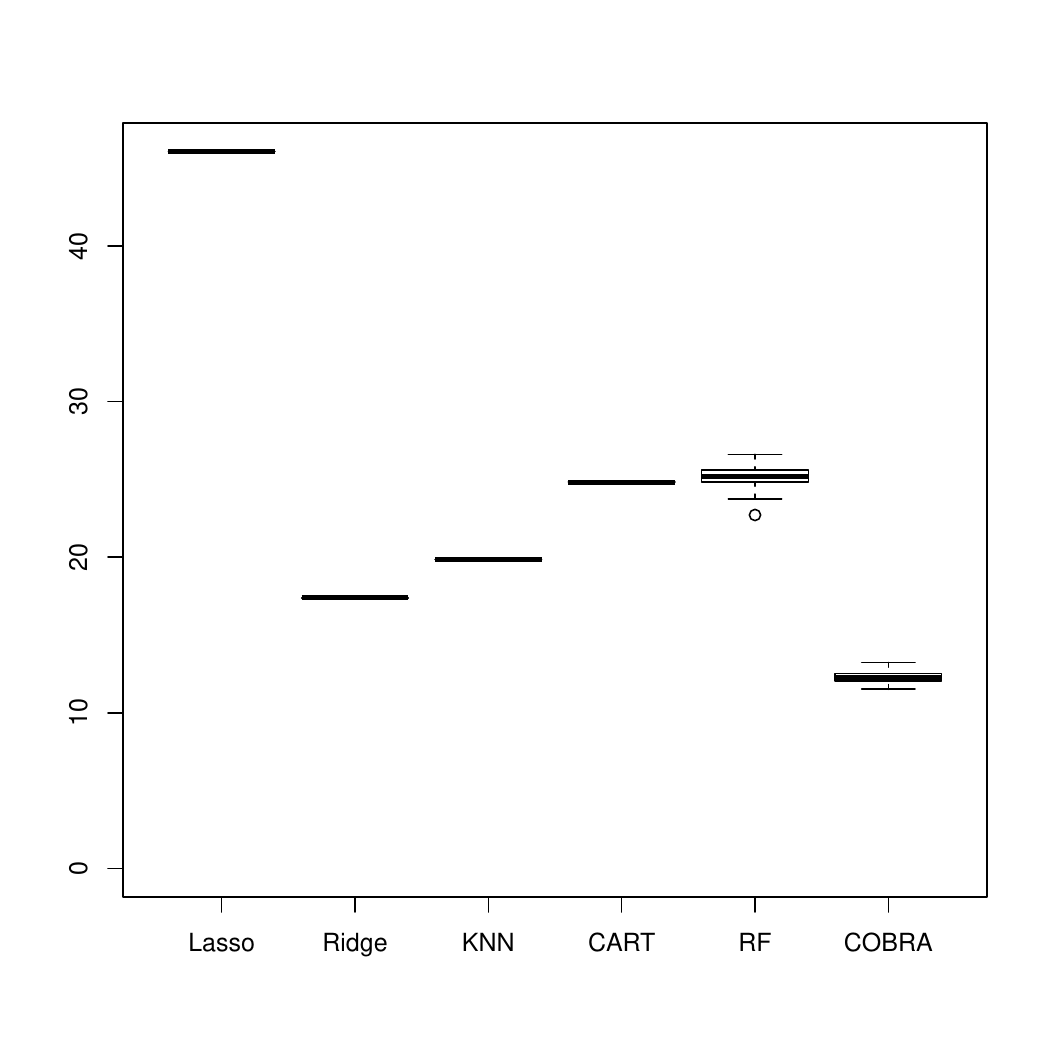}}
  \subfloat[Concrete Compressive Strength.]{\includegraphics[width =
    .24\textwidth]{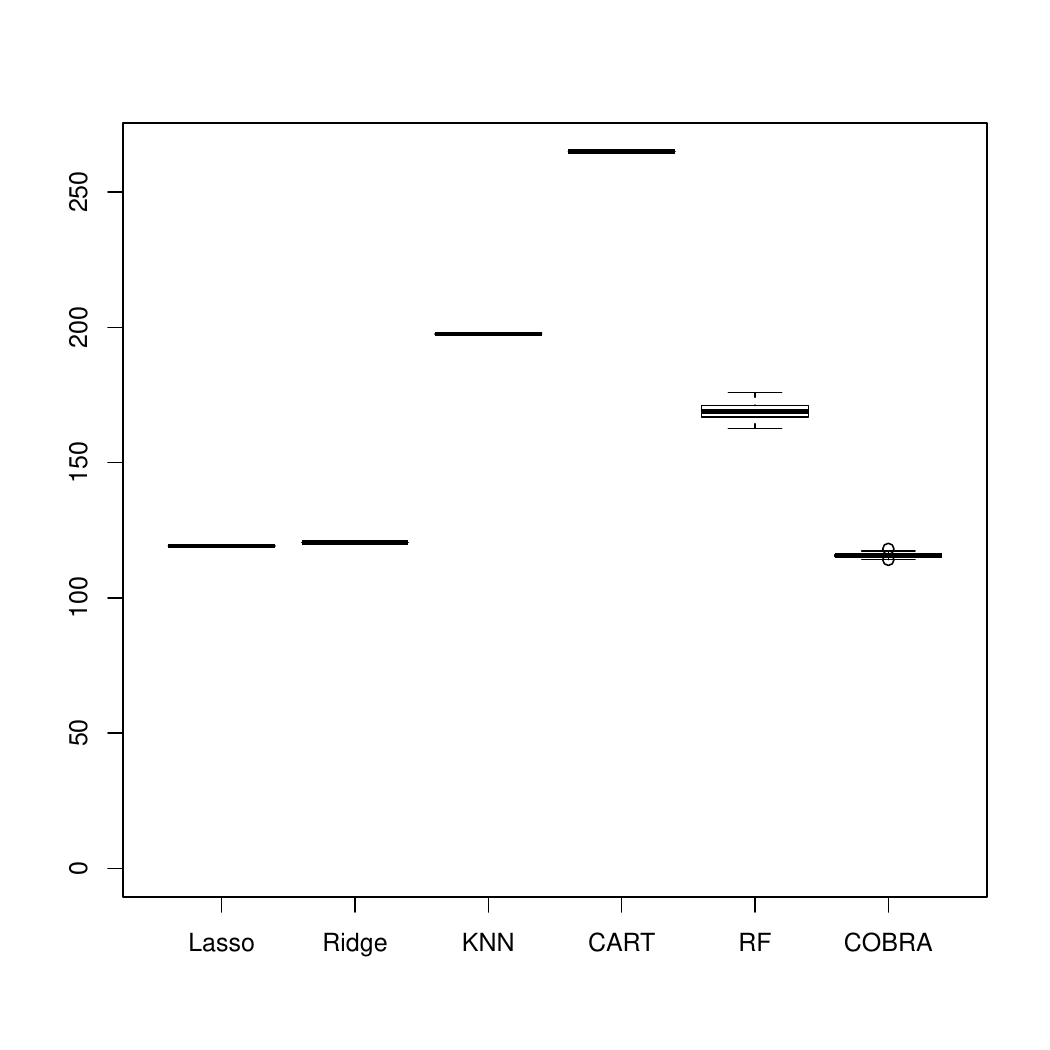}} %\\
  \subfloat[Wine Quality, red wine.]{\includegraphics[width =
    .24\textwidth]{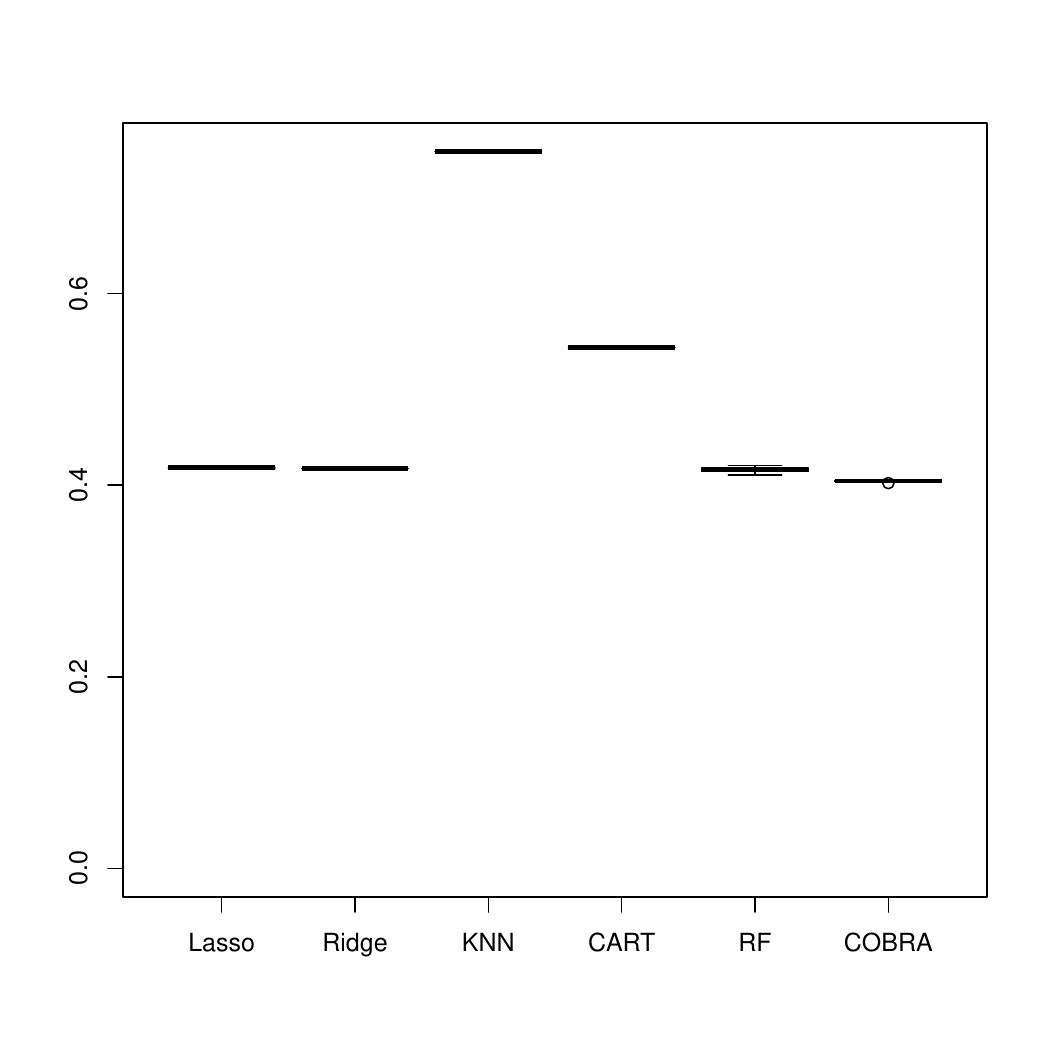}}
  \subfloat[Wine Quality, white wine.]{\includegraphics[width =
    .24\textwidth]{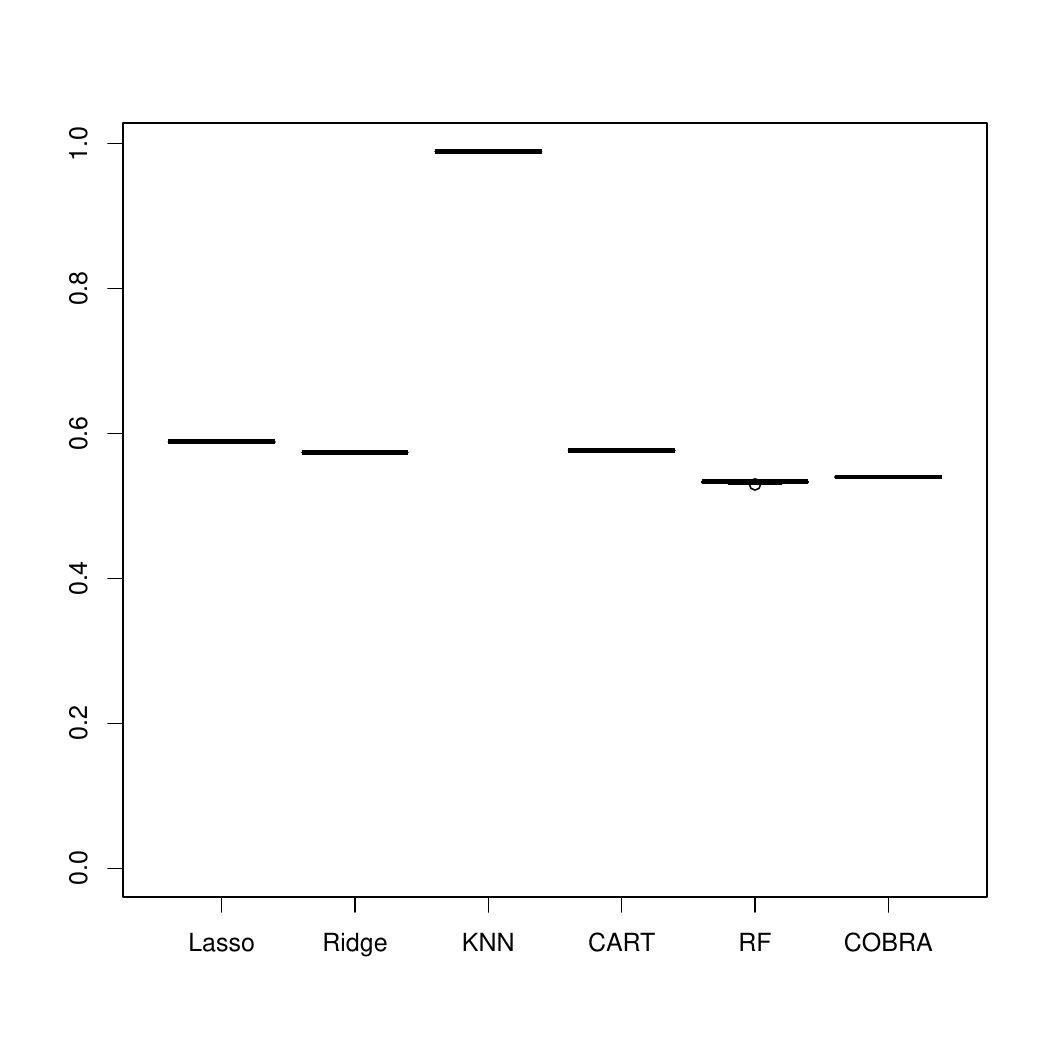}}
\end{figure}

\end{document}